\documentclass{article}

\usepackage{graphicx}
\usepackage{amsmath,amssymb,amsfonts,amsthm}
\usepackage{mathtools}
\usepackage{booktabs}
\usepackage{caption}
\usepackage{algorithmicx}
\newtheorem{thm}{Theorem}[section]
\newtheorem{theorem}[thm]{Theorem}

\newtheorem{lemma}[thm]{Lemma}

\newtheorem{assumption}[thm]{Assumption}
\newtheorem{proposition}[thm]{Proposition}

\newtheorem{definition}[thm]{Definition}

\newtheorem{remark}[thm]{Remark}

\newcommand{\beq}{\begin{equation}}
\newcommand{\eeq}{\end{equation}}
\newcommand{\beqa}{\begin{eqnarray}}
\newcommand{\eeqa}{\end{eqnarray}}
\newcommand{\beqas}{\begin{eqnarray*}}
\newcommand{\eeqas}{\end{eqnarray*}}
\newcommand{\ei}{\end{itemize}}

\newcommand{\vgap}{\vspace{.1in}}

\def\defi{\vcentcolon=}

\newcommand{\R}{\mathbb{R}}

\newcommand{\cA}{{\cal A}}

\newcommand{\A}{{\cal A}}

\newcommand{\inner}[2]{\langle #1,#2\rangle}

\newcommand{\tr}{\mathrm{tr}}

\newcommand{\ignore}[1]{}

\newtheorem{innercustomgeneric}{\customgenericname}
\providecommand{\customgenericname}{}
\newcommand{\newcustomtheorem}[2]{%
  \newenvironment{#1}[1]
  {%
   \renewcommand\customgenericname{#2}%
   \renewcommand\theinnercustomgeneric{##1}%
   \innercustomgeneric
  }
  {\endinnercustomgeneric}
}
\newcustomtheorem{blackbox}{Blackbox}

\usepackage{setspace}
\usepackage{hyperref}
\usepackage{bm}
\usepackage{algorithm}
\usepackage{array}
\usepackage{placeins}
\usepackage{algpseudocode}
\usepackage{xspace}
\usepackage{enumitem}
\makeatletter
\renewcommand{\ALG@beginalgorithmic}{\footnotesize}
\makeatother
\usepackage{xcolor}
\usepackage{natbib}
\bibpunct[, ]{[}{]}{,}{n}{}{,}%
\usepackage{tikz}
\usetikzlibrary{calc,positioning,arrows.meta,quotes,decorations.pathreplacing}
\usepackage{pgfplots}
\pgfplotsset{compat=1.17}

\providecommand{\ourmethodplus}{HALLaR+\xspace}
\providecommand{\cuourmethodplus}{cuHALLaR+\xspace}

\newcommand{\cD}{\mathcal{D}}

\newcommand{\cP}{\mathcal{P}}

\newcommand{\Lip}[1]{({\rm Lip})_{#1}}

\newcommand{\tD}[1]{\widetilde{\mathcal D}_{#1}}

\newcommand{\tLip}[1]{\widetilde{(\operatorname{Lip}})_{#1}}

\newcommand{\tM}[1]{\widetilde{M}_{#1}}

\DeclareMathOperator*{\argmin}{arg\,min}

\allowdisplaybreaks

\title{
A Low-Rank Augmented Lagrangian Frank-Wolfe-Based Method \\ for SDPs with Unbounded Feasible Regions 
\thanks{\textbf{Funding}:
Jacob M. Aguirre is supported by the National Science Foundation Graduate Research Fellowship under Grant No. DGE-2039655.
Renato D.C. Monteiro is supported by AFOSR Grant FA9550-25-1-0131.
}}
\author{
Jacob M. Aguirre\thanks{H. M. Stewart School of Industrial and Systems Engineering, Georgia Tech, Atlanta, GA, 30332-0205.
\protect\href{mailto:aguirre@gatech.edu}{aguirre@gatech.edu}}
\and
Renato D.C. Monteiro\thanks{School of Industrial and Systems Engineering, Georgia Institute of Technology, Atlanta, GA, 30332-0205.
\protect\href{mailto:monteiro@isye.gatech.edu}{monteiro@isye.gatech.edu}}
\and
Arnesh Sujanani\thanks{Department of Computer Science and Engineering, Chinese University of Hong Kong, Hong Kong.
\protect\href{mailto:arneshmurlisujanani@cuhk.edu.hk}{arneshmurlisujanani@cuhk.edu.hk}\\}
}
\date{\today}

\begin{document}
\setcitestyle{square}
\maketitle

\begin{abstract}
In contrast to HALLaR (Monteiro et al. 2026), which is designed for solving problems with a trace-bounding constraint, this work develops a new augmented Lagrangian (AL) method, HALLaR+, for solving large-scale semidefinite programming problems without such constraint.
Even though a trace constraint which preserves an optimal solution can always be
added, choosing its bound is challenging: a value that is too small may
exclude all optimal solutions, whereas one that is too large may
substantially increase the computational effort.
HALLaR+ is also a low-rank inexact AL method but,
in contrast to HALLaR, it doubles the penalty parameter at every AL iteration and its AL subproblems have unbounded feasible regions that
are approximately solved by a relaxed unbounded Frank--Wolfe (RUFW) variant.
Its analysis relies on a uniform bound on the diameters of the sublevel sets
of all generated AL subproblems instead of a bounded feasible set.
Assuming that a primal--dual KKT pair exists and the primal optimal set
is nonempty and bounded, it is shown that HALLaR+ returns a $\varepsilon$-primal-dual feasible solution in  $\tilde {\mathcal O}(\varepsilon^{-2})$
total RUFW iterations.
Numerical results suggest that HALLaR+ has very promising competitive performance against current
state-of-the-art solvers.

\vgap


\end{abstract}

\section{Introduction}

This paper develops 
a hybrid
convex--nonconvex augmented Lagrangian (AL) variant of HALLaR of 
\citep{monteiro2024low}, namely HALLaR+, which can solve  SDPs, in general standard form, i.e.,
\begin{align}
\min_{X\in\mathbb{S}^{n}}\left\{\inner{C}{X}:\mathcal A X=b,\ X\succeq0\right\}
 \label{eq:intro_primal_sdp}
\end{align}
where $\mathbb R^m$ is the $m$-dimensional Euclidean space equipped with the standard inner product, $\mathbb S^n$ is the space of $n\times n$ real symmetric matrices equipped
with the Frobenius inner product, $\succeq0$ denotes the Löwner partial order on $\mathbb{S}^{n}$, so that $X\succeq0$ means $X$ is positive semidefinite (PSD), $C\in\mathbb S^n$, $b\in\mathbb R^m$, and
$\mathcal A:\mathbb S^n\to\mathbb R^m$ is a linear operator with adjoint $\mathcal A^*:\mathbb R^m\to\mathbb S^n$.
In contrast, HALLaR can only solve SDPs where one of the constraints in \eqref{eq:intro_primal_sdp} is a bounded trace constraint, i.e.,
$\tr(X) \le \tau$ for some $\tau >0$.
Even though a trace constraint can always be added to \eqref{eq:intro_primal_sdp} with a constant $\tau$ that preserves at least one optimal solution, choosing such  a $\tau$ in practice is challenging: a value that is too small may exclude all optimal solutions, whereas one that is too large may substantially increase the computational effort. 

\vspace{3mm}

\noindent\textbf{Approach and relation to HALLaR.}
\ourmethodplus is an inexact AL method that generates sequences
$\{Z_t\}$ and $\{p_t\}$ according to the recursion
\begin{subequations}\label{eq:intro_al_updates}
\begin{align}
    Z_{t} \quad &\approx \quad
    \argmin_{X\in \mathbb{S}^{n}}\quad
    \left\{\mathcal{L}_{\beta_{t-1}}(X;p_{t-1})\defi
     \inner{C}{X} + \inner{p_{t-1}}{b-\mathcal{A}X}
     + \frac{\beta_{t-1}}{2} \|b-\mathcal{A}X\|^2\,: \ X\succeq0\right\},
    \label{eq:intro_al_primal_update}
    \\
    p_{t} \quad &= \quad
    p_{t-1} + \beta_{t-1}(b-\mathcal{A}Z_{t})\,,\qquad\beta_{t}=\gamma\beta_{t-1}\,,
    \label{eq:intro_al_multiplier_update}
\end{align}
\end{subequations}
with $\gamma=2$. HALLaR differs from the above updates \eqref{eq:intro_al_primal_update} and \eqref{eq:intro_al_multiplier_update} in two aspects. First, 
as HALLAR solves SDPs over the set
$\{X\in \mathbb S^{n}  : A(X)=b, \ \tr(X)\leq  \tau, \ X \succeq 0\}$,
its AL subproblems are obtained by dualizing $A(X)=b$, while keeping $\ \tr(X) \leq \tau$ and $X \succeq 0$ as constraints.
Hence, in contrast to \eqref{eq:intro_al_primal_update}, their feasible regions are bounded. Second, the theoretical analysis of HALLaR assumes $\gamma=1$, and hence that the penalty parameters are kept constant.  

Both HALLaR and \ourmethodplus use Frank--Wolfe (FW), or conditional-gradient, variants to approximately solve their AL subproblem \eqref{eq:intro_al_primal_update}, initialized at $Z_{t-1}$. However, because the AL subproblems arising in \ourmethodplus have an unbounded feasible region, the FW variant used to solve them relies on substantially different techniques from the one used by HALLaR. Both FW variants generate matrices in factored form, i.e., $X=UU^{\top}$, where
the number of columns $r$ of $U$ varies during the algorithm but is usually much smaller than $n$.
% $U\in \mathbb R^{n\times r}$ and $r\ll n$.  
This low-rank representation reduces memory footprint from $O(n^2)$ to
$O(nr)$ and exploits sparsity in the SDP data \citep{burer2003nonlinear,burer2005local,monteiro2024low,
aguirre2025cuhallar}.
Both $X$ and $U$ are viewed as iterates in our discussion, with $U$ referred to as a factor of $X$.
% ??????.
% i $U$ denotes the factor of $X$.
% exploit the fact that their iterates are generated in factored form. 
% In particular, for an iterate $\tX=\tU\tU^T$, we refer to $\tU$ as a factor of $\tX$ and denote it by $U_{\tX}$. 
For a matrix $V$, $\mathrm{col}(V)$ denotes its number of columns.
Starting from a low-rank factor $\tilde U\in \mathbb R^{n\times r}$, an iteration of both FW variants consists of the following two basic steps:
\begin{enumerate}[label=\roman*),leftmargin=*,itemsep=0pt]
\item find a factor $U\in \mathbb R^{n\times r}$ such that $\mathrm{col}(U) \le \mathrm{col}(\tilde U)$ and
$\mathcal{L}_{\beta_{t-1}}(UU^{\top};p_{t-1})\le
\mathcal{L}_{\beta_{t-1}}(\tilde U\tilde U^{\top};p_{t-1})$;
\item without forming the associated full matrix, check whether $UU^{\top}$ satisfies the prescribed approximate optimality conditions for \eqref{eq:intro_al_primal_update}. If so, terminate and set $Z_t=UU^{\top}$; otherwise, perform an FW-type step based at $UU^{\top}$ to escape from this spurious solution $U$ to obtain a subsequent low-rank factor $\tilde{U}^{+}$ satisfying $\textrm{col}(\tilde U^{+})\le \textrm{col}(U)+1$, and hence $\textrm{col}(\tilde U^{+})\le \textrm{col}(\tilde U)+1$.
\end{enumerate}

We now discuss (i) in greater detail. The choice $U=\tilde U$ trivially satisfies step~i) and has been used by methods such as CGAL and T-CGAL \cite{yurtsever2019conditional, yurtsever2015scalable} (see also the methods in \cite{yurtsever2021scalable, ding2021optimal,shinde2021memory} for similar limited memory or reduced storage variants). Such methods can be viewed as essentially pure Frank-Wolfe methods, but have the disadvantage of producing iterates with high number of columns, since steps~(i)--(ii) must generally be repeated many times. On the other hand, to implement the second part of (i), HALLaR and HALLaR+ both use the nonlinear solver ADAP-AIPP, initialized at $\tilde U$, to find a near-stationary point $U$ of
$\min_{V\in\mathbb{R}^{n\times r}}
\mathcal{L}_{\beta_{t-1}}(VV^{\top};p_{t-1})$ satisfying the descent condition in step~i). Since the objective of this unconstrained subproblem is nonconvex and has only a locally Lipschitz continuous gradient, establishing the required complexity guarantees for ADAP-AIPP in this setting requires a more delicate analysis than the one used in HALLaR. The additional difficulty arises from the fact that, whereas the complexity of HALLaR depends on the diameter of the bounded feasible set, the complexity of HALLaR+ depends on a uniform bound on the diameters of all generated sublevel sets, some of which depend on both $\beta_{t-1}$ and $p_{t-1}$.

We next discuss (ii). Unlike HALLaR, which uses a computable FW gap to estimate whether the primal 
gap of 
$
% \mathcal L_{\beta_{t-1}}(\cdot;p_{t-1})-
\min_{X\in \mathbb S^{n}}\left\{\mathcal L_{\beta_{t-1}}(X;p_{t-1}): {X\succeq 0,\ \textrm{tr}(X)\leq 1}\right\} $ at $X=UU^T$ is small, 
HALLaR+ uses approximate optimality conditions 
tailored to the unbounded subproblem \eqref{eq:intro_al_updates}. 
Specifically, with $S=\nabla_X\mathcal{L}_{\beta_{t-1}}(UU^{\top};p_{t-1})$, 
it checks whether $\lambda_{\min}(S)$ is 
nearly nonnegative and $|\inner{S}{UU^{\top}}|$ is small. 
If either condition fails, HALLaR+ finds a 
rank-one matrix $y^F(y^F)^\top$ by solving $\min_{X\succeq 0} \left\{\inner{S}{X}: \mathrm{tr}(X)\leq 1\right\}$ 
and forms the next iterate as a conic 
combination of this matrix and the current 
iterate, with the two coefficients of the 
conic combination chosen by a two-
dimensional line search.

This approach for solving each AL subproblem 
is closely related to those of Laue \cite{laue2012hybrid} and Duchi et al. \cite{duchi2020conic,li2023conic},
which develop FW variants for the generic static problem of minimizing a convex function $\phi(\cdot)$ over $\mathbb S^{n}_{+}$. The two-dimensional line search for the conic coefficients used in our FW variant is motivated by the scheme of Laue. However, whereas Laue's convergence analysis is based on a primal-gap criterion that is not directly checkable in practice, our method uses the aforementioned tractable termination criteria. Duchi et al., on the other hand, propose a different unbounded FW variant based on two successive one-dimensional line-searches
rather than a single two-dimensional line search. Like our FW variant, their method terminates based on checkable approximate optimality conditions. While closely related to these methods, our FW variant requires a more involved analysis because it must account for the changing AL subproblems in \eqref{eq:intro_al_updates} and the dynamics of their associated sublevel sets, a difficulty absent from the static settings considered by Laue and Duchi et al.

\vspace{3mm}

\noindent\textbf{Our results and computational impact.}
Assuming that a primal--dual KKT pair exists and that the set of optimal solutions of \eqref{eq:intro_primal_sdp} is nonempty and bounded, it is shown that \ourmethodplus obtains an approximate primal--dual solution of \eqref{eq:intro_primal_sdp} and its dual with provable computational complexity bounds expressed in terms of parameters associated with the SDP instance and user-specified tolerances.
In particular, if the tolerances on primal feasibility, dual feasibility, and complementarity are all equal to $\varepsilon\in(0,1]$, then, \ourmethodplus performs $\mathcal O(\log(1/\varepsilon))$ AL iterations, $\mathcal O(\varepsilon^{-2})$ FW iterations/minimum eigenvector computations, and $\mathcal O(\varepsilon^{-3}))$ ADAP-AIPP iterations.

Our GPU implementation \cuourmethodplus is faster than cuHALLaR \citep{aguirre2025cuhallar} on 12 of 14 Hamming and 24 of 28 GSET maximum stable set instances, as well as all 26 matrix completion instances. On a single NVIDIA H200, it reaches $10^{-5}$ relative precision in about 44 seconds on a Hamming SDP ($n\approx 8.4\times10^6$, $m\approx 9.7\times10^7$) and 31 seconds on a matrix completion SDP ($n=2\times10^6$, $m\approx 2.6\times10^8$). On the other hand, cuLoRADS \citep{han2024accelerating} fails on the Hamming instance and takes about $4,100$ seconds on the matrix completion instance.

At $10^{-4}$ relative precision, \cuourmethodplus solves 80 of 81 QAPLIB relaxations, compared with 56 for CSDP \citep{borchers1999csdp}, 17 for HALLaR, and one for cuLoRADS. On the Mittelmann benchmarks \citep{mittelmann2003independent}, \cuourmethodplus solves 27 instances within one hour, cuLoRADS solves 15, and CARDAL \citep{li2026curvatureawarerankadaptivedistributedaugmentedlagrangian} solves 14. The final \cuourmethodplus iterate has rank two on every Hamming instance, at most $r+1$ on every matrix completion instance with ground-truth rank $r$, and 465 on the largest random quadratic assignment instance.

\vspace{3mm}

\noindent\textbf{Related works.}
Interior-point methods are widely used to solve small and medium-scale SDPs to high accuracy. Their applicability to large-scale instances is limited, however, because each Newton-type iteration typically requires solving a large linear system, often through the factorization of a large and potentially dense matrix
\cite{borchers1999csdp, Overton, Mon97, ToddToh}. We therefore focus below on methods specifically designed for large-scale SDPs.

SDPNAL+ \cite{zhao2010newton, yang2015sdpnal+} is an augmented Lagrangian (AL) method that approximately solves each AL subproblem using a semismooth Newton--CG method. Spectral bundle methods, and more generally bundle methods, have also been developed for large-scale SDPs (see, for example, \cite{Helmberg, Helmberg2, Helmberg3, DingGrimmer, Oustry, angell2024fast}). Another substantial line of work applies the alternating direction method of multipliers (ADMM) and related operator-splitting methods to large-scale SDPs
(e.g., \cite{odonoghue2016conic, garstka2021cosmo, deng2022new, madani2015admm, zheng2017fast, Kang, deng2022enhanced, wen2010alternating, o2016conic}).
Since methods based on nonconvex low-rank (LR) factorizations are most closely related to the approach developed in this paper, we discuss them in greater detail below. We refer the reader to the survey \cite{majumdar2020recent} for a broader overview of methods for solving large-scale SDPs.

The nonconvex LR approach introduced in \cite{burer2003nonlinear, burer2005local} underlies many successful algorithms for solving important classes of SDPs. In particular, the SDPLR method developed in these works applies an AL method to the smooth nonconvex reformulation
\begin{align}\label{eq:sdp-lr}
   \tag{$P_r$}
    %U_t \quad \approx\quad
   \min_U\quad
    \left \{ 
    C \bullet U U^T \quad:\quad \mathcal{A}(U U^T) = b, \quad
    \|U\|_F \leq 1, \quad
    U \in \mathbb R^{n\times r} \right\}
\end{align}
and approximately solves its AL subproblems using a limited-memory BFGS method. Although SDPLR is designed to handle equality constraints, it can be modified to accommodate inequality constraints (see, for example, \cite{kulis2007fast}).

Several methods solve \eqref{eq:sdp-lr} without employing an AL framework. For MaxCut SDPs, several specialized algorithms solve \eqref{eq:sdp-lr} using optimization techniques that preserve feasibility throughout their execution
(e.g., \cite{burer2001projected, homer1997design,
erdogdu2022convergence,mei2017solving, journee2010low}).
Riemannian optimization methods have also been employed for structured classes of SDPs for which the feasible set of \eqref{eq:sdp-lr} is a smooth manifold
(e.g., \cite{rosen2019se, journee2010low, huang2017solving, mei2017solving}).
Related feasibility-preserving methods based on manifold-optimization techniques have been developed for more challenging classes of SDPs in which the feasible sets of \eqref{eq:sdp-lr} are algebraic varieties \cite{TangToh, TangToh2, TangToh3}. Alternative feasibility-preserving strategies have also been explored in \cite{rosen2021scalable}.

More recently, several low-rank SDP solvers have combined low-rank factorization with augmented-Lagrangian-type schemes \citep{huang2024suboptimality,brosch2026augmented}. Like SDPLR, LoRADS and its GPU counterpart cuLoRADS \cite{han2024accelerating, han2025low} use the factorization $X=UU^{T}$ in their first phase and apply an inexact AL method to \eqref{eq:sdp-lr}. In their second phase, these methods adopt the bilinear factorization $X=UV^{T}$, penalize the additional constraint $U=V$, and apply an alternating-minimization scheme to the resulting augmented Lagrangian function. The ALORA method of \cite{DiAL} combines gradient descent with a spectral block Frank--Wolfe method to approximately minimize its nonconvex AL subproblems. CARDAL \citep{li2026curvatureawarerankadaptivedistributedaugmentedlagrangian} combines rank-adaptive AL iterations with matrix-free L-BFGS and negative-curvature corrections, provides slack-based approximate optimality conditions, and supports distributed GPU execution. SDPF+ \citep{tang2026preconditioned} combines a feasible low-rank method with an AL phase using a dynamically weighted quadratic penalty. 
Other low-rank methods focus on specially structured SDPs whose factorized nonconvex AL subproblems take the form $\argmin_{U\in \mathcal M} \mathcal{L}_{\beta_{t-1}}(UU^{\top};p_{t-1})$, where $\mathcal M$ is a Riemannian manifold.
Such subproblems are solved using either a Riemannian trust-region method \cite{ManiSDP, wang2025dual} or a Riemannian semismooth Newton method \cite{RiemannianSemiSmooth}. For SDP-RLT and polyhedral-SDP relaxations, respectively, RiNNAL+ and RiNNAL-POP \citep{hou2026rinnalplus,hou2026polyhedral} combine low-rank AL iterations with projected-gradient steps in the original matrix variable to adjust the rank and escape spurious factor-space local minima.

\vspace{3mm}

\noindent\textbf{Structure of the paper.}
This paper is organized into five sections and two appendices.
Section~2 introduces the pair of primal-dual SDP problems and develops a general AL framework for solving them.
Section~3 presents X-AL, an instance of this framework whose iterates are generated in the $\mathbb S_{+}^{n}$ space, and establishes its main complexity result for solving the SDP pair.
Section~4 introduces HALLaR+, a computationally and memory-efficient variant of X-AL implemented in the $U$-space. Theorem~\ref{thm:GeneralHLRComplexity}, the main theorem of the paper, establishes the complexity bounds for HALLaR+.
Section~5 reports computational experiments comparing HALLaR+ with several existing solvers on a large collection of SDPs.
Appendix~A provides technical convex-analysis results used throughout the paper and develops the relaxed unbounded Frank-Wolfe (RUFW) framework used by X-AL. Finally, Appendix~B presents the convergence analyses of an adaptive version of FISTA and an inexact proximal point method, both implemented in the U-space and used as subroutines within HALLaR+.

\subsection{Basic definitions and notation}

Let $\mathbb{E}$ be a finite-dimensional inner-product space. We denote its
inner product by $\inner{X}{Y}$ and the associated norm by
$\|X\|\defi\sqrt{\inner{X}{X}}$. When $\mathbb E=\mathbb S^n$, we use the
trace inner product
$\inner{X}{Y}\defi\operatorname{tr}(XY)$, whose associated norm is the
Frobenius norm. We use the notation
\[
% \log_2^+(t)\defi\max\{0,\log_2 t\}
\log_2^+(t)\defi\max\{0,\log_2 t\},\qquad \log_1^+(t)\defi\max\{1,\log t\},
\]
and write $\mathcal O_1(\Theta)$ as shorthand for
$\mathcal O(1+\Theta)$.
The operator norm of $\mathcal A$ is $\|\mathcal A\|\defi\max\{\|\mathcal AX\|:\|X\|\le1\}$.

We first recall some standard definitions from convex analysis. For a set
$\Omega\subseteq\mathbb E$, its indicator function is denoted by
$\delta_\Omega$, where $\delta_\Omega(X)=0$ if $X\in\Omega$ and
$\delta_\Omega(X)=\infty$ otherwise. If $\Omega$ is closed and convex, its
normal cone at $X\in\Omega$ is
\begin{equation}\label{def: normal cone}
 \mathcal N_\Omega(X)
 \defi
 \{V:\inner{V}{U-X}\le0, \ \forall \,U\in\Omega\}.
\end{equation}
For a proper closed convex function
$f:\mathbb E\to\mathbb R\cup\{\infty\}$, a point $X\in\operatorname{dom}f$,
and a scalar $\varepsilon\ge0$, the $\varepsilon$-subdifferential of $f$ at
$X$ is
\[
 \partial_\varepsilon f(X)
 \defi
 \{V:f(U)\ge f(X)+\inner{V}{U-X}-\varepsilon ,
       \ \forall  U\in\mathbb E\}.
\]
For a nonempty closed convex set $\Omega\subseteq\mathbb E$, the Euclidean
projection of $X\in\mathbb E$ onto $\Omega$ is
\begin{equation}\label{ProjOperator}
 \operatorname{proj}_\Omega(X)
 \defi
 \argmin_{Y\in\Omega}\|X-Y\|
\end{equation}
and the distance from a point $q\in\mathbb E$ to $\Omega$ is
$\operatorname{dist}(q,\Omega)\defi
\inf\{\|q-s\|:s\in\Omega\}$.

Given a set $\Omega\subseteq\mathbb E$, a point $X_0\in\Omega$, and a
function $f:\Omega\to\mathbb R$, we define the sublevel set of $f$ at $X_0$ and its radius by
\begin{align}
 \mathcal S(X_0;f)
 &\defi\{X'\in\Omega:f(X')\le f(X_0)\},                                      \label{SublevelDef}\\
D(X_0;f)
&\defi\sup\{\|X'\|:X'\in\mathcal S(X_0;f)\}.                              \label{DiameterDef}
\end{align}

We next introduce notation related to a differentiable function $g:\mathbb E\to\mathbb R$. Its affine
approximation at $\bar z\in\mathbb E$ is
\begin{equation}\label{eq:defell}
 \ell_g(z;\bar z)
 \defi
 g(\bar z)+\inner{\nabla g(\bar z)}{z-\bar z}
 \qquad\forall z\in\mathbb E.
\end{equation}
Function $g$ is said to have a locally
Lipschitz continuous gradient if, for every $U\in\mathbb E$, there exist
constants $R=R_U>0$ and $\bar L_U>0$ such that
\begin{equation}\label{eq:locally_lipschitz_gradient}
 \|\nabla g(U')-\nabla g(U'')\|
 \le \bar L_U\|U'-U''\|, \quad \forall \, U',U'' \in U + B(R)
\end{equation}
where $B(R)\defi\{X\in\mathbb E:\|X\|\le R\}$. The Lipschitz constant of
$\nabla g(\cdot)$  and the maximum of $\|\nabla g(\cdot)\|$ both relative to the ball $B(R)$ are denoted,
respectively, by
\begin{align}
 [{\rm LIP}](R;g)
 &\defi
 \inf\left\{
 L>0:
 \|\nabla g(X)-\nabla g(X')\|\le L\|X-X'\|
 \quad\forall X,X'\in B(R)
 \right\},                                 \label{SmoothnessDef}\\
 \nabla(R;g)
 &\defi\sup\{\|\nabla g(U)\|:U\in B(R)\}.                                 \label{GradDef}
\end{align}
\section{The SDP problem and an inexact AL framework}
\label{section: AL general framework}

This section develops an inexact augmented Lagrangian (AL) framework for computing approximate primal-dual solutions of an SDP. 
The framework is independent of the particular method used to solve the AL subproblems. At each outer iteration, the inner solver need only return a PSD matrix satisfying prescribed approximate first-order optimality conditions for the current AL subproblem. The multiplier is then updated, and the penalty parameter is increased whenever the primal residual remains above the target tolerance. This separation allows the outer-iteration complexity of the AL framework to be analyzed independently of the specific inner solver used to satisfy the required inexactness conditions.

The section is organized as follows. Subsection~\ref{subsec:SDPproblem} introduces the primal-dual SDP pair, discusses its assumptions, and specifies the approximate primal-dual solution concept targeted by the framework. Subsection~\ref{subsec:ial-framework} formally presents the inexact AL framework. Finally, Subsection~\ref{subsec:hallar-proof} proves Theorem~\ref{thm:complexity_relaxed_corrected},
which gives the main outer-iteration complexity bound for the framework.

\subsection{The SDP problem and assumptions}\label{subsec:SDPproblem}

We consider the primal-dual pair of SDPs
\begin{equation}
\label{eqn:general_probl}
\min_{X \in \mathbb{S}_+^n} \left\{\inner{C}{X} \; : \; \mathcal{A} X = b \right\}
\end{equation}
and
\begin{equation}\label{eqn:lagrangian_dual}
\max\left\{\inner{b}{p}:\ C-\A^*p\succeq0,\ \ p \in \mathbb{R}^{m}\right\}.
\end{equation}
In contrast to the formulation in \cite{monteiro2024low}, which assumes a bounded feasible region, problem~\ref{eqn:general_probl} is posed over the full PSD cone $\mathbb S_{+}^{n}$, without a constraint of the form $\operatorname{tr}(X)\le \tau$. In the absence of such a bound, the analysis instead relies on properties of the SDP solution set.

We now state the assumptions that are made on the SDP data.
\begin{assumption}
\label{ass:problem-setup}
The set of optimal solutions of \eqref{eqn:general_probl} is nonempty and bounded.
\end{assumption}
\begin{assumption}\label{ass:existence-KKT}
There exists an optimal primal dual pair $(X_*,p_*)\in\mathbb{S}_+^n\times\mathbb{R}^{m}$ satisfying
\begin{equation}
\label{eqn:optimality_conditions_general}
0 \ \in \ C - \mathcal{A}^{*} p_* + \mathcal{N}_{\mathbb{S}_+^n}(X_*)\,,\qquad \mathcal{A} X_* = b.
\end{equation}
\end{assumption}
The definition of the normal cone in \eqref{def: normal cone}, applied with $\Omega=\mathbb S_+^n$ and $X=X_*$ implies that \eqref{eqn:optimality_conditions_general} is equivalent to
\begin{equation}\label{eq:KKT}
\A X_*=b,\qquad X_*\succeq 0,\qquad C-\mathcal A^{*}p_{*}\succeq 0,\qquad \inner{C-\mathcal A^{*}p_{*}}{X_*}=0.
\end{equation}

Relation \eqref{eq:KKT} expresses primal feasibility, dual feasibility, and
complementarity. In particular, $(X_{*},p_{*})$
is a primal-dual KKT pair \eqref{eqn:general_probl}--\eqref{eqn:lagrangian_dual}, and the primal and
dual problems have zero duality gap at this pair.

The AL framework will terminate when approximate counterparts of the three
nontrivial conditions in \eqref{eq:KKT} have been obtained. The next definition preciesely characterizes our notion of inexactness.

\begin{definition}\label{def:EPD-solution}
For tolerances $(\varepsilon_p,\varepsilon_c,\varepsilon_d)>0$,
a pair $(\bar Z,\bar p)\in\mathbb S_+^n\times\mathbb R^m$ is called an
$(\varepsilon_p,\varepsilon_c,\varepsilon_d)$-solution of
\eqref{eqn:general_probl}--\eqref{eqn:lagrangian_dual} if
\begin{equation}
\label{eqn:framework-goal}
\|\mathcal A(\bar Z)-b\|\le \varepsilon_p,
\qquad
\lambda_{\min}(\bar S)\ge -\varepsilon_d,
\qquad
\left|\inner{\bar S}{\bar Z}\right|\le \varepsilon_c,
\end{equation}
where $\bar S:=S(\bar p)$ and
\begin{equation}\label{def:dual_slack_map}
S(p)\defi C-\mathcal A^*p.
\end{equation}
\end{definition}

Definition~\ref{def:EPD-solution} is a direct relaxation of the exact KKT
conditions. Indeed, setting
$\varepsilon_p=\varepsilon_c=\varepsilon_d=0$ in
\eqref{eqn:framework-goal} recovers \eqref{eq:KKT}.

\subsection{AL Framework}\label{subsec:ial-framework}
This subsection presents the AL framework for approximately solving the primal-dual SDP pair \eqref{eqn:general_probl}--\eqref{eqn:lagrangian_dual}. For ease of analysis and presentation, the AL framework is presented in the $X$-space, meaning that all iterates are full $n\times n$ positive semidefinite matrices. Our ultimate goal is to later present an algorithm that never works with such dense matrices $X$ and instead stores only low-rank factors $U$, which typically have only a few columns. This memory-efficient low-rank algorithm, HALLaR+, is presented in Section~\ref{sec:u-space-hallar} and can be viewed as a variant of the AL framework.

The AL framework is based on approximately solving a sequence of augmented
Lagrangian subproblems \eqref{eq:intro_al_updates}. For convenience, we define the shorthand \begin{equation}\label{def:primal_residual}
r(X)\defi b-\A X\qquad\forall X\in\mathbb{S}^{n}.
\end{equation}
Hence, for a multiplier $p\in\mathbb R^m$ and a penalty
parameter $\beta\in\mathbb R_+$, the augmented Lagrangian function is defined as

\begin{equation}
\label{eqn:AL_general}
\mathcal{L}_{\beta}(X;p) \ \defi \ \inner{C}{X} + \inner{p}{r(X)} + \frac{\beta}{2} \|r(X)\|^2.
\end{equation}
Moreover, during the $t$-th iteration of the AL framework, the augmented Lagrangian subproblem \eqref{eq:intro_al_primal_update} can be written as
\[
\min_{X \in \mathbb{S}_+^n}\mathcal{L}_{\beta_{t-1}}(X;p_{t-1})
\ =\
\min_{X \in \mathbb{S}_+^n}
\left\{\inner{C}{X}+\inner{p_{t-1}}{r(X)}
+\frac{\beta_{t-1}}{2}\|r(X)\|^2\right\}.
\]
Since the AL objective in \eqref{eq:intro_al_primal_update} is convex and differentiable and
the feasible region is the closed convex cone $\mathbb S_+^n$, its
first-order optimality condition is
\begin{equation}\label{eq:dbeta_optimality_condition}
0\in \nabla_X \mathcal L_{\beta_{t-1}}(X;p_{t-1})+\mathcal N_{\mathbb S_+^n}(X).
\end{equation}
By the normal-cone characterization in \eqref{def: normal cone},
condition~\eqref{eq:dbeta_optimality_condition} is equivalent to
\begin{equation}\label{eq:dbeta_optimality_condition_equiv}
\nabla_X\mathcal L_{\beta_{t-1}}(X;p_{t-1}) \succeq 0,
\qquad
\inner{\nabla_X\mathcal L_{\beta_{t-1}}(X;p_{t-1})}{X}=0.
\end{equation}
The AL framework in this subsection approximately solves each AL subproblem based on stopping criteria that are direct relaxations of the optimality conditions in \eqref{eq:dbeta_optimality_condition_equiv} (see \eqref{eq:inner_comp} below). This differs from the HALLaR and CGAL methods in \cite{monteiro2024low,yurtsever2019conditional}, which assume a compact feasible region and approximately solve their AL subproblems by performing a prescribed number of Frank--Wolfe iterations or by terminating when the Frank--Wolfe gap is sufficiently small. These approaches are not directly applicable to \eqref{eq:intro_al_primal_update}, whose feasible region is unbounded.

The AL framework is now presented.

\noindent
\begin{minipage}[t]{1\columnwidth}
\rule[0.5ex]{1\columnwidth}{1pt}

\noindent \textbf{AL Framework}

\noindent \rule[0.5ex]{1\columnwidth}{1pt}%
\end{minipage}

\noindent \textbf{Input}: Initial pair $(Z_0,p_0)\in\mathbb{S}^n_+\times\R^m$,
tolerance triple $(\varepsilon_p,\varepsilon_c,\varepsilon_d)\in\R_{++}^3$,
and initial penalty parameter $\beta_0\in\R_{++}$.

\noindent \textbf{Output}: Pair $(\bar Z,\bar p)$ satisfying \eqref{eqn:framework-goal}.

\begin{itemize}[leftmargin=2.0em]
\item[{\bf 0.}] Set $t\gets 1$.
\item[{\bf 1.}] Find $Z_t\in\mathbb{S}^{n}_{+}$, which together with 
$S_t=C-\cA^*[p_{t-1}+ \beta_{t-1} r(Z_t)]$,
satisfies
\begin{align}
\lambda_{\min}\!\left(S_{t}\right) \ge -\varepsilon_{d},\qquad
\left|\inner{S_{t}}{Z_{t}}\right|
\le\ \varepsilon_c.
\label{eq:inner_comp}
\end{align}

\item[{\bf 2.}] Set
$p_{t} = p_{t-1}+\beta_{t-1}r(Z_{t})$.

\item[{\bf 3.}] \textbf{If} $\|r(Z_t)\|\le\varepsilon_p$, \textbf{then stop} and output $(\bar Z,\bar p)\gets(Z_t,p_t)$.
\item[{\bf 4.}] set $\beta_{t}\gets 2\beta_{t-1}$, $t\gets t+1$, and go to \textbf{Step}~1.
\end{itemize}
\noindent\rule[0.5ex]{1\columnwidth}{1pt}
\vspace{1em}

First, Step~{\bf 1} computes $Z_t\in\mathbb S_+^n$ satisfying an inexact
optimality condition for the current AL subproblem. Indeed, the identity
\begin{equation}\label{eqn:AL_gradient}
\nabla_X\mathcal L_{\beta}(X;p)
=
C-\A^*[p+\beta r(X)]
\end{equation}
implies that
$\nabla_X\mathcal L_{\beta_{t-1}}(Z_t;p_{t-1})=S_t$.
Thus, the conditions in \eqref{eq:dbeta_optimality_condition_equiv} are a special case of \eqref{eq:inner_comp} with $\varepsilon_c=\varepsilon_d=0$
The AL framework requires only that \eqref{eq:inner_comp} hold and does not
prescribe how such a point is computed. A Frank--Wolfe-type method for
computing such a point is developed in later sections.

Second, Step~{\bf 2} performs the standard AL multiplier update, which, together
with the definition of $S(\cdot)$ in \eqref{def:dual_slack_map},
implies that $S_t=S(p_t)$. The two conditions in \eqref{eq:inner_comp} then
ensure that $(Z_t,p_t)$ satisfies the approximate dual-feasibility and
complementarity requirements in \eqref{eqn:framework-goal}. Third, because $Z_t\in\mathbb S_+^n$ by construction, Step~{\bf 3} only
needs to check primal feasibility. If
$\|r(Z_t)\|\le\varepsilon_p$, then $(Z_t,p_t)$ satisfies all the
conditions in \eqref{eqn:framework-goal} in view of \eqref{eq:inner_comp}. In this case $(Z_t,p_t)$ is returned, otherwise, if $\|r(Z_t)\|>\varepsilon_p$, Step~{\bf 4} doubles the penalty
parameter and proceeds to the next iteration. Consequently, for every
$t\ge1$, the penalty parameters satisfy

\begin{equation}\label{eqn:beta_identities_corrected}
\beta_{t-1}=2^{t-1}\beta_0.
\end{equation}
For the analysis, the same formula defines the penalty parameter at the
terminal index.

Theorem~\ref{thm:complexity_relaxed_corrected} gives the main
outer-iteration complexity bound for the AL framework. Its proof is
presented in Subsection~\ref{subsec:hallar-proof}.

\begin{theorem}\label{thm:complexity_relaxed_corrected}
Let
$\varepsilon\defi(\varepsilon_p,\varepsilon_c,\varepsilon_d)>0$ be given and define
\begin{equation}
\label{eq:bar_beta_tolerance_bound_corrected}
\begin{aligned}
C_1 = C_1(\varepsilon) 
\defi
\frac{\varepsilon_c+\varepsilon_d}{\varepsilon_p}, \qquad 
&\quad\tau_*\ \defi\ \max\{1,2\operatorname{tr}(X_*)\}, 
&\quad\bar\beta(\varepsilon)
&\defi
\beta_0
+\frac{1}{\varepsilon_p}
\left[
4\left(\|p_0\|+\|p_*\|\right)+16C_1\tau_*
\right].
\end{aligned}
\end{equation}
Then the following statements about the AL framework hold:
\begin{itemize}
    \item[i)] it terminates in at most $T$ iterations, where
\begin{equation}
\label{eq:T_def_corrected}
T=T_\varepsilon
 \defi
1+\left\lceil
\log_2^+\!\left(
\frac{
2\left(\|p_0\|+\|p_*\|\right)+8C_1\tau_*
}{\beta_0\varepsilon_p}
\right)\right\rceil;
\end{equation}
\item[ii)] its last generated penalty parameter $\bar\beta$ is bounded by 
\begin{equation}
\label{eq:beta_last_bound_corrected}
\bar\beta
\le
\bar\beta(\varepsilon)
\end{equation}
where the bound $\bar\beta(\varepsilon)$ is defined in \eqref{eq:bar_beta_tolerance_bound_corrected};
\item[iii)] its output pair $(\bar Z,\bar p)$ is an
$(\varepsilon_p,\varepsilon_c,\varepsilon_d)$-solution of
\eqref{eqn:general_probl}--\eqref{eqn:lagrangian_dual}, and every multiplier $p_t$ generated by the AL framework satisfies
\begin{equation}
\label{eq:multiplier_bound_corrected}
\begin{aligned}
\|p_t\|
&\le
\|p_*\|
+\sqrt{
\|p_0-p_*\|^2
+4\varepsilon_p C_1\bar\beta(\varepsilon)\tau_*
}.
\end{aligned}
\end{equation}

\end{itemize}
\end{theorem}

\subsection{Proof of Theorem~\ref{thm:complexity_relaxed_corrected}}
\label{subsec:hallar-proof}
Besides Assumption~\ref{ass:existence-KKT}, the proof below uses only relation \eqref{eq:inner_comp} and Steps~{\bf 2}--{\bf 4} of the AL framework, and hence it does not depend on how Step~{\bf 1} is implemented. For every $t\ge1$, define the gap quantity
\begin{equation}\label{eqn:vartheta_def_new}
\vartheta_t^*
\ \defi\
\max_{U\in\Delta_{\tau_*}}\inner{S_t}{Z_t-U},
\end{equation}
where $\Delta_{\tau}$ is as in \eqref{eq:Delta_tau_def} and $\tau_*$ is as in \eqref{eq:bar_beta_tolerance_bound_corrected}. We first show that the inner product between the dual error $p_t-p_*$ and the primal residual $r(Z_t)$ is controlled by $\vartheta_t^*$. This relation will be instrumental in bounding the sequence of primal residuals $\|r(Z_t)\|$.

\begin{lemma}
\label{lem:inner-product-bound-ial}
For every $t \ge 1$, the iterate pair $(Z_t,p_t)$ satisfies
\begin{equation}
\label{eq:ip-resid-ial}
\inner{p_t - p_*}{r(Z_t)} \ \le \ \vartheta_t^*,
\end{equation}
where $\vartheta_t^*$ is as in \eqref{eqn:vartheta_def_new}.
\end{lemma}

\begin{proof}
Using \eqref{eqn:vartheta_def_new}, the inclusion $X_* \in \Delta_{\tau_*}$, and the fact that the identity $p_t=p_{t-1}+\beta_{t-1}r(Z_t)$ and \eqref{eqn:AL_gradient} imply that $S_{t}=C-\mathcal A^{*}p_t$,
we have
\begin{align}
\vartheta_t^*
&\overset{\eqref{eqn:vartheta_def_new}}{=}\ \max_{U \in \Delta_{\tau_*}} \inner{S_{t}}{Z_t - U}
\geq \ \inner{S_{t}}{Z_t - X_*} \overset{\eqref{eqn:AL_gradient}}{=}\inner{C - \A^* p_t}{Z_t - X_*} \notag\\
&= \inner{\A^*(p_* - p_t)}{Z_t - X_*} + \inner{C - \A^* p_*}{Z_t - X_*} \label{eq:ial_vartheta_identity}.
\end{align}
Using \eqref{eqn:optimality_conditions_general} and the definition in \eqref{def:primal_residual}, the first term in \eqref{eq:ial_vartheta_identity} 
satisfies
\begin{align}
\inner{\A^*(p_* - p_t)}{Z_t - X_*}
&=\ \inner{p_* - p_t}{\A(Z_t - X_*)}
\ \overset{\eqref{eqn:optimality_conditions_general}}{=}\ \inner{p_* - p_t}{\A Z_t - b} \notag\\
&=\ \inner{p_t - p_*}{b - \A Z_t}
\ \overset{\eqref{def:primal_residual}}{=}\ \inner{p_t - p_*}{r(Z_t)}.
\label{eq:ial_first_term_identity}
\end{align}
It is then easy to see from 
\eqref{eqn:optimality_conditions_general}
that $C-\mathcal A^{*}p_{*}\in -\mathcal{N}_{\mathbb S_{+}^{n}}(X_{*})$. This inclusion, the fact that $Z_t\in \mathbb S_{+}^{n}$, and the definition of normal cone in \eqref{def: normal cone} with $u=Z_t$, $v=-(C-\mathcal A^{*}p_{*})$, $\Omega=\mathbb S_{+}^{n}$, and $x=X_{*}$ imply that
\begin{equation}\label{eq:ial_second_term_nonnegative}
\inner{C - \A^* p_*}{Z_t - X_*}\ge\ 0.
\end{equation}
Relation \eqref{eq:ip-resid-ial} then follows immediately from relations \eqref{eq:ial_vartheta_identity}, \eqref{eq:ial_first_term_identity}, and \eqref{eq:ial_second_term_nonnegative}.
\end{proof}

The following lemma bounds the quantity $\vartheta_t^*$ by a fixed quantity that does not vary with the iteration index $t$.

\begin{lemma}\label{lem:fw-gap-translation-corrected}
For every $t\ge1$, we have
\begin{equation}\label{eq:fw_gap_lp_translation}
\vartheta_t^*
\ \le\
\varepsilon_c+\varepsilon_d\tau_*,
\end{equation}
where $\vartheta_t^*$ is as in \eqref{eqn:vartheta_def_new}.
\end{lemma}
\begin{proof}
Definition \eqref{eqn:vartheta_def_new} and
Proposition~\ref{prop:trace_ball_gap_residual_bounds}(a), applied with
$(X,G,\tau)=(Z_t,S_t,\tau_*)$, imply that
\begin{align*}
\vartheta_t^*
&\le
|\inner{S_t}{Z_t}|
+\max\{\lambda_{\max}(-S_t),0\}\tau_* \\
&\le
\varepsilon_c+\varepsilon_d\tau_*,
\end{align*}
where the last inequality follows from \eqref{eq:inner_comp} and the
identity $\lambda_{\max}(-S_t)=-\lambda_{\min}(S_t)$.
\end{proof}

The next lemma establishes key bounds on the sequence of Lagrange multipliers and primal feasibility residuals.

\begin{lemma}\label{lem:main_bounds_relaxed-ial_corrected}
For every generated iteration index $t\ge1$, we have
\begin{equation}\label{eq:key_iterate_bound_lp_corrected}
\|p_t-p_*\|^2+\beta_{t-1}^2\|r(Z_t)\|^2
\le \|p_0-p_*\|^2+4\varepsilon_pC_1\beta_{t-1}\tau_*.
\end{equation}
\end{lemma}
\begin{proof}
The multiplier update in Step~{\bf 2} implies
\begin{align}
\|p_t-p_*\|^2+\beta_{t-1}^2\|r(Z_t)\|^2
&=\|p_{t-1}-p_*\|^2
+2\beta_{t-1}\inner{p_t-p_*}{r(Z_t)}\notag\\
&\le\|p_{t-1}-p_*\|^2
+2\varepsilon_pC_1\beta_{t-1}\tau_*,
\label{eq:ial_one_step_energy_corrected}
\end{align}
where the inequality follows from \eqref{eq:ip-resid-ial},
\eqref{eq:fw_gap_lp_translation}, and
\eqref{eq:bar_beta_tolerance_bound_corrected}.
For $t=1$, relation \eqref{eq:key_iterate_bound_lp_corrected} follows directly
from \eqref{eq:ial_one_step_energy_corrected}. For $t\ge2$, applying
\eqref{eq:key_iterate_bound_lp_corrected} at index $t-1$ to the right-hand
side of \eqref{eq:ial_one_step_energy_corrected} and using
$\beta_{t-2}=\beta_{t-1}/2$ proves the induction step.
\end{proof}

\begin{proof}[Proof of Theorem~\ref{thm:complexity_relaxed_corrected}]
For convenience, we define
\begin{equation}\label{eq:Psi_def_corrected}
\Psi_* = \Psi_*(\varepsilon)
\defi
2\left(\|p_0\|+\|p_*\|\right)+8C_1\tau_*.
\end{equation}
Using the fact that $y\le2^{\lceil\log_2^+y\rceil}\le1+2y$ for every $y>0$
with $y=\Psi_*/(\beta_0\varepsilon_p)$, the definition of $T$ in
\eqref{eq:T_def_corrected}, relation
\eqref{eqn:beta_identities_corrected}, and the definition of
$\bar\beta(\varepsilon)$ in \eqref{eq:bar_beta_tolerance_bound_corrected}, we obtain
\begin{equation}\label{eq:T_beta_bounds_corrected}
\bar\beta(\varepsilon)
\ =\
\beta_0+\frac{2\Psi_*}{\varepsilon_p}
\ \ge\
\beta_02^{T-1}
\ =\
\beta_{T-1}
\ \ge\
\frac{\Psi_*}{\varepsilon_p}
\ \overset{\eqref{eq:Psi_def_corrected}}{\ge}\
\frac{1}{\varepsilon_p}
\max\left\{
2\left(\|p_0\|+\|p_*\|\right),\ 8C_1\tau_*
\right\}.
\end{equation}

Suppose now, for the sake of contradiction, that the AL framework has not
terminated at the end of iteration $T$, i.e., that $\|r(Z_T)\|>\varepsilon_p$
in Step~{\bf 3}. Then, dropping the nonnegative term $\|p_T-p_*\|^2$ from the
left-hand side of \eqref{eq:key_iterate_bound_lp_corrected} with $t=T$ and
dividing the resulting inequality by $\beta_{T-1}^2\varepsilon_p^2$, we
obtain
\begin{align}
1
\ <\
\frac{\|r(Z_T)\|^2}{\varepsilon_p^2}
&\overset{\eqref{eq:key_iterate_bound_lp_corrected}}{\le}
\frac{\|p_0-p_*\|^2}{\varepsilon_p^2\beta_{T-1}^2}
+\frac{4C_1\tau_*}{\varepsilon_p\beta_{T-1}}
\ \le\
\frac{\left(\|p_0\|+\|p_*\|\right)^2}{\varepsilon_p^2\beta_{T-1}^2}
+\frac{4C_1\tau_*}{\varepsilon_p\beta_{T-1}}
\overset{\eqref{eq:T_beta_bounds_corrected}}{\le}
\frac14+\frac12
\ =\
\frac34,
\label{eq:contradiction_target_corrected}
\end{align}
where the second inequality uses the triangle inequality and the last one
bounds the two fractions by means of the first and second terms in the
maximum in \eqref{eq:T_beta_bounds_corrected}, respectively. This
contradiction shows that the AL framework performs at most $T$ outer
iterations.

Let $\bar t\le T$ be the terminal outer index, so that
$(\bar Z,\bar p)=(Z_{\bar t},p_{\bar t})$ and $\bar\beta=\beta_{\bar t-1}$.
Since $\{\beta_t\}$ is nondecreasing by relation
\eqref{eqn:beta_identities_corrected}, we have
$\bar\beta=\beta_{\bar t-1}\le\beta_{T-1}$, and hence the upper bound in
\eqref{eq:T_beta_bounds_corrected} yields
\eqref{eq:beta_last_bound_corrected}. Moreover, it is easy to see that
\eqref{eqn:AL_gradient} implies that
$\bar S=C-\A^*\bar p=S_{\bar t}$, and thus the two inequalities in
\eqref{eq:inner_comp} with $t=\bar t$ together with the stopping test in
Step~{\bf 3} give the three relations in \eqref{eqn:framework-goal}. Hence,
$(\bar Z,\bar p)$ is an $(\varepsilon_p,\varepsilon_c,\varepsilon_d)$-solution of
\eqref{eqn:general_probl}--\eqref{eqn:lagrangian_dual}.

Finally, let $t\le\bar t$ be any generated iteration index. Dropping the
nonnegative term $\beta_{t-1}^2\|r(Z_t)\|^2$ from the left-hand side of
\eqref{eq:key_iterate_bound_lp_corrected} and using the facts that
$\beta_{t-1}\le\beta_{\bar t-1}=\bar\beta$ and
$\bar\beta\le\bar\beta(\varepsilon)$, we conclude that
\[
\|p_t-p_*\|^2
\ \overset{\eqref{eq:key_iterate_bound_lp_corrected}}{\le}\
\|p_0-p_*\|^2+4\varepsilon_pC_1\beta_{t-1}\tau_*
\ \overset{\eqref{eq:beta_last_bound_corrected}}{\le}\
\|p_0-p_*\|^2+4\varepsilon_pC_1\bar\beta(\varepsilon)\tau_*.
\]
Relation \eqref{eq:multiplier_bound_corrected} now follows by taking square
roots in the above inequality and using the triangle inequality
$\|p_t\|\le\|p_*\|+\|p_t-p_*\|$.
\end{proof}

\section{X-AL Method}\label{subsec:xspace-step-one}

This section presents a concrete method for approximately solving the
linearly constrained SDP \eqref{eqn:general_probl} under
Assumption~\ref{ass:problem-setup}, and establishes its total
inner-iteration complexity. The method is a special instance of the AL
framework of Section~\ref{section: AL general framework} in which the
Relaxed Unbounded Frank--Wolfe routine of Appendix~\ref{sec:RUFW} implements Step~{\bf 1}.
A feature of the resulting bound is that the number of inner iterations
performed during an outer iteration does not depend on how suboptimal the warm
start $Z_{t-1}$ is for the current AL subproblem.

\subsection{X-AL Method and Main Complexity Result}

We now make Step~{\bf 1} of the AL framework clear while keeping the
description entirely in the $X$-space. At outer iteration $t\ge1$, the
$X$-space implementation of Step~{\bf 1} applies the relaxed unbounded
Frank--Wolfe framework (RUFW) of Appendix~\ref{sec:RUFW} directly to the AL
subproblem
\begin{equation}
\label{eq:xspace_inner_subproblem}
\min\left\{\mathcal L_{\beta_{t-1}}(X;p_{t-1}):
X\in\mathbb S_+^n\right\},
\end{equation}
with initial point $Z_{t-1}$ and tolerance pair
$(\varepsilon_d,\varepsilon_c)$. The output returned by this call is denoted by the matrix $Z_t$.

We are now ready to present the X-Space AL Method.

\vspace{0.75em}
\noindent
\begin{minipage}[t]{1\columnwidth}
\rule[0.5ex]{1\columnwidth}{1pt}

\noindent \textbf{$X$-AL Method}

\noindent \rule[0.5ex]{1\columnwidth}{1pt}
\end{minipage}

\noindent \textbf{Input}: Initial pair $(Z_0,p_0)\in\mathbb S_+^n\times\R^m$,
tolerances $(\varepsilon_p,\varepsilon_c,\varepsilon_d)\in\R_{++}^3$, and
initial penalty parameter $\beta_0\in\R_{++}$.

\noindent \textbf{Output}: Pair $(\bar Z,\bar p)$ satisfying
\eqref{eqn:framework-goal}.

\begin{itemize}[leftmargin=2.0em]
\item[{\bf 0.}] Set $t\gets1$.
\item[{\bf 1.}] Perform the following $X$-space relaxed unbounded
Frank--Wolfe (RUFW) call:
\begin{enumerate}[leftmargin=3.2em]
\item[{\bf 1.a)}] Set $\tilde X\gets Z_{t-1}$.
\item[{\bf 1.b)}] Compute $X\in\mathbb S_+^n$ such that
\begin{equation}\label{eq:xspace_RUFW_correction}
\mathcal L_{\beta_{t-1}}(X;p_{t-1})
\le
\mathcal L_{\beta_{t-1}}(\tilde X;p_{t-1})
\end{equation}
and set
\[
S\gets\nabla_X\mathcal L_{\beta_{t-1}}(X;p_{t-1}).
\]
\item[{\bf 1.c)}] Compute a minimum
eigenpair $(\lambda_{\min}(S),v_{\min}(S))$ of $S$ with
$\|v_{\min}(S)\|=1$.
\item[{\bf 1.d)}] \textbf{If}
\begin{equation}\label{eq:xspace_RUFW_stop}
\lambda_{\min}(S)\ge-\varepsilon_d,
\qquad
\left|\inner{S}{X}\right|\le\varepsilon_c,
\end{equation}
\textbf{then} set $Z_t\gets X$ and go to Step~{\bf 2}.
\item[{\bf 1.e)}] Set
\begin{equation}\label{eq:xspace_RUFW_oracle}
y^F
\ \defi\
\begin{cases}
v_{\min}(S), & \lambda_{\min}(S)<0,\\
0, & \lambda_{\min}(S)\ge0,
\end{cases}
\qquad
X^F\ \defi\ y^F(y^F)^\top 
\end{equation}
and compute
\begin{equation}\label{eq:xspace_RUFW_linesearch}
(\alpha,\omega)
\in
\argmin\left\{
\mathcal L_{\beta_{t-1}}(\alpha' X+\omega' X^F;p_{t-1}):
\alpha'\ge0,\ \omega'\ge0
\right\}.
\end{equation}
\item[{\bf 1.f)}] Set
\begin{equation}\label{eq:xspace_RUFW_update}
\tilde X\gets
\alpha X+\omega X^F,
\end{equation}
and return to Step~{\bf 1.b)}.
\end{enumerate}

\item[{\bf 2.}] Set
\begin{equation}\label{eq:xspace_multiplier_update}
p_t\gets p_{t-1}+\beta_{t-1}r(Z_t),
\qquad
S_t\gets C-\A^*p_t.
\end{equation}
\item[{\bf 3.}] \textbf{If} $\|r(Z_t)\|\le\varepsilon_p$, \textbf{then} set
$(\bar Z,\bar p)\gets(Z_t,p_t)$, output $(\bar Z,\bar p)$, and \textbf{stop}.
\item[{\bf 4.}] Set
\begin{equation}\label{eq:xspace_penalty_update}
\beta_t\gets2\beta_{t-1},
\qquad
t\gets t+1,
\end{equation}
and return to Step~{\bf 1}.
\end{itemize}
\noindent\rule[0.5ex]{1\columnwidth}{1pt}
\vspace{0.75em}

Step~{\bf 1} approximately solves the current AL subproblem by means of
an $X$-space relaxed unbounded Frank--Wolfe (RUFW) procedure described in Appendix~\ref{sec:RUFW}. At the
beginning of the $t$-th outer iteration, Step~{\bf 1.a)} initializes the
RUFW iterate at $Z_{t-1}$. Given the current point $\widetilde X$,
Step~{\bf 1.b)} computes a point $X\in\mathbb S_+^n$ whose AL value is
no larger than that of $\widetilde X$, and evaluates the corresponding
gradient
\[
S=\nabla_X\mathcal L_{\beta_{t-1}}(X;p_{t-1}).
\]
Step~{\bf 1.c)} then computes a minimum eigenpair of $S$, and Step~{\bf 1.e)} forms the rank-one matrix
$X^F=y^F(y^F)^\top$ associated with a direction of negative curvature
of the linearized objective. Step~{\bf 1.d)} tests whether $X$ satisfies
the given approximate optimality conditions
\eqref{eq:xspace_RUFW_stop}. If so, the inner loop terminates and
$Z_t$ is set equal to $X$. Otherwise, Step~{\bf 1.e)} minimizes the AL
objective over nonnegative combinations of $X$ and $X^F$, and
Step~{\bf 1.f)} uses the resulting point as the next RUFW iterate. The inner loop in Step~{\bf 1} exits at Step~{\bf 1.d)} precisely when
the termination criterion \eqref{eq:xspace_RUFW_stop} holds, at which
point the method sets $Z_t=X$ and proceeds to Step~{\bf 2}. Hence, upon
exiting Step~{\bf 1},
$
S=\nabla_X\mathcal L_{\beta_{t-1}}(Z_t;p_{t-1}).
$

Step~{\bf 2} then performs the standard AL multiplier update. Combining the fact that relations \eqref{eqn:AL_gradient} and
\eqref{eq:xspace_multiplier_update} imply that
\begin{equation}
\label{eq:xspace_gradient_identity}
\nabla_X\mathcal L_{\beta_{t-1}}(Z_t;p_{t-1})
=C-\A^*\!\left[p_{t-1}+\beta_{t-1}r(Z_t)\right]
=C-\A^*p_t
=S_t.
\end{equation}
with the stopping conditions in
\eqref{eq:xspace_RUFW_stop} shows that the pair $(Z_t,p_t)$ produced
after Step~{\bf 2} always satisfies \eqref{eq:inner_comp}. Hence, Step~{\bf 1} implements the approximate
AL-subproblem solve required in Step~{\bf 1} of the AL framework
presented in Subsection~\ref{subsec:ial-framework}. The remaining steps
carry out the corresponding operations of that framework.
Step~{\bf 2} updates the multiplier and defines $S_t$, Step~{\bf 3}
applies the primal-feasibility stopping test, and Step~{\bf 4} increases
the penalty parameter before the next outer iteration.

The following result provides an iteration-complexity bound for the X-AL
method and will be proved in the next subsection.
\begin{theorem}\label{thm:total_AL_X_complexity}
Set $\varepsilon\defi(\varepsilon_p,\varepsilon_c,\varepsilon_d)$.
Let $C_1=C_1(\varepsilon)$ and $\tau_*$ be as in
\eqref{eq:bar_beta_tolerance_bound_corrected}, and let $T=T_\varepsilon$
and $\bar\beta(\varepsilon)$ be as in \eqref{eq:T_def_corrected} and
\eqref{eq:bar_beta_tolerance_bound_corrected}, respectively. Define
$\bar\kappa$ and $\bar\tau$ by
\begin{align}
\bar\kappa
&\defi
\mathcal L_{\beta_0}(Z_0;p_0)
+\frac{9\left(\|p_0\|^2+\|p_*\|^2\right)}{\beta_0}
+(6T+4)\varepsilon_pC_1\tau_*.
\label{eq:xspace_bar_alpha_def}\\
\bar\tau
&\defi
\max\left\{
1,
2\sup\left\{
\operatorname{tr}(Z):
Z\in\mathbb S_+^n,
\inner{C}{Z}+\frac{\beta_0}{4}\|r(Z)\|^2\le\bar\kappa
\right\}
\right\}.\label{eq:xspace_uniform_trace_radius_def}
\end{align}
Then, the X-AL method returns a pair $(\bar Z,\bar p)$ satisfying \eqref{eqn:framework-goal} in at most $T$ outer iterations and at most
\begin{equation}\label{eq:xspace_total_inner_complexity_bound}
2T
+24\bar\tau^2\|\A\|^2\bar\beta(\varepsilon)
\left(
\frac{1}{\varepsilon_c}
+\frac{2}{\varepsilon_d}
\right)
\end{equation}
inner iterations.
\end{theorem}

We now make two remarks about Theorem~\ref{thm:total_AL_X_complexity}.
First, the bound in \eqref{eq:xspace_total_inner_complexity_bound} does not
involve the objective gaps
$\mathcal L_{\beta_{t-1}}(Z_{t-1};p_{t-1})-\min_{X \in \mathbb{S}_+^n}\mathcal{L}_{\beta_{t-1}}(X;p_{t-1})$
of the warm starts $Z_{t-1}$ for the AL subproblems. This is a consequence of
Proposition~\ref{Prop:RUFW}(b), whose bound is independent of the initial
objective gap because the first conic search of every RUFW call contains the
full Frank--Wolfe step to the linear-oracle point over the trace ball
$\Delta_{\tau_t}$ (see Lemma~\ref{lem:RUFW_gap_complexity}(a)). The initial
pair $(Z_0,p_0)$ enters \eqref{eq:xspace_total_inner_complexity_bound} only
through $T$, $\bar\beta(\varepsilon)$, and the radius $\bar\tau$. Second,
suppose that $\varepsilon_p=\varepsilon_c=\varepsilon_d=\varepsilon\in(0,1]$.
Then $C_1=2$, $\bar\beta(\varepsilon)=\mathcal O(1/\varepsilon)$, and
$T=\mathcal O(\log(1/\varepsilon))$. Moreover, since
$\varepsilon T$ is bounded for $\varepsilon\in(0,1]$, the scalar $\bar\kappa$ in
\eqref{eq:xspace_bar_alpha_def} is bounded above by a constant that does not
depend on $\varepsilon$, and hence so is $\bar\tau$ in view of
\eqref{eq:xspace_uniform_trace_radius_def} and
Proposition~\ref{prop:AL-sublevel-set-bounded}. Thus, for fixed problem data,
initial pair, and $\beta_0$, the bound
\eqref{eq:xspace_total_inner_complexity_bound} is
$\mathcal O(\varepsilon^{-2})$.

\subsection{Proof of Theorem~\ref{thm:total_AL_X_complexity}}

The proof of Theorem~\ref{thm:total_AL_X_complexity} relies on the following lemma.

\begin{lemma}\label{prop:xspace_one_call}
For every iteration index $t$ generated by the $X$-space AL method, define
\begin{equation}\label{eq:xspace_tau_t_def}
\tau_t
\ \defi\
\max\left\{1,\ 2\sup\left\{\operatorname{tr}(X):
X\in\mathcal S\!\left(
Z_{t-1};\mathcal L_{\beta_{t-1}}(\cdot;p_{t-1})
\right)\right\}\right\}.
\end{equation}
Then
the following statements about the $t$-th outer iteration hold:
\begin{itemize}
\item[(a)] the quantity $\tau_t$ is finite;

\item[(b)] the RUFW call in Step~{\bf 1} terminates in at most
\begin{equation}\label{eq:xspace_one_call_bound}
\left\lfloor
2+12\tau_t^2\beta_{t-1}\|\A\|^2
\left(
\frac{1}{\varepsilon_c}
+\frac{2}{\varepsilon_d}
\right)
\right\rfloor
\end{equation}
inner iterations with a pair $(Z_t,p_t)$ satisfying \eqref{eq:inner_comp} and the relation
\begin{equation}\label{eq:xspace_inner_descent}
\mathcal L_{\beta_{t-1}}(Z_t;p_{t-1})
\ \le\
\mathcal L_{\beta_{t-1}}(Z_{t-1};p_{t-1}).
\end{equation}
\end{itemize}
\end{lemma}

\begin{proof}
\textit{(a)} Fix a generated iteration index $t$. Definition
\eqref{SublevelDef} implies
\[
\mathcal S\!\left(
Z_{t-1};\mathcal L_{\beta_{t-1}}(\cdot;p_{t-1})
\right)
=
\left\{X\in\mathbb S_+^n:
\mathcal L_{\beta_{t-1}}(X;p_{t-1})
\le \mathcal L_{\beta_{t-1}}(Z_{t-1};p_{t-1})
\right\}.
\]
Proposition~\ref{prop:AL-sublevel-set-bounded}, applied with
$(q,\beta,\kappa)=(p_{t-1},\beta_{t-1},
\mathcal L_{\beta_{t-1}}(Z_{t-1};p_{t-1}))$, shows that this set is
bounded. It is nonempty because it contains
$Z_{t-1}$. Hence, $\tau_t<\infty$ by
\eqref{eq:xspace_tau_t_def}.

\textit{(b)} Relations \eqref{eqn:AL_general} and
\eqref{eqn:AL_gradient} verify RUFW assumption~\textbf{(A1)} with
$L=\beta_{t-1}\|\A\|^2$, while part~\textit{(a)} verifies
assumption~\textbf{(A2)}. Apply Proposition~\ref{Prop:RUFW} to the
objective $\mathcal L_{\beta_{t-1}}(\cdot;p_{t-1})$ with initial point
$Z_{t-1}$, Lipschitz constant $\beta_{t-1}\|\A\|^2$, and tolerance pair
$(\epsilon,\delta)=(\varepsilon_d,\varepsilon_c)$. Indeed, it follows from
\eqref{eq:RUFW_level_set}, \eqref{SublevelDef}, and
\eqref{eq:xspace_tau_t_def} that $\bar\tau=\tau_t$, and hence the bound in
\eqref{eq:RUFW_two_tolerance_complexity} becomes the bound in
\eqref{eq:xspace_one_call_bound}.
It then follows from these conclusions, Proposition~\ref{Prop:RUFW}(a)-(b), the stopping criteria in \eqref{eq:xspace_RUFW_stop}, and relation \eqref{eq:xspace_gradient_identity} that RUFW terminates with a pair $(Z_t,p_t)$ satisfying \eqref{eq:inner_comp} and \eqref{eq:xspace_inner_descent} in at most \eqref{eq:xspace_one_call_bound} inner iterations.
\end{proof}

The remainder of this subsection is dedicated to bounding the local radii $\tau_t$ uniformly. For every $\kappa\in\R$, define
\begin{equation}\label{eq:xspace_fixed_penalized_level_set_def}
\cP_{\beta_0}(\kappa)
\ \defi\
\left\{U\in\mathbb S_+^n:
\inner{C}{U}+\frac{\beta_0}{4}\|\A U-b\|^2\le\kappa
\right\}.
\end{equation}
We first observe that Proposition~\ref{prop:AL-sublevel-set-bounded},
applied with $(q,\beta)=(0,\beta_0/2)$, implies that
$\cP_{\beta_0}(\kappa)$ is bounded for every $\kappa\in\R$.
Our strategy for showing that the local sublevel sets are uniformly bounded
is to find $\bar\kappa\in\R$, independent of $t$, such that
\[
\mathcal S\!\left(
Z_{t-1};\mathcal L_{\beta_{t-1}}(\cdot;p_{t-1})
\right)
\subseteq
\cP_{\beta_0}(\bar\kappa)
\]
for every generated iteration index $t$.

\begin{lemma}\label{lem:xspace_combined_multiplier_dynamic_AL_bounds}
For every iteration index $t \ge 1$, we have:
\begin{align}
\label{eq:xspace_uniform_multiplier_bound}
\frac{\|p_t\|^2}{\beta_t}
&\le
\frac{3}{\beta_0}\left(\|p_0\|^2+\|p_*\|^2\right)
+4\varepsilon_pC_1\tau_*,\\
\mathcal L_{\beta_t}(Z_t;p_t)
&\le
\mathcal L_{\beta_0}(Z_0;p_0)
+\frac{3\|p_0-p_*\|^2}{\beta_0}
+6t\varepsilon_pC_1\tau_*.\label{eq:xspace_dynamic_AL_bound}
\end{align}
where $C_1$ and $\tau_*$ are defined in
\eqref{eq:bar_beta_tolerance_bound_corrected}.
\end{lemma}

\begin{proof}
The triangle inequality, \eqref{eq:key_iterate_bound_lp_corrected}, and
$\beta_t=2\beta_{t-1}$ imply
\begin{align*}
\frac{\|p_t\|^2}{\beta_t}
&\le \frac{\|p_*\|^2+\|p_t-p_*\|^2}{\beta_{t-1}}\\
&\overset{\eqref{eq:key_iterate_bound_lp_corrected}}{\le}
\frac{\|p_*\|^2+\|p_0-p_*\|^2}{\beta_{t-1}}
+4\varepsilon_pC_1\tau_*\\
&\le \frac{3(\|p_0\|^2+\|p_*\|^2)}{\beta_0}
+4\varepsilon_pC_1\tau_*,
\end{align*}
where the last inequality uses $\beta_{t-1}\ge\beta_0$ and
$\|p_0-p_*\|^2\le2(\|p_0\|^2+\|p_*\|^2)$.
This proves \eqref{eq:xspace_uniform_multiplier_bound}.

Relations \eqref{eq:xspace_inner_descent}, \eqref{eqn:AL_general},
\eqref{eq:xspace_multiplier_update}, and \eqref{eq:xspace_penalty_update}
also imply
\begin{align}
\mathcal L_{\beta_t}(Z_t;p_t)
&\le \mathcal L_{\beta_{t-1}}(Z_{t-1};p_{t-1})
+\frac{3\beta_{t-1}}2\|r(Z_t)\|^2\notag\\
&\overset{\eqref{eq:key_iterate_bound_lp_corrected}}{\le}
\mathcal L_{\beta_{t-1}}(Z_{t-1};p_{t-1})
+\frac{3\|p_0-p_*\|^2}{2\beta_{t-1}}
+6\varepsilon_pC_1\tau_*.
\label{eq:xspace_dynamic_AL_step}
\end{align}
Induction using \eqref{eq:xspace_dynamic_AL_step} and
\eqref{eqn:beta_identities_corrected} yields
\[
\mathcal L_{\beta_t}(Z_t;p_t)
\le \mathcal L_{\beta_0}(Z_0;p_0)
+\frac{3\|p_0-p_*\|^2}{\beta_0}(1-2^{-t})
+6t\varepsilon_pC_1\tau_*,
\]
which implies \eqref{eq:xspace_dynamic_AL_bound}.
\end{proof}

\begin{lemma}\label{lem:xspace_local_fixed_level}
For every $t\ge0$, we have
\begin{equation}\label{eq:xspace_local_fixed_level}
\mathcal S\!\left(Z_t;\mathcal L_{\beta_t}(\cdot;p_t)\right)
\subseteq
\cP_{\beta_0}(\bar\kappa),
\end{equation}
where $\mathcal S(\cdot;\cdot)$ is as in \eqref{SublevelDef}, and
$\bar\kappa$ is as in \eqref{eq:xspace_bar_alpha_def}.
\end{lemma}

\begin{proof}
Assume that $U\in\mathcal S\!\left(Z_t;\mathcal L_{\beta_t}(\cdot;p_t)\right)$, which clearly implies that $\mathcal L_{\beta_t}(U;p_t)\le \mathcal L_{\beta_t}(Z_t;p_t)$, in view of the definition of $\mathcal{S}(\cdot;\cdot)$ in  \eqref{SublevelDef}.
This observation together with \eqref{eqn:AL_general}, the fact that $\beta_t\ge\beta_0$, and Young's inequality, then imply that
\begin{align}\label{eqn:identity2}
\inner{C}{U}+\frac{\beta_0}{4}\|r(U)\|^2 & \le \inner{C}{U}+\frac{\beta_t}{4}\|r(U)\|^2
\overset{\eqref{eqn:AL_general}}{=}
\mathcal L_{\beta_t}(U;p_t)-\inner{p_t}{r(U)}-\frac{\beta_t}{4}\|r(U)\|^2\notag\\
&\le
\mathcal L_{\beta_t}(Z_t;p_t)-\inner{p_t}{r(U)}-\frac{\beta_t}{4}\|r(U)\|^2
\le
\mathcal L_{\beta_t}(Z_t;p_t)
+\frac{\|p_t\|^2}{\beta_t}\le\bar\kappa,
\end{align}
where the last inequality is due to \eqref{eq:xspace_uniform_multiplier_bound}, \eqref{eq:xspace_dynamic_AL_bound}, and the definition of $\bar\kappa$ in \eqref{eq:xspace_bar_alpha_def}.
The conclusion of the lemma now follows from
\eqref{eq:xspace_fixed_penalized_level_set_def}.
\end{proof}

\begin{proof}[Proof of Theorem~\ref{thm:total_AL_X_complexity}]
Lemma~\ref{prop:xspace_one_call}(b) implies that X-AL is an instance of the
AL framework. Hence, Theorem~\ref{thm:complexity_relaxed_corrected} implies
that it terminates at an outer iteration index $\bar t\le T$ with a pair
satisfying \eqref{eqn:framework-goal}.
For every $1\le t\le\bar t$, Lemma~\ref{lem:xspace_local_fixed_level} at
index $t-1$ and definitions \eqref{eq:xspace_tau_t_def} and
\eqref{eq:xspace_uniform_trace_radius_def} imply $\tau_t\le\bar\tau$.
Thus, induction using \eqref{eq:xspace_one_call_bound} and
$\beta_{t-1}=2\beta_{t-2}$ shows that the first $t$ outer iterations
perform at most
\begin{equation}\label{eq:xspace_prefix_inner_bound}
2t+12\bar\tau^2\|\A\|^2
\left(\frac1{\varepsilon_c}+\frac2{\varepsilon_d}\right)
(2\beta_{t-1}-\beta_0)
\end{equation}
inner iterations. Applying \eqref{eq:xspace_prefix_inner_bound} with
$t=\bar t$ and using $\bar t\le T$ and
\eqref{eq:beta_last_bound_corrected} yields
\eqref{eq:xspace_total_inner_complexity_bound}.
\end{proof}

The following lemma bounds the total decrease of the AL objectives over
the inner solves. It will be used to bound the ADAP-AIPP work of HALLaR+.

\begin{lemma}\label{lem:xspace_total_AL_decrease}
Suppose that the iterates generated by the AL framework satisfy
\eqref{eq:xspace_inner_descent} for every $t\ge1$.
Then the total decrease of the AL objectives over all inner solves is at most
\begin{equation}\label{eq:xspace_total_AL_decrease_bound}
\bar\kappa+\|b\|\|p_*\|,
\end{equation}
where $\bar\kappa$ is as in \eqref{eq:xspace_bar_alpha_def}.
\end{lemma}
\begin{proof}
Let $t\ge1$ be an iteration index.
The multiplier update, penalty doubling, and \eqref{eqn:AL_general} imply
\begin{align}
\mathcal L_{\beta_t}(Z_t;p_t)
+\frac{3\|p_t-p_*\|^2}{2\beta_t}
&=\mathcal L_{\beta_{t-1}}(Z_t;p_{t-1})
+\frac{3\|p_{t-1}-p_*\|^2}{2\beta_{t-1}}\notag\\
&\quad+3\inner{p_t-p_*}{r(Z_t)}
-\frac{3\|p_t-p_*\|^2}{4\beta_{t-1}}\notag\\
&\le\mathcal L_{\beta_{t-1}}(Z_t;p_{t-1})
+\frac{3\|p_{t-1}-p_*\|^2}{2\beta_{t-1}}
+3\varepsilon_pC_1\tau_*.
\label{eq:xspace_AL_decrease_step}
\end{align}
The last inequality follows from \eqref{eq:ip-resid-ial},
\eqref{eq:fw_gap_lp_translation}, and
\eqref{eq:bar_beta_tolerance_bound_corrected}.
Moreover, \eqref{eq:KKT} and completion of the square imply
\begin{equation}\label{eq:xspace_AL_decrease_lower}
\mathcal L_{\beta_t}(Z_t;p_t)
+\frac{3\|p_t-p_*\|^2}{2\beta_t}
\ge \inner{b}{p_*}+\frac{\|p_t-p_*\|^2}{\beta_t}
\ge \inner{b}{p_*}.
\end{equation}
Applying \eqref{eq:xspace_AL_decrease_step} successively over the outer
iterations and using \eqref{eq:xspace_AL_decrease_lower} and $\bar t\le T$
bounds the total decrease over the inner solves by
\begin{align*}
\mathcal L_{\beta_0}(Z_0;p_0)
+\frac{3\|p_0-p_*\|^2}{2\beta_0}
+3T\varepsilon_pC_1\tau_*-\inner{b}{p_*}\overset{\eqref{eq:xspace_bar_alpha_def}}{\le}
\bar\kappa+\|b\|\|p_*\|,
\end{align*}
where the last inequality uses
$\|p_0-p_*\|^2\le2(\|p_0\|^2+\|p_*\|^2)$ and Cauchy--Schwarz.
\end{proof}

\section{ HALLaR+ in the \texorpdfstring{$U$}{U}-Space}
\label{sec:u-space-hallar}

This section presents HALLaR+, the factorized implementation of the $X$-AL
method of Section~\ref{subsec:xspace-step-one}. It contains two subsections. Subsection~\ref{subsec:u-space-hallar-description} states the HALLaR+ method and presents its main complexity result, namely Theorem~\ref{thm:GeneralHLRComplexity}. Subsection~\ref{ProofUAL} is devoted to the proof of Theorem~\ref{thm:GeneralHLRComplexity}.

 \subsection{Description of HALLaR+ and Main Theorem}
 \label{subsec:u-space-hallar-description}

 We start by stating the HALLaR+ method.

\vspace{0.75em}
\noindent
\begin{minipage}[t]{1\columnwidth}
\rule[0.5ex]{1\columnwidth}{1pt}

\noindent \textbf{HALLaR+ Method}

\noindent \rule[0.5ex]{1\columnwidth}{1pt}
\end{minipage}

\noindent \textbf{Input}: Initial pair
$(Y_0,p_0)\in\R^{n\times r_0}\times\R^m$ for some $r_0\ge1$,
tolerance triple
$(\varepsilon_p,\varepsilon_c,\varepsilon_d)\in\R_{++}^3$,
initial penalty parameter $\beta_0\in\R_{++}$, and hybrid low-rank (HLR) parameters
$\eta_0>0$, $M_0>1/2$, and $\bar\rho>0$.

\noindent \textbf{Output}: Pair $(\bar Y,\bar p)$ such that
$(\bar Y\bar Y^\top,\bar p)$ satisfies \eqref{eqn:framework-goal}.

\begin{itemize}[leftmargin=2.0em]
\item[{\bf 0.}] Set $t\gets1$.
\item[{\bf 1.}] Perform the following $U$-space HLR (U-HLR) call to approximately minimize $t$-th AL subproblem:
\begin{enumerate}[leftmargin=3.2em]
\item[{\bf 1.a)}] Set $(\tilde U_0,\hat r_1)\gets (Y_{t-1},r_{t-1})$
and $k\gets 1$.
\item[{\bf 1.b)}] Call ADAP-AIPP in Appendix~\ref{app:aipp} with input 
\begin{equation}\label{AIPPCallHLR}
    W_0=\tilde U_{k-1}, \quad  g(U)=\mathcal L_{\beta_{t-1}}(UU^{\top};p_{t-1}), \quad \Omega=\mathbb R^{n\times \hat r_k}, \quad \eta_0=\eta_0, \quad M_0=M_0, \quad \bar \rho=\bar \rho
\end{equation}
to obtain output pair $(U_k,R_k)\in \mathbb R^{n\times \hat r_{k}}\times \mathbb R^{n\times \hat r_k}$ and set
\[S_k\gets C-\mathcal A^{*}[p_{t-1}+\beta_{t-1}r(U_kU_k^{\top})].\]
\item[{\bf 1.c)}] 
Compute a minimum eigenpair $(\lambda_{\min}(S_k),v_{\min}(S_k))$ of $S_k$ with $\|v_{\min}(S_k)\|=1$, and set
\begin{equation}\label{eigenv}
y^{F}_k\ \defi\ \begin{cases} v_{\min}(S_k), & \lambda_{\min}(S_k)<0,\\ 0, & \text{otherwise}.\end{cases}
\end{equation}

\item[{\bf 1.d)}] \textbf{If}
\begin{equation}\label{eq:hlr_stop_test}
\lambda_{\min}(S_k)\ \ge\ -\varepsilon_d,\qquad \left|\inner{S_kU_k}{U_k}\right|\ \le\ \varepsilon_c,
\end{equation}
then \textbf{set $(Y_t,r_t)\gets(U_k, \hat r_{k})$} and go to {\bf Step 2}.

\item[{\bf 1.e)}] Compute
\begin{equation}\label{eq:Uspace_RUFW_linesearch}
(\alpha_k,\omega_k)
\in
\argmin\left\{
\mathcal L_{\beta_{t-1}}\left(\alpha U_kU_k^{\top}+\omega y^{F}_k{(y_k^{F})}^\top;p_{t-1}\right):
\alpha\ge0,\ \omega\ge0
\right\}
\end{equation}
and set
\begin{equation}\label{TildeY}
(\tilde U_k,\hat r_{k+1})\ \defi\ \begin{cases}
(\sqrt{\omega_k}\,y_k^{F},\ 1), & \alpha_k=0,\\
\left(\left[\sqrt{\alpha_k}\,U_k,\ \sqrt{\omega_k}\,y^{F}_k\right],\ \hat r_{k}+1\right), & \alpha_k>0.
\end{cases}
\end{equation}
\item[{\bf 1.f)}] Set $k \gets k+1$
and return to Step~{\bf 1.b)}.
\end{enumerate}

\item[{\bf 2.}] Set
\begin{equation}\label{eq:Uspace_multiplier_update}
p_t\gets p_{t-1}+\beta_{t-1}r(Y_tY_t^{\top}).
\end{equation}
\item[{\bf 3.}] \textbf{If} $\|r(Y_tY_t^{\top})\|\le\varepsilon_p$, then \textbf{output} 
$(\bar Y,\bar p)\gets(Y_t,p_t)$ and \textbf{stop}.
\item[{\bf 4.}] Set
\begin{equation}\label{eq:Uspace_penalty_update}
\beta_t\gets2\beta_{t-1},
\qquad
t\gets t+1,
\end{equation}
and return to Step~{\bf 1}.
\end{itemize}
\noindent\rule[0.5ex]{1\columnwidth}{1pt}
\vspace{0.75em}

\begin{theorem}\label{thm:GeneralHLRComplexity}
Set $\varepsilon\defi(\varepsilon_p,\varepsilon_c,\varepsilon_d)$. Let
$C_1$, $\tau_*$, and $\bar\beta(\varepsilon)$ be as in
\eqref{eq:bar_beta_tolerance_bound_corrected}, let $T$ be as in
\eqref{eq:T_def_corrected}, and let $\bar\kappa$ and $\bar\tau$ be as in
\eqref{eq:xspace_bar_alpha_def} and \eqref{eq:xspace_uniform_trace_radius_def},
respectively. Define
\begin{equation}\label{barCeps}
    \bar C(\varepsilon)\defi24\bar\beta(\varepsilon)\|\mathcal A\|^2
(3\sqrt{2\bar\tau}+1)^2+8\left[\|C\|+\|\mathcal A\|
\left(\|p_*\|
+\sqrt{\|p_0-p_*\|^2+4\varepsilon_pC_1\bar\beta(\varepsilon)\tau_*}
+\bar\beta(\varepsilon)\|b\|\right)\right].
\end{equation}
The following statements about HALLaR+ hold:
\begin{itemize}
    \item[(a)] it returns a pair $(\bar Y,\bar p)$ such that $(\bar Y\bar Y^{\top},\bar p)$ is an $(\varepsilon_p,\varepsilon_c,\varepsilon_d)$-solution of \eqref{eqn:general_probl}--\eqref{eqn:lagrangian_dual} in at most $T$ outer iterations and at most
\begin{equation}\label{eq:Uspace_total_inner_complexity_bound}
2T
+24\bar\tau^2\|\A\|^2\bar\beta(\varepsilon)
\left(
\frac{1}{\varepsilon_c}
+\frac{2}{\varepsilon_d}
\right)
\end{equation}
inner iterations/minimum eigenvector (MEV) computations;
\item[(b)] the total number of ADAP-AIPP iterations/ADAP-FISTA calls it performs in its Step~{\bf 1.b)} is at most
\begin{equation}\label{eq:Uspace_total_aipp_bound}
\left(1+\left\lceil\log_2^+\left(\eta_0\bar C(\varepsilon)\right)\right\rceil\right)
\left[
2T
+24\bar\tau^2\|\A\|^2\bar\beta(\varepsilon)
\left(
\frac{1}{\varepsilon_c}
+\frac{2}{\varepsilon_d}
\right)
\right]
+\frac{C_{\sigma}\max\left\{\eta_0^{-1},\, \bar C(\varepsilon)\right\}}{\bar\rho^2}\left(\bar\kappa+\|b\|\|p_*\|\right),
\end{equation}
where $\eta_0$ and $\bar\rho$ are inputs to HALLaR+ and $C_\sigma$ is as in \eqref{Csigma}.
\end{itemize}
\end{theorem}

\subsection{Proof of Theorem~\ref{thm:GeneralHLRComplexity}}\label{ProofUAL}
We first state a key proposition, without proof, which will be instrumental in proving Theorem~\ref{thm:GeneralHLRComplexity}. Due to the length of its proof, we defer the proof to Subsection~\ref{prop:hlrProof}.
\begin{proposition}\label{prop:hlr-general}
Consider the $t$-th iteration of HALLaR+ and let
\begin{equation}\label{Cbartau}
    C(\tau_t,\beta_{t-1})\defi24\beta_{t-1}\|\mathcal A\|^2(3\sqrt{2\tau_t}+1)^2+8\|C-\mathcal A^{*}(p_{t-1}+\beta_{t-1}b)\|
\end{equation}
where $\tau_t$ is as in \eqref{eq:xspace_tau_t_def} with
$Z_{t-1}\defi Y_{t-1}Y_{t-1}^\top$.
The following statements about the $U$-HLR call in Step~{\bf 1} hold:
\begin{itemize}
    \item[(a)] it performs at most
    \begin{equation}\label{GeneralHLRComplexity}
    \left\lfloor
    2+12\tau_t^2\beta_{t-1}\|\mathcal A\|^2
    \left(\frac{1}{\varepsilon_c}+\frac{2}{\varepsilon_d}\right)
    \right\rfloor
    \end{equation}
iterations/minimum eigenvector (MEV) computations to find a pair $(Y_t,r_t)$ satisfying
\begin{align}
\lambda_{\min}\left(\nabla \mathcal L_{\beta_{t-1}}(Y_tY_t^{\top};p_{t-1})\right)&\ge -\varepsilon_d,
\qquad
\left|\left\langle \nabla \mathcal L_{\beta_{t-1}}(Y_tY_t^{\top};p_{t-1})Y_t,Y_t\right\rangle\right|
\le \varepsilon_c,\label{SolutionTypeHLR}\\
\mathcal L_{\beta_{t-1}}(Y_tY_t^{\top};p_{t-1})&\leq \mathcal L_{\beta_{t-1}}(Y_{t-1}Y_{t-1}^{\top};p_{t-1}),\label{descentHLR}
\end{align}
where $(\varepsilon_c,\varepsilon_d)$ are input tolerances;
\item[(b)] the total number of ADAP-AIPP iterations/ADAP-FISTA calls it performs is at most
\begin{align*}
&\left(1+\left\lceil\log_2^+\left(\eta_0C(\tau_t,\beta_{t-1})\right)\right\rceil\right)
\left(
2+12\tau_t^2\beta_{t-1}\|\mathcal A\|^2
\left(\frac{1}{\varepsilon_c}+\frac{2}{\varepsilon_d}\right)
\right)\\
&\qquad+\frac{C_{\sigma}\max\left\{\eta_0^{-1},\, C(\tau_t,\beta_{t-1})\right\}}{\bar\rho^2}\left[\mathcal L_{\beta_{t-1}}(Y_{t-1}Y_{t-1}^{\top};p_{t-1})-\mathcal L_{\beta_{t-1}}(Y_t Y_t^{\top};p_{t-1})\right]
\end{align*}
where $\eta_0$ and $\bar \rho$ are inputs to HALLaR+ and $C_{\sigma}$ and $C(\tau_t,\beta_{t-1})$ are as in \eqref{Csigma} and \eqref{Cbartau}, respectively.
\end{itemize}
\end{proposition}

We are now ready to prove Theorem~\ref{thm:GeneralHLRComplexity}.
\begin{proof}[Proof of Theorem~\ref{thm:GeneralHLRComplexity}]
For every generated iteration index $t\ge0$, set $Z_t\defi Y_tY_t^\top$.

(a) Let $t\ge1$ be a generated iteration index. Proposition~\ref{prop:hlr-general}(a) implies that the $U$-HLR call in Step~{\bf 1} terminates with a matrix $Y_t$ satisfying \eqref{SolutionTypeHLR} and \eqref{descentHLR}. Moreover, relations \eqref{eqn:AL_gradient} and \eqref{eq:Uspace_multiplier_update} imply that
\[
\nabla_X \mathcal L_{\beta_{t-1}}(Z_t;p_{t-1})
=C-\A^*\left[p_{t-1}+\beta_{t-1}r(Z_t)\right]
=C-\A^*p_t
=S_t,
\]
and the identity $Z_t=Y_tY_t^\top$ implies that $\inner{S_tY_t}{Y_t}=\inner{S_t}{Z_t}$. Hence, \eqref{SolutionTypeHLR} implies that $(Z_t,p_t)$ satisfies \eqref{eq:inner_comp}, and \eqref{descentHLR} is precisely \eqref{eq:xspace_inner_descent}. Thus, Step~{\bf 1} of HALLaR+ implements Step~{\bf 1} of the AL framework, and Steps~{\bf 2}--{\bf 4} of HALLaR+ coincide with Steps~{\bf 2}--{\bf 4} of the AL framework with $Z_t=Y_tY_t^\top$. It then follows from Theorem~\ref{thm:complexity_relaxed_corrected} that HALLaR+ terminates at an outer iteration index $\bar t\le T$ with a pair $(\bar Y,\bar p)$ such that $(\bar Y\bar Y^\top,\bar p)$ is an $(\varepsilon_p,\varepsilon_c,\varepsilon_d)$-solution of \eqref{eqn:general_probl}--\eqref{eqn:lagrangian_dual}, that $\bar\beta=\beta_{\bar t-1}\le\bar\beta(\varepsilon)$, and that every multiplier $p_t$ with $t\ge1$ satisfies \eqref{eq:multiplier_bound_corrected}.

The proofs of Lemmas~\ref{lem:xspace_combined_multiplier_dynamic_AL_bounds}
and \ref{lem:xspace_local_fixed_level} use only \eqref{eq:inner_comp},
\eqref{eq:xspace_inner_descent}, and the outer updates. Hence, they apply
with $Z_t=Y_tY_t^\top$ and imply $\tau_t\le\bar\tau$.
Applying the induction used in \eqref{eq:xspace_prefix_inner_bound} to
\eqref{GeneralHLRComplexity} and using $\bar t\le T$ and
\eqref{eq:beta_last_bound_corrected} yields
\eqref{eq:Uspace_total_inner_complexity_bound}. Since every inner iteration
performs one MEV computation, this proves (a).

(b) For every $1\le t\le\bar t$, the triangle inequality,
\eqref{eq:multiplier_bound_corrected}, and
\eqref{eq:beta_last_bound_corrected} imply
\begin{align*}
\|C-\A^*(p_{t-1}+\beta_{t-1}b)\|
&\le\|C\|+\|\A\|
\left(\|p_{t-1}\|+\beta_{t-1}\|b\|\right)\\
&\le\|C\|+\|\A\|\left(\|p_*\|
+\sqrt{\|p_0-p_*\|^2+4\varepsilon_pC_1\bar\beta(\varepsilon)\tau_*}
+\bar\beta(\varepsilon)\|b\|\right),
\end{align*}
where, for $t=1$, the last inequality uses
$\|p_0\|\le\|p_*\|+\|p_0-p_*\|$.
Together with $\tau_t\le\bar\tau$ and definitions
\eqref{barCeps} and \eqref{Cbartau}, this implies
\begin{equation}\label{eq:hlr_uniform_C_bound}
C(\tau_t,\beta_{t-1})\le\bar C(\varepsilon).
\end{equation}
Lemma~\ref{lem:xspace_total_AL_decrease}, applied with $Z_t=Y_tY_t^\top$,
bounds the total AL decrease by \eqref{eq:xspace_total_AL_decrease_bound}.
Combining this bound with Proposition~\ref{prop:hlr-general}(b),
\eqref{eq:hlr_uniform_C_bound}, and part~(a) yields
\eqref{eq:Uspace_total_aipp_bound}.
\end{proof}

\subsubsection{Proof of Proposition~\ref{prop:hlr-general}}\label{prop:hlrProof}

This subsection is devoted to proving Proposition~\ref{prop:hlr-general}, which establishes the iteration complexity of the U-HLR call performed in Step~\textbf{1} of the $t$-th iteration of HALLaR+.
Accordingly, throughout this subsection, we fix an arbitrary iteration index $t\geq 1$ and consider just the $t$-th iteration of HALLaR+.
We now define several quantities which arise in the technical complexity analysis of the U-HLR call.
Let
\begin{equation}\label{DefPhiT}
\phi_t(\cdot)\defi\mathcal L_{\beta_{t-1}}(\cdot; p_{t-1}), \quad g_t(U)\defi\phi_t(UU^{\top}).
\end{equation}
For every iteration index $k\geq 1$ of the U-HLR call in Step~1, define
\begin{equation}\label{DefsHLR}
\tD{k}\defi6D(\tilde U_{k-1};g_t)+1, \quad \tLip{k}\defi[{\rm LIP}](\tD{k};g_t),
\end{equation}
\begin{equation}\label{def:tildeM}
 \tM{k}\defi24\beta_{t-1}\|\mathcal A\|^2\tD{k}^2+8\|C-\mathcal A^{*}(p_{t-1}+\beta_{t-1}b)\|
 \end{equation}
where $D(\cdot; \cdot)$ is as in \eqref{DiameterDef}, and $[{\rm LIP}](\cdot;\cdot)$ is as in \eqref{SmoothnessDef}.

We now start off with a lemma that shows that the call to ADAP-AIPP in Step~\textbf{1.b)} is made with inputs that satisfy assumptions \textbf{(C1)}, $\textbf{(C2)}$, and $\textbf{(C4)}$ in Appendix~\ref{app:aipp}. It is easy to see that assumption $\textbf{(C3)}$ is always satisfied.

\begin{lemma}\label{AIPPHLR}
   Let $k\geq 1$ and consider the $k$-th iteration of the U-HLR call in Step~\textbf{1}. The input of the ADAP-AIPP call made in Step~\textbf{1.b)} of U-HLR satisfies
    \begin{itemize}
        \item[(a)] the function $g_t(U)=\phi_t(UU^{\top})$ is differentiable with locally Lipschitz continuous gradient. Moreover,
        \begin{equation}\label{eq:tildeLipBound}
           \tLip{k}\leq 6\beta_{t-1}\|\mathcal A\|^2\tD{k}^2+2\|C-\mathcal A^{*}(p_{t-1}+\beta_{t-1}b)\|;
        \end{equation}
        \item[(b)] $g_t(U)$ is bounded below on $\mathbb R^{n\times \hat r_k}$;
        \item[(c)] the sublevel set $\mathcal S(\tilde U_{k-1};g_t)$ has finite diameter $D(\tilde U_{k-1};g_t)$,
        where $\mathcal S(\cdot;\cdot)$ and $D(\cdot;\cdot)$
        are as in \eqref{SublevelDef} and \eqref{DiameterDef}, respectively.
    \end{itemize}
\end{lemma}
\begin{proof}
    (a) It is easy to see that $g_t(U)=\phi_t(UU^{\top})$ is differentiable. To show $g_t(U)$ has locally Lipschitz continuous gradient, let $R>0$ be an arbitrary scalar and $B[R]$ be a ball centered at $0$ with radius $R$. It suffices to show that $g_t(U)$ has Lipschitz continuous gradient on $B[R]$. Thus, let $U$ and $V$ be in $B[R]$.
    First, it is easy to see that
   \begin{align}
       \|UU^{\top}-VV^{\top}\|&=\|U(U-V)^{\top}-(V-U)V^{\top}\| \nonumber\\
       &\leq \|U\|\|U-V\|+\|V\|\|U-V\|\leq 2R\|U-V\|\label{UUBound}.
   \end{align}
Moreover, it follows from the fact that $\phi_t(\cdot)$ is $(\beta_{t-1}\|\mathcal A\|^2)$-smooth that
\begin{align}
    \|\nabla \phi_t(VV^{\top})\|=\|\nabla \phi_t(VV^{\top})-\nabla \phi_t(0)+\nabla \phi_t(0)\|&\leq \beta_{t-1}\|\mathcal A\|^2\|VV^{\top}\|+\|\nabla \phi_t(0)\| \nonumber\\
    &\leq \beta_{t-1}\|\mathcal A\|^2R^2+\|\nabla \phi_t(0)\| \label{UpBound}.
\end{align}
    It then follows from relations \eqref{UUBound} and  \eqref{UpBound}, the fact that $\phi_t$ is $(\beta_{t-1}\|\mathcal A\|^2)$-smooth, and the fact that $\nabla g_t(U)=2\nabla \phi_t(UU^{\top})U$ that
\begin{align}
    &\|\nabla g_t(U)-\nabla g_t(V)\|
    \leq \|2\nabla \phi_t(UU^{\top})U-2\nabla \phi_t(VV^{\top})U
    +2\nabla \phi_t(VV^{\top})U-2\nabla\phi_t(VV^{\top})V\|\nonumber\\
    &\leq 2\|U\|\|\nabla\phi_t(UU^{\top})-\nabla\phi_t(VV^{\top})\|+2\|\nabla\phi_t(VV^{\top})\|\|U-V\|\nonumber\\
    &\leq 2\beta_{t-1}\|\mathcal A\|^2\|U\|\|UU^{\top}-VV^{\top}\|+2\|\nabla\phi_t(VV^{\top})\|\|U-V\|\nonumber\\
    &\overset{\eqref{UUBound},\eqref{UpBound}}{\leq} 4\beta_{t-1}\|\mathcal A\|^2R^2\|U-V\|+2(\beta_{t-1}\|\mathcal A\|^2 R^2+\|\nabla \phi_t(0)\|)\|U-V\|\nonumber \\
    &=(6\beta_{t-1}\|\mathcal A\|^2R^2+2\|\nabla \phi_t(0)\|)\|U-V\|\label{LipBoundg}.
\end{align}
Hence, it follows that $\nabla g_t(U)$ is ($6\beta_{t-1}\|\mathcal A\|^2R^{2}+2\|\nabla \phi_t(0)\|$)-Lipschitz continuous on $B[R]$. Since $R$ was arbitrary, it follows that $\nabla g_t(U)$ is locally Lipschitz continuous. Moreover, it follows from the definitions of $\tLip{k}$ and $[{\rm LIP}](\cdot;\cdot)$ in \eqref{DefsHLR} and \eqref{SmoothnessDef}, respectively, and replacing $R$ in \eqref{LipBoundg} with $\tD{k}$ that
$\tLip{k}\leq 6\beta_{t-1}\|\mathcal A\|^2\tD{k}^2+2\|\nabla \phi_t(0)\|$. The bound in \eqref{eq:tildeLipBound} then follows from this conclusion and the fact that the definition of $\phi_t$ in \eqref{DefPhiT} implies that $\nabla \phi_t(0)=C-\mathcal A^{*}(p_{t-1}+\beta_{t-1}b)$.

(b) It is easy to see that
\begin{align}\label{gboundedbelow}
    \inf_{U\in \mathbb R^{n\times \hat r_k}}g_t(U)=\inf_{U\in \mathbb R^{n\times \hat r_k}}\phi_t(UU^{\top})=\inf_{X\succeq0,\,\textrm{rank}(X)\leq \hat r_k}\phi_t(X)\geq \inf_{X\succeq0}\phi_t(X)>-\infty.
\end{align}
Hence, it follows from \eqref{gboundedbelow} that $g_t$ is bounded below on $\mathbb R^{n\times \hat r_k}$.

(c) Let $U\in \mathcal S(\tilde U_{k-1};g_t)$. It follows by the definition of $\mathcal S(\cdot;\cdot)$ in \eqref{SublevelDef} that
$g_t(U)\leq g_t(\tilde U_{k-1})$, which by the definition of $g_t$ in \eqref{DefPhiT} implies that $\phi_t(UU^{\top})\leq \phi_t(\tilde U_{k-1}\tilde U^{\top}_{k-1})$. Thus, $UU^{\top}\in \mathcal S(\tilde U_{k-1}\tilde U_{k-1}^{\top};\phi_t)$.
Now observe that it follows from Proposition~\ref{prop:AL-sublevel-set-bounded} and the definition of $\phi_t$ in \eqref{DefPhiT} that the sublevel sets of $\phi_t$ are bounded. Hence $\mathcal S(\tilde U_{k-1}\tilde U^{\top}_{k-1};\phi_t)$ is bounded. Thus it follows that $\|UU^{\top}\|\leq M$ for some scalar $M>0$. Now observe from the Cauchy-Schwarz inequality that \[\|U\|^2=\inner{U}{U}=\inner{I_n}{UU^{\top}}\leq \sqrt{n}\|UU^{\top}\|\leq \sqrt{n}M.\]
Thus, we have $\|U\|\leq n^{1/4}\sqrt{M}$. Since $U$ was arbitrary, it follows that $\mathcal S(\tilde U_{k-1};g_t)$ is bounded with finite diameter $D(\tilde U_{k-1};g_t)$.
\end{proof}

The following lemma establishes properties of the ADAP-AIPP call made in Step~\textbf{1.b)}.

\begin{lemma}\label{AIPPHLRLemma}
Consider the $k$-th iteration of the U-HLR call in Step~1. The following statements about the ADAP-AIPP call in Step~\textbf{1.b)} of U-HLR hold:
\begin{itemize}
\item[(a)] it terminates with a point $U_k\in \mathbb R^{n\times \hat r_k}$ satisfying
\begin{align}
    &R_k=2\nabla \phi_t(U_kU_k^{\top})U_k, \quad \|R_k\|\leq \bar \rho \label{HLRaippterm1}\\
    & \phi_t(U_kU_k^{\top})\leq \phi_t(\tilde U_{k-1}\tilde U^{\top}_{k-1})\label{HLRaippterm2};
\end{align}

\item[(b)]
its total number of ADAP-AIPP iterations/ADAP-FISTA calls performed within the ADAP-AIPP call is bounded by
\[
    1+\frac{C_{\sigma}\max\left\{\eta_0^{-1},\, \tM{k}\right\}}{\bar\rho^2}\left[\phi_t(\tilde U_{k-1}\tilde U^{\top}_{k-1})-\phi_t(U_kU_k^{\top})\right]+ \left\lceil\log_2^+\left(\eta_0\tM{k}\right)\right\rceil
    \]
where $C_{\sigma}$ is as in \eqref{Csigma}, $\tM{k}$ is as in \eqref{def:tildeM}, and $\eta_0$ and $\bar \rho$ are inputs to ADAP-AIPP.
\end{itemize}
\end{lemma}
\begin{proof}
(a)
Let $k\geq1$ and consider the $k$-th iteration of the U-HLR call performed in Step~\textbf{1}. By
Lemma~\ref{AIPPHLR}(a)--(c), and because Step~{\bf 1.b)} calls ADAP-AIPP with
$\Omega=\mathbb R^{n\times \hat r_k}$, a closed convex set with easy projection, the inputs of this call satisfy assumptions \textbf{(C1)}--\textbf{(C4)}
of Appendix~\ref{app:aipp}.
It then follows from Proposition~\ref{prop:aipp_generic_complexity}(a), the fact that the call to ADAP-AIPP is made with initial point $W_0=\tilde U_{k-1}$ and function $g_t(U)=\phi_t(UU^{\top})$ as in \eqref{DefPhiT}, and the fact that ADAP-AIPP terminates with $(\bar W,\bar R)=(U_k,R_k)$ that the pair $(U_k,R_k)$ satisfies
\begin{align}
    &R_k\in \nabla g_t(U_k)+\mathcal{N}_{\mathbb R^{n\times \hat r_k}}(U_k),\qquad \|R_k\|\le \bar\rho\label{aippHLR1}\\
    & \phi_t(U_kU_k^{\top})\leq \phi_t(\tilde U_{k-1}\tilde U_{k-1}^{\top}).\label{aippHLR2}
\end{align}
Relations \eqref{HLRaippterm1} and \eqref{HLRaippterm2} then follow immediately from \eqref{aippHLR1}, \eqref{aippHLR2}, and the facts that $\nabla g_t(U_k)=2\nabla \phi_t(U_kU_k^{\top})U_k$ and $\mathcal{N}_{\mathbb R^{n\times \hat r_k}}(U_k)=\left\{0\right\}$.

(b)
Let $k\geq 1$ and consider the $k$-th iteration of the U-HLR call performed in Step~\textbf{1} of HALLaR+.
It follows
from Lemma~\ref{AIPPHLR} that the
input of the ADAP-AIPP call
satisfies assumptions \textbf{(C1)}-\textbf{(C4)} in Appendix~\ref{app:aipp}. It then follows from this observation, Proposition~\ref{prop:aipp_generic_complexity}(b), and the definition of $\tLip{k}$ in \eqref{DefsHLR} that the total number of ADAP-AIPP iterations/ADAP-FISTA calls performed within the ADAP-AIPP call is upper bounded by
\[1+\frac{C_{\sigma}\max\left\{\eta_0^{-1},\,4\tLip{k}\right\}}{\bar\rho^2}\left[\phi_t(\tilde U_{k-1}\tilde U^{\top}_{k-1})-\phi_t(U_kU_k^{\top})\right]+ \left\lceil\log_2^+\left(4\eta_0\tLip{k}\right)\right\rceil.\]
The complexity bound in part (b) then follows immediately from this bound and the fact that the bound on $\tLip{k}$ in Lemma~\ref{AIPPHLR}(a)  and the definition of $\tM{k}$ in \eqref{def:tildeM} imply that
$4\tLip{k}\leq \tM{k}$.
\end{proof}

\begin{lemma}\label{HLRInstance}
The following statements about the U-HLR call performed in Step~\textbf{1} hold:
\begin{itemize}
    \item[(a)] it is an instance of the RUFW framework of Appendix~\ref{sec:RUFW} applied to
    \begin{equation}
    \min_{X\in \mathbb S^{n}_{+}} \phi_t(X)
    \end{equation}
with initial point $\tilde X_0=Y_{t-1}Y_{t-1}^\top$ and tolerance pair $(\epsilon,\delta)=(\varepsilon_d,\varepsilon_c)$,
where $\phi_t$ is as in \eqref{DefPhiT};
    \item[(b)] for every iteration index $k\ge1$ at which Step~\textbf{1.e)} is performed, it holds that
\begin{equation}\label{eq:hlr_local_descent}
\phi_t(\tilde U_k\tilde U_k^{\top})\le\phi_t(U_kU_k^{\top})\le\phi_t(\tilde U_{k-1}\tilde U^{\top}_{k-1}).
\end{equation}
Hence, it follows that for every iteration index $k\geq 1$,
$\phi_t(\tilde U_{k-1}\tilde U_{k-1}^{\top})\leq \phi_t(\tilde U_{0}\tilde U^{\top}_{0}).$
\end{itemize}
\end{lemma}
\begin{proof}
(a) Let $k\geq 1$ be a generated iteration index of the U-HLR call and associate with the U-HLR iterates the matrices
\begin{equation}\label{eq:hlr_factor_matrices}
X_k\ \defi\ U_kU_k^{\top},\qquad X_k^{F}\ \defi\ y^{F}_k{(y^{F}_k)}^{\top}.
\end{equation}
Set $\tilde X_0\defi\tilde U_0\tilde U_0^{\top}=Y_{t-1}Y_{t-1}^\top$, where the last identity is due to Step~\textbf{1.a)}, and
whenever Step~\textbf{1.e)} is executed, also define
$\tilde X_k\defi\tilde U_k\tilde U_k^\top$.
Relation~\eqref{HLRaippterm2} gives
$\phi_t(X_k)\le\phi_t(\tilde X_{k-1})$, which is
\eqref{eq:RUFW_correction_step}. Moreover,
$
S_k=C-\A^*(p_{t-1}+\beta_{t-1}r(X_k))
=\nabla_X\phi_t(X_k)
$
and
$
\inner{S_kU_k}{U_k}=\inner{S_k}{X_k}.
$
Thus, Steps~\textbf{1.c)}--\textbf{1.d)} of the U-HLR call coincide with Steps~2 and 3 of the RUFW framework with $G_k=S_k$ and tolerance pair $(\epsilon,\delta)=(\varepsilon_d,\varepsilon_c)$.
Relation \eqref{eq:Uspace_RUFW_linesearch} is the RUFW conic search \eqref{eqn:alpha-linesearch-appendix}, and direct
multiplication in \eqref{TildeY} yields
\[
\tilde X_k=\alpha_kX_k+\omega_kX_k^F,
\]
which is \eqref{eqn:fw-short-step-appendix}. Hence the U-HLR call is an instance of the RUFW framework applied to $\phi_t$ with initial point $\tilde X_0=Y_{t-1}Y_{t-1}^\top$ and tolerance pair $(\epsilon,\delta)=(\varepsilon_d,\varepsilon_c)$.

(b) By part~(a), Lemma~\ref{lem:RUFW_factor_monotonicity} applies to the iterates $X_k=U_kU_k^{\top}$, $X_k^{F}=y^F_k(y^F_k)^{\top}$, and $\tilde X_k=\tilde U_k\tilde U_k^{\top}$. Thus, \eqref{eq:RUFW_factor_monotonicity} gives $\phi_t(\tilde U_k\tilde U_k^{\top})\le\phi_t(U_kU_k^{\top})\le\phi_t(\tilde U_{k-1}\tilde U^{\top}_{k-1})$, which is \eqref{eq:hlr_local_descent}. Applying these relations recursively implies the last conclusion of the lemma.
\end{proof}

\begin{lemma}\label{TildeDLemma}
For any iteration index $k\geq 1$ of the U-HLR call in Step~1, it holds that
\begin{equation}\label{eq:tildeM_uniform_bound}
\widetilde{\mathcal M}_k\le C(\tau_t,\beta_{t-1})
\end{equation}
where $C(\tau_t,\beta_{t-1})$ and $\tM{k}$ are as in \eqref{Cbartau} and \eqref{def:tildeM}, respectively.
\end{lemma}
\begin{proof}
Let $k\geq 1$ be an iteration index of the U-HLR call in Step~1. We claim that
\begin{equation}\label{TildeDBound}
   {\tD{k}}\leq 3\sqrt{2\tau_t}+1
\end{equation}
where $\tau_t$ and $\tD{k}$ are as in \eqref{eq:xspace_tau_t_def} and \eqref{DefsHLR}, respectively.
Hence, consider $U\in \mathcal S(\tilde U_{k-1};g_t)$ where $\mathcal S(\cdot; \cdot)$ is as in \eqref{SublevelDef}. It then follows from this inclusion, the definition of $g_t$ in \eqref{DefPhiT}, the last conclusion of Lemma~\ref{HLRInstance}(b), and the fact that $\tilde U_0= Y_{t-1}$ in Step~1.a) that
\[
\phi_t(UU^\top)\le\phi_t(\tilde U_{k-1}\tilde U^{\top}_{k-1})
\leq \phi_t(\tilde U_0\tilde U_0^{\top})=\phi_t(Y_{t-1}Y_{t-1}^{\top})
\]
which implies that $UU^{\top}\in \mathcal S(Z_{t-1};\phi_t)$ with $Z_{t-1}\defi Y_{t-1}Y_{t-1}^{\top}$. It then follows from this inclusion and the
definition of $\tau_t$ in \eqref{eq:xspace_tau_t_def} that $\|U\|^2=\operatorname{tr}(UU^\top)\le\tau_t/2$. This bound together with the definition of $\tD{k}$ in \eqref{DefsHLR} gives
\[
\tD{k}\le6\sqrt{\tau_t/2}+1=3\sqrt{2\tau_t}+1
\]
which proves \eqref{TildeDBound}.
It then follows from \eqref{TildeDBound} and the definitions of $C(\tau_t,\beta_{t-1})$ and $\tM{k}$ in \eqref{Cbartau} and \eqref{def:tildeM}, respectively, that relation \eqref{eq:tildeM_uniform_bound} must hold.
\end{proof}

We are now ready to prove Proposition~\ref{prop:hlr-general}.
\begin{proof}[Proof of Proposition~\ref{prop:hlr-general}]
(a) It follows from Lemma~\ref{HLRInstance}(a)
that the U-HLR call in Step~1 is an instance of the RUFW framework applied to $\phi_t$ with initial point $\tilde X_0=Y_{t-1}Y_{t-1}^\top=Z_{t-1}$ and tolerance pair $(\epsilon,\delta)=(\varepsilon_d,\varepsilon_c)$. Relations \eqref{eqn:AL_general} and \eqref{eqn:AL_gradient} imply that $\phi_t$ satisfies RUFW assumption~\textbf{(A1)} with $L=\beta_{t-1}\|\A\|^2$, while Proposition~\ref{prop:AL-sublevel-set-bounded} implies that $\phi_t$ satisfies RUFW assumption~\textbf{(A2)}. Moreover, the definitions of $\bar\tau$ in \eqref{eq:RUFW_level_set} and $\tau_t$ in \eqref{eq:xspace_tau_t_def}, together with \eqref{SublevelDef}, imply that $\bar\tau=\tau_t$. It then follows from Proposition~\ref{Prop:RUFW}(b) with $L=\beta_{t-1}\|\A\|^2$ and $\bar\tau=\tau_t$ that the U-HLR call terminates in at most \eqref{GeneralHLRComplexity} iterations with an iterate $X_k=U_kU_k^\top$ satisfying \eqref{SolutionTypeRUFW} and \eqref{descentRUFW} with $\phi=\phi_t$ and $\tilde X_0=Y_{t-1}Y_{t-1}^\top$. Since Step~\textbf{1.d)} sets $Y_t=U_k$, and since $S_k=\nabla\phi_t(U_kU_k^\top)$ and $\inner{S_kU_k}{U_k}=\inner{S_k}{U_kU_k^\top}$, relations \eqref{SolutionTypeRUFW} and \eqref{descentRUFW} are precisely \eqref{SolutionTypeHLR} and \eqref{descentHLR}. Statement (a) now follows from the fact that every iteration of the U-HLR call performs exactly one MEV computation in Step~\textbf{1.c)}.

(b) Let $k\geq 1$ be an iteration index of the U-HLR call in Step~1.
It follows from the complexity bound in Lemma~\ref{AIPPHLRLemma}(b) and the bound on $\tM{k}$ in \eqref{eq:tildeM_uniform_bound}
that the ADAP-AIPP call in Step~\textbf{1.b)} of the $k$-th iteration of U-HLR performs at most
\begin{equation}\label{AIPPGlobalHLR}
1+\frac{C_{\sigma}\max\left\{\eta_0^{-1},\,C(\tau_t,\beta_{t-1})\right\}}{\bar\rho^2}\left[\phi_t(\tilde U_{k-1}\tilde U^{\top}_{k-1})-\phi_t(U_kU_k^{\top})\right]+ \left\lceil\log_2^+\left(\eta_0C(\tau_t,\beta_{t-1})\right)\right\rceil
\end{equation}
iterations/ADAP-FISTA calls.

Let $\hat k$ denote the final iteration index of the U-HLR call.
Relation \eqref{eq:hlr_local_descent} implies that the total decrease
produced by its ADAP-AIPP calls is at most
\[
\phi_t(Y_{t-1}Y_{t-1}^\top)-\phi_t(Y_tY_t^\top).
\]
Combining this bound with \eqref{AIPPGlobalHLR},
\eqref{GeneralHLRComplexity}, and the definition of $\phi_t$ in
\eqref{DefPhiT} proves (b).
\end{proof}

\section{Numerical Results}
\label{sec:numerical-experiments}

In this section we report computational results for \cuourmethodplus,
together with comparisons with existing SDP solvers.
The problem classes are maximum stable set, matrix
completion, quadratic assignment problem (QAP), and general SDPs from the Mittelmann
benchmark collection \citep{mittelmann2003independent}. We first describe the implementations and the
comparison criteria. 

\paragraph{Implementation details}

The GPU implementation \cuourmethodplus uses C/CUDA. It is compiled for the NVIDIA H200 GPU architecture with compute capability 9.0 using CMake 3.30.2, GCC 12.3.0, and NVIDIA \texttt{nvcc} 12.6.68 (CUDA Toolkit 12.6.1), and linked against OpenBLAS 0.3.23. The GPU experiments use one NVIDIA H200 with 142 GB of GPU memory on a PACE node equipped with an Intel Xeon Platinum 8469C CPU.

For maximum stable set, \cuourmethodplus uses hybrid sparse low-rank
(HSLR) input \citep{aguirre2025user}. All cuHALLaR baselines are the
on-the-fly results reported in \citep{aguirre2025cuhallar}, using the
settings reported there.
These implementations, as well as the QAP implementations of
\cuourmethodplus and HALLaR, evaluate $CU$, $\mathcal A(UU^\top)$,
and $\mathcal A^*(p)U$ directly from the problem data without assembling
the SDP coefficient matrices.
cuLoRADS, CARDAL, and CSDP use SDPA input \citep{sdpa1}.
All solvers use this format for the Mittelmann collection.

\paragraph{Competing methods and comparison criteria}

The GPU experiments compare cuHALLaR+ with cuHALLaR, cuLoRADS v1.0.1
\citep{han2024accelerating} and the single-GPU implementation of CARDAL
v0.1.0 \citep{li2026curvatureawarerankadaptivedistributedaugmentedlagrangian}.

Given positive tolerances $(\varepsilon_p,\varepsilon_g,\varepsilon_d)$,
cu\ourmethodplus stops at $(\bar X,\bar p)$, with $\bar X\succeq0$ and
$\bar S\defi C-\mathcal A^*\bar p$, when
\begin{equation}\label{eq:numerical-hallar-stopping-test}
\frac{\|\mathcal A\bar X-b\|_2}{1+\|b\|_2}
\le\varepsilon_p,\qquad
\frac{|\inner{C}{\bar X}-\inner{b}{\bar p}|}
{1+|\inner{C}{\bar X}|+|\inner{b}{\bar p}|}
\le\varepsilon_g,\qquad
\frac{\max\{0,-\widehat\lambda_{\min}(\bar S)\}}
{1+|\inner{C}{\bar X}|+|\inner{b}{\bar p}|}\le\varepsilon_d,
\end{equation}
where $\widehat\lambda_{\min}$ is the computed smallest eigenvalue.
We set $\varepsilon_p=\varepsilon_g=\varepsilon_d=\epsilon$.
HALLaR and cuHALLaR use the trace-constrained stopping tests of
\citep{monteiro2024low} and \citep{aguirre2025cuhallar}, respectively.
The comparison target $\epsilon$ is $10^{-4}$ for QAP and $10^{-5}$ for
every other problem class. Solver-specific stopping tests and phase settings
are retained. The time limit is one hour for the Mittelmann SDPs, six hours
for the QAPLIB instances, and eight hours
for the random QAP instances.

Times are in seconds, and an entry time/rank also reports the final rank.
The cuLoRADS v1.0.1 times include its dual-infeasibility evaluations,
which cuLoRADS performs on the CPU after each phase and reports
separately from its solver time, whereas the Mittelmann benchmark
site\footnote{\url{https://plato.asu.edu/ftp/sparse_sdp.html}} reports
only the solver time. The dual residual is part of the stopping
test~\eqref{eq:numerical-hallar-stopping-test}, and the \cuourmethodplus
times include every evaluation of it. Therefore, we include the corresponding
cuLoRADS dual evaluation time within its final solve time. We use $\mathrm{TL}$ to denote an instance that hit the time-limit and \texttt{**} to denote a crash or an instance too large for the solver. An entry $*/N$ denotes an instance with reported maximum relative error $N$.

\subsection{Maximum stable set}
Given a graph $\mathcal{G} = ([n],E)$,
the maximum stable set problem consists of finding
a subset of vertices of largest cardinality
such that no two vertices are connected by an edge.
Lov\'asz \cite{lovasz1979shannon} introduced a constant, the $\vartheta$-function, which upper bounds the value of the maximum stable set.
The $\vartheta$-function is the value of the SDP
\begin{align*}
    \max \quad \{{e}{e}^{\top} \bullet X \; : \;
    X_{ij} = 0, \; ij \in E, \;
    \tr X = 1, \;
    X \succeq 0,\quad
    X \in \mathbb S^n(\mathbb{R})
    \}
\end{align*}
where $e = (1,1,\dots,1) \in \mathbb{R}^n$ is the all ones vector.
It was shown in \cite{grotschel1984polynomial} that the $\vartheta$-function agrees exactly with the stable set number for perfect graphs.

Table~\ref{tab:stable_set_hamming_hallar_plus} compares cuHALLaR+
with cuHALLaR from \citep{aguirre2025cuhallar}, cuLoRADS, and CARDAL on
the Hamming instances.

Both cuHALLaR+ and cuHALLaR solve all 14 Hamming instances and return
rank-two solutions throughout. cuHALLaR+ is faster on 12 instances;
cuHALLaR is faster on $H_{18,2}$ and $H_{19,2}$. cuLoRADS and CARDAL
report solves on the first five instances, and cuHALLaR+ is faster than
both on each of these. CARDAL is faster than cuLoRADS only on $H_{10,2}$.
On the largest graph, with over eight million vertices, cuHALLaR+
takes $44.13$ seconds, compared with $83.905$ seconds for cuHALLaR.

\begin{table}[!htbp]
\captionsetup{font=small}
\centering
\scriptsize
\begin{tabular}{@{}c|c|ccc|cc@{}}
\toprule
Graph $(n;|E|)$ & cuHALLaR+ & cuHALLaR & cuLoRADS & CARDAL & \shortstack{cuLoRADS/\\cuHALLaR+} & \shortstack{CARDAL/\\cuHALLaR+} \tabularnewline
\midrule
$H_{10,2}$(1,024; 5,120) & 0.552/2 & 1.735/2 & 1.05/14 & $0.83\,/\,19$ & {\textbf{\color{blue}1.90}} & {\textbf{\color{blue}1.50}} \tabularnewline
$H_{11,2}$(2,048; 11,264) & 0.48/2 & 1.627/2 & 1.23/16 & $1.96\,/\,21$ & {\textbf{\color{blue}2.56}} & {\textbf{\color{blue}4.08}} \tabularnewline
$H_{12,2}$(4,096; 24,576) & 1.57/2 & 2.033/2 & 1.91/17 & $7.01\,/\,23$ & {\textbf{\color{blue}1.22}} & {\textbf{\color{blue}4.46}} \tabularnewline
$H_{13,2}$(8,192; 53,248) & 0.95/2 & 3.106/2 & 5.56/19 & $54.72\,/\,24$ & {\textbf{\color{blue}5.85}} & {\textbf{\color{blue}57.60}} \tabularnewline
$H_{14,2}$(16,384; 114,688) & 0.93/2 & 3.170/2 & 15.44/20 & $539.45\,/\,26$ & {\textbf{\color{blue}16.60}} & {\textbf{\color{blue}580.05}} \tabularnewline
\midrule
$H_{15,2}$(32,768; 245,760) & 2.52/2 & 3.223/2 & \texttt{**} & \texttt{**} & -- & -- \tabularnewline
$H_{16,2}$(65,536; 524,288) & 3.69/2 & 4.457/2 & \texttt{**} & \texttt{**} & -- & -- \tabularnewline
$H_{17,2}$(131,072; 1,114,112) & 2.25/2 & 3.904/2 & \texttt{**} & \texttt{**} & -- & -- \tabularnewline
$H_{18,2}$(262,144; 2,359,296) & 3.05/2 & 2.627/2 & \texttt{**} & \texttt{**} & -- & -- \tabularnewline
$H_{19,2}$(524,288; 4,980,736) & 7.81/2 & 5.281/2 & \texttt{**} & \texttt{**} & -- & -- \tabularnewline
$H_{20,2}$(1,048,576; 10,485,760) & 5.19/2 & 8.109/2 & \texttt{**} & \texttt{**} & -- & -- \tabularnewline
$H_{21,2}$(2,097,152; 22,020,096) & 8.73/2 & 16.966/2 & ** & \texttt{**} & -- & -- \tabularnewline
$H_{22,2}$(4,194,304; 46,137,344) & 14.91/2 & 37.145/2 & \texttt{**} & \texttt{**} & -- & -- \tabularnewline
$H_{23,2}$(8,388,608; 96,468,992) & 44.13/2 & 83.905/2 & \texttt{**} & \texttt{**} & -- & -- \tabularnewline
\bottomrule
\end{tabular}
\caption{Runtimes for Hamming stable-set SDPs.}
\label{tab:stable_set_hamming_hallar_plus}
\end{table}

Table~\ref{tab:stable_set_gset_hallar_plus} reports the GSET results.
The cuHALLaR runs from \citep{aguirre2025cuhallar} use trace bound $\tau=1$.
Ratios are shown only for cuLoRADS
runs that pass the accuracy check.

Both cuHALLaR+ and cuHALLaR solve all 28 GSET instances, compared
with 26 for CARDAL and 19 for cuLoRADS. cuHALLaR+ is faster than cuHALLaR
on 24 instances. On the largest graph, G81, cuHALLaR+ returns a rank-two
solution in $48.91$ seconds, compared with a rank-three solution in
$64.20$ seconds for cuHALLaR.

\begin{table}[!htbp]
\captionsetup{font=small}
\centering
\scriptsize
\resizebox{\textwidth}{!}{%
\begin{tabular}{@{}c|c|ccc|c@{}}
\toprule
Graph $(n;|E|)$ & cuHALLaR+ & cuHALLaR & cuLoRADS & CARDAL & \shortstack{cuLoRADS/\\cuHALLaR+} \tabularnewline
\midrule
G1(800; 19,176) & 46.97/101 & 45.19/98 & 49.29/108 & $1232.91\,/\,121$ & {\textbf{\color{blue}1.05}} \tabularnewline
G10(800; 19,176) & 44.81/101 & 40.52/97 & $*/2.388\times10^{-5}$ & $1153.66\,/\,96$ & -- \tabularnewline
G11(800; 1,600) & 4.36/2 & 4.23/2 & 113.20/32 & $1.36\,/\,15$ & {\textbf{\color{blue}25.96}} \tabularnewline
G12(800; 1,600) & 3.19/2 & 2.59/2 & 1.68/14 & $0.79\,/\,15$ & 0.53 \tabularnewline
G14(800; 4,694) & 124.21/58 & 214.73/74 & 25.25/32 & $4.47\,/\,17$ & 0.20 \tabularnewline
G20(800; 4,672) & 165.73/170 & 359.40/166 & $*/2.646\times10^{-5}$ & $260.89\,/\,19$ & -- \tabularnewline
G43(1000; 9,990) & 27.30/63 & 35.84/60 & 41.15/72 & $968.02\,/\,68$ & {\textbf{\color{blue}1.51}} \tabularnewline
G51(1,000; 5,909) & 174.10/134 & 273.84/128 & 67.89/72 & $14.39\,/\,19$ & 0.39 \tabularnewline
G23(2,000; 19,990) & 35.91/79 & 52.59/76 & 119.17/81 & $*/6.536\times10^{-6}$ & {\textbf{\color{blue}3.32}} \tabularnewline
G31(2,000; 19,990) & 33.94/68 & 58.23/76 & 138.56/81 & $*/1.393\times10^{-5}$ & {\textbf{\color{blue}4.08}} \tabularnewline
G32(2,000; 4,000) & 22.89/3 & 25.34/3 & 137.63/36 & $2.79\,/\,19$ & {\textbf{\color{blue}6.01}} \tabularnewline
G34(2,000; 4,000) & 22.43/3 & 25.29/3 & 7.50/16 & $2.52\,/\,17$ & 0.33 \tabularnewline
G35(2,000; 11,778) & 203.53/117 & 254.68/108 & 127.86/54 & $83.46\,/\,19$ & 0.63 \tabularnewline
G41(2,000; 11,785) & 231.95/148 & 286.82/139 & 83.58/36 & $23.15\,/\,20$ & 0.36 \tabularnewline
G48(3,000; 6,000) & 8.31/2 & 18.95/2 & 19.43/17 & $4.80\,/\,19$ & {\textbf{\color{blue}2.34}} \tabularnewline
G55(5,000; 12,498) & 64.68/48 & 81.11/45 & 445.35/62 & $127.89\,/\,24$ & {\textbf{\color{blue}6.89}} \tabularnewline
G56(5,000; 12,498) & 67.33/43 & 80.83/45 & 428.62/41 & $426.77\,/\,22$ & {\textbf{\color{blue}6.37}} \tabularnewline
G57(5,000; 10,000) & 3.94/2 & 42.55/2 & $*/1.647\times10^{4}$ & $16.77\,/\,19$ & -- \tabularnewline
G58(5,000; 29,570) & 381.04/127 & 504.55/130 & $*/1.918\times10^{-4}$ & $3249.83\,/\,23$ & -- \tabularnewline
G59(5,000; 29,570) & 379.55/125 & 539.28/121 & $*/2.218\times10^{-4}$ & $3575.02\,/\,24$ & -- \tabularnewline
G60(7,000; 17,148) & 68.41/48 & 96.89/50 & 647.88/41 & $637.34\,/\,25$ & {\textbf{\color{blue}9.47}} \tabularnewline
G62(7,000; 14,000) & 7.03/3 & 25.18/3 & 2017.50/168 & $59.86\,/\,20$ & {\textbf{\color{blue}286.98}} \tabularnewline
G64(7,000; 41,459) & 295.47/55 & 334.12/51 & $*/2.316\times10^{-4}$ & $3235.03\,/\,22$ & -- \tabularnewline
G66(9,000; 18,000) & 22.38/2 & 27.91/3 & 1787.71/29 & $74.35\,/\,23$ & {\textbf{\color{blue}79.88}} \tabularnewline
G67(10,000; 20,000) & 4.36/2 & 5.83/2 & 1815.13/29 & $207.55\,/\,20$ & {\textbf{\color{blue}416.31}} \tabularnewline
G72(10,000; 20,000) & 5.03/2 & 5.81/2 & $*/7.960\times10^{-4}$ & $157.89\,/\,20$ & -- \tabularnewline
G77(14,000; 28,000) & 9.54/2 & 11.81/2 & $*/3.825\times10^{-4}$ & $636.46\,/\,21$ & -- \tabularnewline
G81(20,000; 40,000) & 48.91/2 & 64.20/3 & $*/3.724\times10^{-1}$ & $929.70\,/\,22$ & -- \tabularnewline
\bottomrule
\end{tabular}
}
\caption{Runtimes for GSET stable-set SDPs.}
\label{tab:stable_set_gset_hallar_plus}
\end{table}

\FloatBarrier

\subsection{Matrix completion}

Consider the problem of retrieving a low rank matrix $M \in \mathbb{R}^{n_1 \times n_2}$, where $n_1 \leq n_2$,
by observing a subset of its entries: $M_{ij}$, $ij\in \mathcal I$.
A standard approach to tackle this problem is by considering the nuclear norm relaxation:
\[
\min_Y \quad \left\{\|Y\|_*\quad
:\quad Y_{ij} = M_{ij}, \ \forall\, ij \in \mathcal I,\quad
Y \in \mathbb{R}^{n_1\times n_2}
\right \}
\]
The above problem can be rephrased as the following SDP:
\begin{align*}
\min_{X} \quad
\left\{ \,\frac12 \tr(X) \quad : \quad
X= \begin{pmatrix} W_1 & Y \\ Y^{\top} & W_2 \end{pmatrix} \succeq 0, \quad
Y_{i,j} = M_{i,j} \; \forall\, ij \in \mathcal I,\quad
X \in \mathbb S^{n_1+n_2}(\mathbb{R})
\right\}
\end{align*}
The nuclear norm relaxation was introduced in \cite{fazel2002matrix}.
It was shown in \cite{candes2012exact} it provably completes the matrix when $m = |\mathcal I|$ is sufficiently large
and the indices of the observations are independent and uniform.

Matrix completion instances are generated randomly using the protocol in \citep{aguirre2025cuhallar}: given $r\le n_1\le n_2$, the hidden solution matrix is $M=M_1M_2^\top$, where $M_1\in\mathbb{R}^{n_1\times r}$ and $M_2\in\mathbb{R}^{n_2\times r}$ have independent standard Gaussian entries. Then $m=\lceil r^2(n_1+n_2)\log(n_1+n_2)\rceil$ entries are sampled uniformly at random; the oversampling ratio is $r\log(n_1+n_2)$.

In Table~\ref{tab:matrix_completion_hallar_plus}, both cuHALLaR+ and
cuHALLaR solve all 26 instances, while cuLoRADS and CARDAL solve 25
and 18, respectively. cuHALLaR+ is faster than cuHALLaR on all 26
instances and faster than cuLoRADS and CARDAL on every instance for which
they report a solve. The final ranks of both cuHALLaR+ and cuHALLaR
agree with the planted rank on 24 instances and exceed it by one on
the remaining two. On the largest instance, with over 261 million
observations, cuHALLaR+ takes $30.82$ seconds, compared with $53.3$
seconds for cuHALLaR and $4{,}124.44$ seconds for cuLoRADS.

\begin{table}[!htbp]
\captionsetup{font=small}
\centering
\scriptsize
\resizebox{\textwidth}{!}{%
\begin{tabular}{@{}cc|c|ccc|c@{}}
\toprule
$(n_1,n_2;m)$ & $r$ & cuHALLaR+ & cuHALLaR & cuLoRADS & CARDAL & \shortstack{cuLoRADS/\\cuHALLaR+} \tabularnewline
\midrule
3,000, 7,000; 828,931 & 3 & 0.53/3 & 0.77/3 & 3.12/19 & $0.79\,/\,28$ & {\textbf{\color{blue}5.89}} \tabularnewline
3,000, 7,000; 828,931 & 3 & 0.64/3 & 0.70/3 & 3.14/19 & $*/1.137\times10^{-4}$ & {\textbf{\color{blue}4.91}} \tabularnewline
3,000, 7,000; 2,302,586 & 5 & 0.83/5 & 1.28/5 & 3.79/19 & $1.15\,/\,35$ & {\textbf{\color{blue}4.57}} \tabularnewline
3,000, 7,000; 2,302,586 & 5 & 0.79/5 & 1.09/5 & 3.42/19 & $1.08\,/\,35$ & {\textbf{\color{blue}4.33}} \tabularnewline
10,000, 21,623; 2,948,996 & 3 & 1.04/3 & 1.22/3 & 9.36/21 & $2.92\,/\,30$ & {\textbf{\color{blue}9.00}} \tabularnewline
10,000, 21,623; 2,948,996 & 3 & 1.02/3 & 1.11/3 & 9.21/21 & $2.90\,/\,30$ & {\textbf{\color{blue}9.03}} \tabularnewline
10,000, 21,623; 8,191,654 & 5 & 1.68/5 & 2.05/5 & 12.64/21 & $487.26\,/\,40$ & {\textbf{\color{blue}7.52}} \tabularnewline
10,000, 21,623; 8,191,654 & 5 & 1.73/5 & 1.95/5 & 10.92/21 & $8.40\,/\,37$ & {\textbf{\color{blue}6.31}} \tabularnewline
25,000, 50,000; 3,367,574 & 2 & 1.17/2 & 1.20/2 & 16.72/23 & $5.90\,/\,31$ & {\textbf{\color{blue}14.29}} \tabularnewline
25,000, 50,000; 7,577,040 & 3 & 1.41/3 & 1.82/3 & 21.14/23 & $11.72\,/\,35$ & {\textbf{\color{blue}14.99}} \tabularnewline
30,000, 70,000; 4,605,171 & 2 & 1.59/2 & 1.63/2 & 25.79/24 & $8.59\,/\,31$ & {\textbf{\color{blue}16.22}} \tabularnewline
30,000, 70,000; 10,361,633 & 3 & 2.01/3 & 2.57/3 & 27.13/24 & $19.06\,/\,34$ & {\textbf{\color{blue}13.50}} \tabularnewline
50,000, 100,000; 7,151,035 & 2 & 1.42/2 & 1.89/2 & 40.20/24 & $933.59\,/\,33$ & {\textbf{\color{blue}28.31}} \tabularnewline
50,000, 100,000; 16,089,828 & 3 & 2.57/3 & 3.68/3 & 40.13/24 & $2672.44\,/\,37$ & {\textbf{\color{blue}15.61}} \tabularnewline
80,000, 120,000; 9,764,859 & 2 & 2.09/2 & 2.45/2 & 50.41/25 & $27.80\,/\,35$ & {\textbf{\color{blue}24.12}} \tabularnewline
80,000, 120,000; 21,970,931 & 3 & 4.86/3 & 5.02/3 & 38.35/25 & $63.03\,/\,37$ & {\textbf{\color{blue}7.89}} \tabularnewline
120,000, 180,000; 34,051,152 & 3 & 5.43/4 & 7.35/4 & 87.41/26 & $127.07\,/\,38$ & {\textbf{\color{blue}16.10}} \tabularnewline
160,000, 240,000; 46,437,192 & 3 & 4.85/4 & 8.79/4 & 409.53/26 & $253.83\,/\,37$ & {\textbf{\color{blue}84.44}} \tabularnewline
240,000, 360,000; 71,845,299 & 3 & 10.74/3 & 12.82/3 & 744.89/27 & $2797.54\,/\,37$ & {\textbf{\color{blue}69.36}} \tabularnewline
\midrule
320,000, 480,000; 97,865,043 & 3 & 14.43/3 & 17.41/3 & 1187.79/28 & $*/9.668\times10^{-1}$ & {\textbf{\color{blue}82.31}} \tabularnewline
400,000, 600,000; 124,339,596 & 3 & 18.13/3 & 22.31/3 & 2037.88/28 & $*/3.377\times10^{-5}$ & {\textbf{\color{blue}112.40}} \tabularnewline
480,000, 720,000; 151,176,587 & 3 & 24.41/3 & 28.13/3 & $*/1.555\times10^{-6}$ & $*/3.338\times10^{-1}$ & -- \tabularnewline
560,000, 840,000; 178,314,984 & 3 & 30.69/3 & 34.72/3 & 460.49/29 & $*/9.782\times10^{-1}$ & {\textbf{\color{blue}15.00}} \tabularnewline
640,000, 960,000; 205,711,405 & 3 & 30.36/3 & 39.25/3 & 3102.16/29 & $*/9.802\times10^{-1}$ & {\textbf{\color{blue}102.18}} \tabularnewline
720,000, 1,080,000; 233,333,416 & 3 & 30.59/3 & 48.85/3 & 482.13/29 & $*/9.951\times10^{-1}$ & {\textbf{\color{blue}15.76}} \tabularnewline
800,000, 1,200,000; 261,155,840 & 3 & 30.82/3 & 53.3/3 & 4124.44/30 & $*/9.999\times10^{-1}$ & {\textbf{\color{blue}133.82}} \tabularnewline
\bottomrule
\end{tabular}%
}
\caption{Runtimes for matrix-completion SDPs. The cuHALLaR results are
reproduced from \citep[Table~2]{aguirre2025cuhallar}, using
tolerance $10^{-5}$ and a three-hour time limit.}
\label{tab:matrix_completion_hallar_plus}
\end{table}

\FloatBarrier

\subsection{Quadratic assignment problem}
\label{subsec:qap}

Given matrices $A,B\in\mathbb S^\nu$, the quadratic assignment problem
(QAP) consists of minimizing $\tr(APBP^{\top})$ over the $\nu\times\nu$
permutation matrices $P$. With column-wise vectorization, the identity
\[
\tr(APBP^{\top})=\operatorname{vec}(P)^{\top}(B\otimes A)\operatorname{vec}(P)
\]
motivates the following SDP relaxation:
\begin{equation}\label{eq:qap-sdp-relaxation}
\min_X\quad\left\{(B\otimes A)\bullet X\;:\;
\begin{array}{l}
\operatorname{tr}_1X=I_\nu,\quad \operatorname{tr}_2X=I_\nu,\\
e_{\nu^2}e_{\nu^2}^{\top}\bullet X=\nu^2,\quad
X\succeq0,\quad X\in\mathbb S^{\nu^2}(\mathbb R)
\end{array}\right\}.
\end{equation}
Here $e_{\nu^2}\in\mathbb R^{\nu^2}$ is the all-ones vector. The index
$(i,j)$ denotes the entry $i+(j-1)\nu$ of a vectorized $\nu\times\nu$ matrix,
and the partial traces are defined by
\[
[\operatorname{tr}_1X]_{ij}=\sum_{\ell=1}^\nu X_{(i,\ell),(j,\ell)},
\qquad
[\operatorname{tr}_2X]_{ij}=\sum_{\ell=1}^\nu X_{(\ell,i),(\ell,j)}.
\]
For every permutation matrix $P$, the rank-one matrix
$X=\operatorname{vec}(P)\operatorname{vec}(P)^{\top}$ is feasible for
\eqref{eq:qap-sdp-relaxation} and has the same objective value. Thus,
\eqref{eq:qap-sdp-relaxation} is a relaxation of QAP and is equivalent to
that in \citep[Section~4.6]{monteiro2024low}, up to the ordering of the
Kronecker factors. Its matrix dimension is $n=\nu^2$. For $\nu\ge2$, we omit
the redundant final diagonal equation of $\operatorname{tr}_2X=I_\nu$,
leaving $m=\nu(\nu+1)$ independent scalar constraints. The partial-trace
constraints imply $\tr X=\nu$, but \ourmethodplus does not require this trace bound as an input. HALLaR uses the trace bound $\tau=\nu$.

For nonsymmetric QAPLIB data, the objective matrix is replaced by
$0.5(B\otimes A+B^{\top}\otimes A^{\top})$. This symmetrization preserves the
quadratic form and leaves the constraints in
\eqref{eq:qap-sdp-relaxation} unchanged.

Table~\ref{tab:qap_qaplib_hallar_plus} covers 32 QAPLIB instances, including
the 11 \texttt{tai} instances in \citep[Table~10]{monteiro2024low}.
All methods use a six-hour time limit.

cuHALLaR+ starts from $U_0=e_{\nu^2}/\sqrt{\nu}$ and the zero multiplier.
The objective is scaled
by $1/(\|A\|_F\|B\|_F)$, and the same scaling is used in the stopping
test~\eqref{eq:numerical-hallar-stopping-test}.
For all instances after the first 11, CSDP receives the independent equality
constraints and the redundant constraint $\tr X+s=\nu$, $s\ge0$. 

CSDP reaches the target accuracy on 23 of the 32 instances, including the
nine smallest \texttt{tai} instances. HALLaR reaches the target on seven
instances. The final errors reported by cuLoRADS exceed the target on all
32 instances.

\begin{table}[!htbp]
\centering
\captionsetup{font=small}
\small
\setlength{\tabcolsep}{4pt}
\begin{tabular}{@{}lr|rr|rr@{}}
\toprule
Instance & $n$ & HALLaR & CSDP & cuHALLaR+ & cuLoRADS\\
\midrule
\texttt{tai12a} & 144 & $*/4.6\times10^{-4}$ & 0.333 & 2.753/17 & $*/1.464\times10^{-2}$\\
\texttt{tai15a} & 225 & 96.376/3 & 0.858 & 7.582/21 & $*/1.254\times10^{4}$\\
\texttt{tai20a} & 400 & 188.321/2 & 3.979 & 6.391/28 & $*/9.061\times10^{-3}$\\
\texttt{tai25a} & 625 & $*/3.5\times10^{-4}$ & 13.998 & 9.228/35 & $*/4.462\times10^{-2}$\\
\texttt{tai30a} & 900 & 695.807/3 & 43.371 & 6.541/42 & $*/3.809\times10^{-2}$\\
\texttt{tai35a} & 1,225 & $*/4.4\times10^{-4}$ & 121.273 & 9.114/49 & $*/5.694\times10^{-2}$\\
\texttt{tai40a} & 1,600 & $*/8.9\times10^{-4}$ & 335.869 & 12.812/56 & $*/4.236\times10^{-2}$\\
\midrule
\texttt{tai50a} & 2,500 & $*/1.9\times10^{-4}$ & 1,894.414 & 18.169/70 & $*/3.573\times10^{-2}$\\
\texttt{tai60a} & 3,600 & $*/4.7\times10^{-3}$ & 6,653.269 & 31.262/85 & $*/4.454\times10^{-2}$\\
\texttt{tai80a} & 6,400 & $*/6.0\times10^{-3}$ & $*/5.260\times10^{-1}$ & 117.579/113 & $*/1.581\times10^{-2}$\\
\texttt{tai100a} & 10,000 & $*/6.9\times10^{-3}$ & $*/5.721$ & 174.861/141 & $*/4.399\times10^{-2}$\\
\midrule
\texttt{esc32a} & 1,024 & $*/9.3\times10^{-4}$ & 150.310 & 6.423/45 & $*/2.285\times10^{-2}$\\
\texttt{esc32b} & 1,024 & $*/3.8\times10^{-3}$ & $*/2.021\times10^{-3}$ & 2.406/45 & $*/2.651\times10^{-1}$\\
\texttt{esc32c} & 1,024 & $*/2.9\times10^{-3}$ & 171.951 & 3.735/45 & $*/7.398\times10^{-2}$\\
\texttt{esc32d} & 1,024 & $*/8.6\times10^{-4}$ & 175.683 & 6.130/45 & $*/3.270\times10^{-1}$\\
\texttt{esc32e} & 1,024 & $*/2.3\times10^{-3}$ & $*/2.990\times10^{-4}$ & 3.022/45 & $*/7.810\times10^{-3}$\\
\texttt{esc32g} & 1,024 & $*/2.3\times10^{-3}$ & $*/1.856\times10^{-4}$ & 2.565/45 & $*/9.339\times10^{-3}$\\
\texttt{esc32h} & 1,024 & $*/3.8\times10^{-4}$ & 185.051 & 3.408/45 & $*/8.705\times10^{-2}$\\
\texttt{kra30a} & 900 & $*/1.1\times10^{-4}$ & $*/3.935\times10^{-4}$ & 11.682/42 & $*/9.852\times10^{-2}$\\
\texttt{kra30b} & 900 & $*/5.4\times10^{-4}$ & 46.562 & 8.673/42 & $*/1.164\times10^{-1}$\\
\texttt{kra32} & 1,024 & $*/1.3\times10^{-3}$ & $*/2.050\times10^{-4}$ & 14.626/45 & $*/2.988\times10^{-1}$\\
\texttt{lipa30a} & 900 & $*/1.8\times10^{-3}$ & 38.018 & 6.417/42 & $*/1.371\times10^{-2}$\\
\texttt{lipa30b} & 900 & 2,748.815/24 & 39.653 & 7.326/42 & $*/9.520\times10^{-3}$\\
\texttt{nug21} & 441 & $*/6.0\times10^{-4}$ & 5.230 & 6.462/29 & $*/2.180\times10^{-2}$\\
\texttt{nug22} & 484 & $*/4.3\times10^{-4}$ & 7.235 & 9.317/31 & $*/3.750\times10^{-2}$\\
\texttt{nug24} & 576 & $*/2.0\times10^{-4}$ & 11.002 & 14.059/34 & $*/9.982\times10^{-1}$\\
\texttt{nug25} & 625 & 1,732.228/35 & $*/1.495\times10^{-4}$ & 4.876/35 & $*/1.938\times10^{-2}$\\
\texttt{nug27} & 729 & $*/2.1\times10^{-3}$ & 23.242 & 18.794/38 & $*/2.282\times10^{-2}$\\
\texttt{nug28} & 784 & 2,696.431/39 & $*/1.246\times10^{-4}$ & 20.518/39 & $*/1.407$\\
\texttt{nug30} & 900 & 5,652.199/42 & 41.891 & 30.640/42 & $*/1.492\times10^{-2}$\\
\texttt{tai25b} & 625 & $*/2.8\times10^{-3}$ & 14.694 & 12.814/35 & $*/3.355\times10^{-1}$\\
\texttt{tho30} & 900 & $*/9.9\times10^{-4}$ & 41.601 & 16.543/42 & $*/5.090\times10^{-2}$\\
\bottomrule
\end{tabular}
\caption{QAPLIB SDP results.}
\label{tab:qap_qaplib_hallar_plus}
\par\smallskip
\end{table}

In Table~\ref{tab:qap_qaplib_remaining1}, CSDP reaches the target on 15 of
18 instances, while HALLaR does so on four. cuLoRADS finishes above the target
on all 18 instances. On the two largest instances, CSDP takes about
$10{,}400$ seconds on \texttt{esc64a} and has no completed result for
\texttt{esc128}.

\begin{table}[!htbp]
\centering
\captionsetup{font=small}
\small
\setlength{\tabcolsep}{4pt}
\begin{tabular}{@{}lr|rr|rr@{}}
\toprule
Instance & $n$ & HALLaR & CSDP & cuHALLaR+ & cuLoRADS\\
\midrule
\texttt{chr12b} & 144 & $*/7.7\times10^{-4}$ & $*/2.539\times10^{-4}$ & 2.088/17 & $*/1.832\times10^{-1}$\\
\texttt{chr12c} & 144 & $*/4.5\times10^{-4}$ & 0.319 & 2.854/17 & $*/2.103$\\
\texttt{chr15a} & 225 & $*/2.3\times10^{-4}$ & 0.815 & 3.569/21 & $*/8.525\times10^{-2}$\\
\texttt{chr15b} & 225 & $*/5.0\times10^{-4}$ & 0.789 & 4.608/21 & $*/7.327\times10^{-2}$\\
\texttt{chr15c} & 225 & 462.838/18 & 0.792 & 8.163/21 & $*/4.291\times10^{-2}$\\
\texttt{chr18a} & 324 & $*/8.0\times10^{-4}$ & 2.084 & 4.675/25 & $*/5.505\times10^{-2}$\\
\texttt{chr18b} & 324 & $*/1.1\times10^{-3}$ & $*/2.776\times10^{-4}$ & 4.298/25 & $*/1.070\times10^{-1}$\\
\texttt{chr20a} & 400 & $*/1.0\times10^{-3}$ & 3.786 & 6.732/28 & $*/4.751\times10^{-2}$\\
\texttt{chr20b} & 400 & $*/1.5\times10^{-4}$ & 3.379 & 3.979/28 & $*/3.675\times10^{-2}$\\
\texttt{chr20c} & 400 & $*/1.9\times10^{-4}$ & 3.851 & 5.675/28 & $*/1.228\times10^{-1}$\\
\texttt{els19} & 361 & $*/8.0\times10^{-4}$ & 3.076 & 5.099/27 & $*/8.363\times10^{-2}$\\
\texttt{esc64a} & 4,096 & $*/1.2\times10^{-3}$ & 10,400.143 & 328.053/90 & $*/6.617\times10^{-2}$\\
\texttt{esc128} & 16,384 & $*/7.7\times10^{-3}$ & \texttt{**} & 96.435/181 & $*/9.880\times10^{-1}$\\
\texttt{had12} & 144 & 81.567/6 & 0.348 & 4.673/17 & $*/4.059\times10^{-2}$\\
\texttt{had14} & 196 & 94.103/2 & 0.614 & 6.275/20 & $*/6.711\times10^{-2}$\\
\texttt{had16} & 256 & 193.668/8 & 1.324 & 5.277/22 & $*/6.728\times10^{-2}$\\
\texttt{had18} & 324 & $*/1.0\times10^{-4}$ & 2.297 & 5.746/25 & $*/6.602\times10^{-2}$\\
\texttt{had20} & 400 & $*/2.0\times10^{-3}$ & 3.973 & 7.221/28 & $*/2.495\times10^{-2}$\\
\bottomrule
\end{tabular}
\caption{Remaining QAPLIB instances of \citep[Table~11]{monteiro2024low},
part~1 of~2.}
\label{tab:qap_qaplib_remaining1}
\end{table}

CSDP reaches the target on 18 of the 31 instances in
Table~\ref{tab:qap_qaplib_remaining2}, while cuLoRADS does so only on
\texttt{sko64}. CSDP misses the accuracy target on all six \texttt{sko100}
instances. On \texttt{tho150}, the largest instance in this group,
cuHALLaR+ reaches the time limit.

\begin{table}[!htbp]
\centering
\captionsetup{font=small}
\small
\setlength{\tabcolsep}{4pt}
\begin{tabular}{@{}lr|rr|rr@{}}
\toprule
Instance & $n$ & HALLaR & CSDP & cuHALLaR+ & cuLoRADS\\
\midrule
\texttt{nug12} & 144 & $*/2.8\times10^{-4}$ & 0.312 & 2.720/17 & $*/2.045\times10^{-2}$\\
\texttt{nug14} & 196 & $*/3.9\times10^{-4}$ & 0.618 & 2.683/20 & $*/3.877\times10^{-2}$\\
\texttt{nug15} & 225 & 151.787/7 & 0.817 & 2.162/21 & $*/4.266\times10^{-2}$\\
\texttt{nug16a} & 256 & $*/3.9\times10^{-4}$ & 1.269 & 4.111/22 & $*/8.999\times10^{-2}$\\
\texttt{nug16b} & 256 & $*/2.9\times10^{-4}$ & 1.199 & 4.089/22 & $*/7.426\times10^{-2}$\\
\texttt{nug17} & 289 & 254.685/7 & 1.538 & 3.864/24 & $*/8.621\times10^{-2}$\\
\texttt{nug18} & 324 & 417.714/14 & 2.097 & 3.865/25 & $*/8.741\times10^{-2}$\\
\texttt{nug20} & 400 & $*/6.5\times10^{-4}$ & 3.641 & 6.658/28 & $*/1.245\times10^{-1}$\\
\texttt{rou12} & 144 & 70.956/3 & 0.329 & 1.825/17 & $*/5.214\times10^{1}$\\
\texttt{rou15} & 225 & $*/1.8\times10^{-4}$ & 0.800 & 6.946/21 & $*/6.922\times10^{-2}$\\
\texttt{rou20} & 400 & 194.810/4 & 3.818 & 6.890/28 & $*/4.727\times10^{-2}$\\
\texttt{scr12} & 144 & $*/1.1\times10^{-3}$ & $*/8.053\times10^{-4}$ & 6.117/17 & $*/5.735\times10^{-1}$\\
\texttt{scr15} & 225 & $*/6.5\times10^{-4}$ & $*/1.031\times10^{-3}$ & 3.267/21 & $*/4.563\times10^{-1}$\\
\texttt{scr20} & 400 & $*/2.3\times10^{-3}$ & $*/8.460\times10^{-4}$ & 4.588/28 & $*/4.714\times10^{-1}$\\
\texttt{sko42} & 1,764 & $*/7.3\times10^{-3}$ & 425.000 & 20.370/59 & $*/1.270\times10^{-2}$\\
\texttt{sko49} & 2,401 & $*/3.0\times10^{-3}$ & 1,446.965 & 27.557/69 & $*/3.462\times10^{-2}$\\
\texttt{sko56} & 3,136 & $*/2.5\times10^{-3}$ & 3,845.140 & 30.136/79 & $*/1.000$\\
\texttt{sko64} & 4,096 & $*/5.1\times10^{-3}$ & 12,333.661 & 70.466/90 & 134.090\\
\texttt{sko100a} & 10,000 & $*/6.8\times10^{-3}$ & $*/2.184\times10^{4}$ & 202.732/141 & $*/2.454\times10^{-2}$\\
\texttt{sko100b} & 10,000 & $*/5.4\times10^{-3}$ & $*/3.144\times10^{3}$ & 205.820/143 & $*/2.401\times10^{-2}$\\
\texttt{sko100c} & 10,000 & $*/9.7\times10^{-3}$ & $*/1.955\times10^{4}$ & 209.715/142 & $*/5.477\times10^{-2}$\\
\texttt{sko100d} & 10,000 & $*/9.7\times10^{-3}$ & $*/2.773\times10^{3}$ & 216.082/142 & $*/8.624\times10^{-3}$\\
\texttt{sko100e} & 10,000 & $*/9.7\times10^{-3}$ & $*/2.625\times10^{3}$ & 215.098/141 & $*/1.416\times10^{-2}$\\
\texttt{sko100f} & 10,000 & $*/9.7\times10^{-3}$ & $*/2.873\times10^{3}$ & 217.526/143 & $*/3.214\times10^{-2}$\\
\texttt{ste36a} & 1,296 & $*/1.4\times10^{-3}$ & 139.864 & 170.791/51 & $*/9.751\times10^{-2}$\\
\texttt{ste36b} & 1,296 & $*/2.5\times10^{-4}$ & $*/1.985\times10^{-4}$ & 100.028/51 & $*/6.575\times10^{-2}$\\
\texttt{ste36c} & 1,296 & 5,445.863/33 & 141.673 & 100.777/51 & $*/1.368\times10^{-1}$\\
\texttt{tho40} & 1,600 & $*/2.9\times10^{-3}$ & $*/1.669\times10^{-4}$ & 13.467/56 & $*/4.046$\\
\texttt{tho150} & 22,500 & $*/2.1\times10^{-3}$ & \texttt{**} & $\mathrm{TL}$ & $*/3.218\times10^{-1}$\\
\texttt{wil50} & 2,500 & $*/4.3\times10^{-3}$ & 1,874.664 & 202.683/70 & $*/3.859\times10^{-2}$\\
\texttt{wil100} & 10,000 & $*/7.4\times10^{-3}$ & $*/2.975\times10^{4}$ & 220.789/141 & $*/2.453\times10^{-2}$\\
\bottomrule
\end{tabular}
\caption{Remaining QAPLIB instances of \citep[Table~11]{monteiro2024low},
part~2 of~2.}
\label{tab:qap_qaplib_remaining2}
\end{table}

Table~\ref{tab:qap_random_hallar_plus} reports results at the 11 random
QAP dimensions in \citep[Table~12]{monteiro2024low}. For our runs, $A$
and $B$ are symmetric with zero diagonal; each strict upper-triangular
entry is obtained by rounding an independent uniform draw on $[0,10]$.
We generate one instance per size using Julia~1.10.1 and seed~0, and
solve \eqref{eq:qap-sdp-relaxation} with a tolerance of $10^{-4}$
and an eight-hour solver limit. HALLaR and cuHALLaR+ use the same
Gaussian unit-norm initial factor.

HALLaR, CSDP, and cuLoRADS record no solves at the target accuracy on the
random instances. At $\nu=400$, the SDP has dimension $160{,}000$ and
the reported HALLaR error is $2.5\times10^{-3}$.
For cuLoRADS, the $\nu=200$ instance with $n=40{,}000$ exhausts GPU memory
on the H200 during solver initialization.

\begin{table}[!htbp]
\centering
\captionsetup{font=small}
\small
\setlength{\tabcolsep}{4pt}
\begin{tabular}{@{}rr|rr|rr@{}}
\toprule
$\nu$ & $n$ & HALLaR & CSDP & cuHALLaR+ & cuLoRADS\\
\midrule
100 & 10,000 & $*/3.6\times10^{-3}$ & $*/1.159\times10^{2}$ & 210.226/141 & $*/2.103\times10^{-2}$\\
120 & 14,400 & $*/3.1\times10^{-3}$ & $*/1.452\times10^{5}$ & 199.461/169 & $*/2.972\times10^{-2}$\\
150 & 22,500 & $*/2.9\times10^{-3}$ & \texttt{**} & 425.046/192 & $*/1.000$\\
175 & 30,625 & $*/4.8\times10^{-3}$ & \texttt{**} & 605.946/227 & $*/1.000$\\
200 & 40,000 & $*/3.5\times10^{-3}$ & \texttt{**} & 826.315/263 & \texttt{**}\\
250 & 62,500 & $*/3.6\times10^{-3}$ & \texttt{**} & 1,381.325/253 & \texttt{**}\\
275 & 75,625 & $*/3.6\times10^{-3}$ & \texttt{**} & 1,724.889/289 & \texttt{**}\\
300 & 90,000 & $*/2.6\times10^{-3}$ & \texttt{**} & 2,107.618/324 & \texttt{**}\\
350 & 122,500 & $*/2.5\times10^{-3}$ & \texttt{**} & 3,011.551/395 & \texttt{**}\\
375 & 140,625 & $*/2.0\times10^{-3}$ & \texttt{**} & 3,529.608/330 & \texttt{**}\\
400 & 160,000 & $*/2.5\times10^{-3}$ & \texttt{**} & 4,096.032/465 & \texttt{**}\\
\bottomrule
\end{tabular}
\caption{Random QAP SDP results.}
\label{tab:qap_random_hallar_plus}
\end{table}

\FloatBarrier

\subsection{General SDP benchmark}
\label{subsec:general-sdp-benchmark}

Table~\ref{tab:general-sdp-benchmark} reports results for the same 27
Mittelmann SDP instances \citep{mittelmann2003independent}
considered in \citep[Table~7]{han2024accelerating}.
CARDAL solves 14 of the 27 instances and reaches the one-hour limit
on the remaining 13. cuLoRADS has reported runtimes on 15 instances.

\begin{table}[htbp]
\centering
\captionsetup{font=small}
\scriptsize
\setlength{\tabcolsep}{2.5pt}
\renewcommand{\arraystretch}{0.96}
\begin{tabular}{@{}lrr|ccc@{}}
\toprule
Instance & $n$ & $m$ & cuHALLaR+ & cuLoRADS & CARDAL \\
\midrule
1dc.1024 & 1024 & 24064 & 43.92 & $*/2.485\times10^{-6}$ & 1513.56 \\
1et.2048 & 2048 & 22529 & 65.62 & $*/4.858\times10^{-6}$ & 908.64 \\
1tc.2048 & 2048 & 18945 & 44.28 & 180.96 & 763.98 \\
1zc.1024 & 1024 & 16641 & 19.95 & $*/1.970\times10^{-6}$ & 276.82 \\
AlH & 5990 & 7230 & 6.60 & $*/1.621\times10^{-5}$ & $*/6.011\times10^{-3}$ \\
BH2 & 2166 & 1743 & 7.12 & $*/9.283\times10^{-5}$ & $*/1.017\times10^{-3}$ \\
cancer\_100 & 569 & 10469 & 42.79 & $*/1.578\times10^{-6}$ & 624.05 \\
CH2 & 2166 & 1743 & 67.89 & $*/2.434\times10^{-5}$ & $*/1.306\times10^{-3}$ \\
checker\_1.5 & 3970 & 3970 & 7.60 & 28.90 & 210.92 \\
cphil12 & 364 & 12376 & 40.69 & $*/4.428\times10^{-6}$ & $*/2.969\times10^{-3}$ \\
fap09 & 28224 & 30276 & 242.79 & 548.29 & $*/2.506\times10^{-2}$ \\
foot & 2208 & 2209 & 95.39 & 74.28 & $*/1.150\times10^{-6}$ \\
G40\_mb & 2000 & 2001 & 15.12 & 59.20 & 1456.10 \\
G48\_mb & 3000 & 3001 & 4.82 & 15.35 & 542.90 \\
G60\_mb & 7000 & 7001 & 35.51 & 141.42 & $*/2.775\times10^{-5}$ \\
H3O & 3162 & 2964 & 36.95 & $*/1.642\times10^{-5}$ & $*/6.406\times10^{-4}$ \\
hand & 1296 & 1297 & 3.09 & 4.22 & 94.26 \\
ice\_2.0 & 8113 & 8113 & 51.47 & 72.39 & 236.10 \\
NH2 & 2166 & 1743 & 127.89 & $*/1.364\times10^{-5}$ & $*/1.377\times10^{-3}$ \\
NH3 & 3162 & 2964 & 126.95 & $*/2.189\times10^{-5}$ & $*/7.693\times10^{-4}$ \\
NH4 & 4426 & 4743 & 46.20 & $*/8.608\times10^{-6}$ & $*/2.636\times10^{-2}$ \\
p\_auss2\_3.0 & 9115 & 9115 & 28.75 & 68.08 & 215.86 \\
rendl1 & 2000 & 2001 & 10.39 & 11.18 & $*/3.668\times10^{-6}$ \\
sensor\_1000 & 6528 & 5549 & 54.97 & 131.24 & $*/4.356\times10^{-2}$ \\
theta12 & 600 & 17979 & 29.86 & 32.28 & 919.99 \\
theta102 & 500 & 37467 & 15.28 & 17.54 & 1560.72 \\
theta123 & 600 & 90020 & 46.37 & 52.91 & 2428.30 \\
\bottomrule
\end{tabular}
\caption{Runtimes for Mittelmann SDPs.}
\label{tab:general-sdp-benchmark}
\end{table}

\FloatBarrier

\section{Concluding remarks}\label{sec:concluding-remarks}
This paper develops HALLaR+, a low-rank inexact augmented Lagrangian method for SDPs with possibly unbounded feasible regions. A key practical advantage is that the method does not require a user-specified trace bound $\tau$. Choosing such a bound can be difficult: a value that is too small may exclude all optimal solutions, whereas an overly conservative value may substantially increase the computational effort. Under the stated assumptions, we establish a uniform bound on the diameters of the relevant augmented Lagrangian sublevel sets and use it to derive iteration-complexity guarantees for obtaining approximate primal–dual solutions.

The numerical experiments demonstrate strong computational performance on several classes of SDPs, with substantial speedups over existing solvers on many instances.
On the Mittelmann benchmarks, cuHALLaR+ achieves speedups of up to approximately $4.09\times$ over cuLoRADS and $112.63\times$ over CARDAL among instances successfully solved by the respective methods.
It also solves the largest random QAP relaxation tested, with matrix dimension $160,000$, in approximately 68 minutes.

Future work includes extending HALLaR+ to SDPs with inequality constraints and multiple positive semidefinite blocks, as well as to doubly nonnegative programs. For the latter, an important challenge is to handle the entrywise nonnegativity constraints while preserving the memory advantages of the low-rank factorization.

\newpage
\appendix
\section{\texorpdfstring{$X$}{X}-Space and Convex-Analysis Machinery}
\label{app:x-space-machinery}

This appendix collects the technical convex-analysis facts and describes the relaxed
unbounded Frank--Wolfe method used in our $X$-space analysis.

\subsection{Technical Convex-Analysis Facts}\label{app:technical-al-facts}

This subsection collects the convex-analysis facts used in the X-AL analysis
of Section~\ref{subsec:xspace-step-one}.

The next proposition establishes the boundedness property used for the
AL sublevel sets in the X-AL analysis.

\begin{proposition}\label{prop:AL-sublevel-set-bounded}
Assume that Assumption~\ref{ass:problem-setup} holds. For every
$q\in\R^m$, $\beta>0$, and $\kappa\in\R$, the set
\begin{equation}\label{eq:general-AL-sublevel-set}
\left\{X\in\mathbb S_+^n:
\mathcal L_\beta(X;q)\le\kappa
\right\}
\end{equation}
is bounded.
\end{proposition}

\begin{proof}
The proof can be found in 
\citep[Theorem~8.4, p.~64]{rockafellar1970convex}
(see also
\citep[Proposition~A.2.2.3, p.~40]{hiriarturruty2001fundamentals}).
\end{proof}

For every $\tau>0$, let
\begin{equation}\label{eq:Delta_tau_def}
\Delta_\tau
\ \defi\
\left\{X\in\mathbb S_+^n:\ \operatorname{tr}(X)\le\tau\right\}.
\end{equation}
For every $G\in\mathbb S^n$, the spectral decomposition of $G$ gives
\begin{equation}\label{eq:LMO_value_on_Delta_taut}
\min_{X\in\Delta_\tau}\inner{G}{X}
\ =\
\tau\min\left\{\lambda_{\min}(G),0\right\}.
\end{equation}

\begin{proposition}\label{prop:trace_ball_gap_residual_bounds}
Let $\tau>0$, $X\in\mathbb S_+^n$, and $G\in\mathbb S^n$, and define
\begin{equation}\label{eq:trace_ball_gap_def}
\epsilon_\tau(X;G)
\ \defi\
\max_{Y\in\Delta_\tau}\inner{G}{X-Y}.
\end{equation}
Then the following statements hold:
\begin{itemize}
\item[(a)] We have
\begin{equation}\label{eq:trace_ball_gap_upper_bound}
\epsilon_\tau(X;G)
\ \le\
|\inner{G}{X}|
+ \tau\max\{\lambda_{\max}(-G),0\}.
\end{equation}

\item[(b)] If $X\in\Delta_{\tau/2}$, then
\begin{equation}\label{eq:trace_ball_gap_residual_bounds}
3 \epsilon_\tau(X;G)
\ \ge\
|\inner{G}{X}|
+ \tau\max\{\lambda_{\max}(-G),0\}.
\end{equation}
\end{itemize}
\end{proposition}

\begin{proof}
a) Relations \eqref{eq:LMO_value_on_Delta_taut} and \eqref{eq:trace_ball_gap_def} imply that
\begin{equation}\label{eq:trace_ball_gap_identity}
\epsilon_\tau(X;G)
=
\inner{G}{X}
-\tau\min\{\lambda_{\min}(G),0\} = \inner{G}{X} + \tau\max\{\lambda_{\max}(-G),0\}.
\end{equation}
and hence that \eqref{eq:trace_ball_gap_upper_bound} holds.

b) Let $X\in\Delta_{\tau/2}$ be given. Since $0$ and $2X$ are both in $\Delta_\tau$ due to the definition of $\Delta_{\tau}$ in \eqref{eq:Delta_tau_def}, it follows 
from \eqref{eq:trace_ball_gap_def} with $Y=0$ and $Y=2X$ that
\begin{equation}\label{eq:trace_ball_gap_inner_bound}
\epsilon_\tau(X;G)
\ge
\max\big\{
\inner{G}{X-0},
\inner{G}{X-2X}
\big\}
=
|\inner{G}{X}|.
\end{equation}
Combining \eqref{eq:trace_ball_gap_identity} and
\eqref{eq:trace_ball_gap_inner_bound}, we obtain
\[
\tau\max\{\lambda_{\max}(-G),0\}
=
\epsilon_\tau(X;G)-\inner{G}{X}
\le
\epsilon_\tau(X;G)+|\inner{G}{X}|
\le
2\epsilon_\tau(X;G).
\]
Adding this inequality to \eqref{eq:trace_ball_gap_inner_bound} proves
\eqref{eq:trace_ball_gap_residual_bounds}.
\end{proof}

\subsection{\texorpdfstring{$X$}{X}-Space Relaxed Unbounded Frank--Wolfe Material}\label{sec:RUFW}
Given an initial point $\tilde X_0\in \mathbb S_{+}^{n}$, this subsection presents the Relaxed Unbounded Frank-Wolfe (RUFW) framework for solving
\begin{equation}\label{eq:RUFW_problem}
\phi_*
\ \defi\
\min\{\phi(X):\ X\in\mathbb S_+^n\},
\end{equation}
where the function $\phi:\mathbb S_+^n\to\mathbb R$ is assumed to satisfy the following conditions:

\begin{itemize}
\item[\textbf{(A1)}] $\phi$ is convex and differentiable on $\mathbb S_{+}^{n}$, and its gradient is $L$-Lipschitz continuous for some $L>0$, i.e.,
\[
\|\nabla \phi(X')-\nabla \phi(X)\|\leq L\|X'-X\|\qquad \forall X,X'\in \mathbb S_{+}^{n};
\]

\item[\textbf{(A2)}] the sublevel set $\mathcal S(\tilde X_0;\phi)$ is bounded, where $\mathcal S(\cdot;\cdot)$ is defined as in \eqref{SublevelDef} with $\Omega=\mathbb S_+^n$.
\end{itemize}

Assumption \textbf{(A1)} implies that
\begin{equation}\label{eq:RUFW_factor_curvature}
\phi(X')
\ \le\
\phi(X)+\inner{\nabla\phi(X)}{X'-X}
+\frac{L}{2}\|X'-X\|^2
\qquad
\forall\,X,X'\in\mathbb{S}_{+}^n.
\end{equation}
Here, the inner product is the Froebnius one and the norm is the one induced by the Frobenius inner product, i.e., $\inner{U}{V}=\textrm{tr}(U^{T}V)$ and $\|U\|^2=\textrm{tr}(U^TU)$.

Under \textbf{(A1)}-\textbf{(A2)}, it is also easy to see the minimum in~\eqref{eq:RUFW_problem} is attained since \textbf{(A1)}--\textbf{(A2)} imply that $\phi$ is a continuous coercive function.
Likewise, these assumptions imply that step 4 in the RUFW framework below is well-defined. We now describe the type of solution that RUFW aims to find.\\

\noindent{\bf Problem:}
Given a function $\phi$ and initial point $\tilde X_0\in \mathbb S_{+}^{n}$ satisfying \textbf{(A1)--(A2)} and a tolerance pair $(\epsilon,\delta)\in \mathbb R^{2}_{++}$, RUFW aims to find a point $\bar X\in \mathbb S_{+}^{n}$ satisfying
\begin{align}
\lambda_{\min}\left(\nabla \phi(\bar X)\right)&\ge -\epsilon,
\qquad
\left|\left\langle \nabla \phi(\bar X),\bar X\right\rangle\right|
\le \delta \label{SolutionTypeRUFW};\\
&\phi(\bar X)\leq \phi(\tilde X_0)\label{descentRUFW}.
\end{align}

The formal description of the relaxed unbounded Frank--Wolfe (RUFW) framework is given below.

\vspace{0.75em}
\noindent
\begin{minipage}[t]{1\columnwidth}
\rule[0.5ex]{1\columnwidth}{1pt}

\noindent \textbf{Relaxed Unbounded Frank--Wolfe (RUFW) Framework}

\noindent \rule[0.5ex]{1\columnwidth}{1pt}
\end{minipage}

\noindent \textbf{Input}: tolerance pair $(\epsilon,\delta)\in\R_{++}^2$, function $\phi:\mathbb S_{+}^n\to\mathbb R$, and initial point $\tilde X_0\in\mathbb S_+^n$.

\smallskip
\noindent \textbf{Step 0.} Set $k\gets 1$.

\smallskip
\noindent \textbf{Step 1.} Choose $X_k\in\mathbb S_+^n$ such that
\begin{equation}\label{eq:RUFW_correction_step}
\phi(X_k)\ \le\ \phi(\tilde X_{k-1}).
\end{equation}

\smallskip
\noindent \textbf{Step 2.} Set $G_k\gets \nabla\phi(X_k)$. Compute a minimum eigenpair $(\lambda^{\min}_{k}(G_k),v^{\min}_{k}(G_k))$ of $G_k$ with $\|v_k^{\min}(G_k)\|=1$, and set
\begin{align}
y^{F}_k
&\defi
\begin{cases}
v^{\min}_{k}(G_k), & \text{if }\lambda^{\min}_{k}(G_k)<0,\\
0, & \text{if }\lambda^{\min}_{k}(G_k)\ge0,
\end{cases}
\qquad
X_k^F\ \defi\ y^{F}_k(y^{F}_k)^\top;
\label{eq:RUFW_oracle_step}
\end{align}

\smallskip
\noindent \textbf{Step 3.} If
\begin{equation}\label{eq:RUFW_stop_test}
\lambda^{\min}_{k}(G_k)\ \ge\ -\epsilon,
\qquad
\left|\inner{G_k}{X_k}\right|\ \le\ \delta,
\end{equation}
then set $\bar X\gets X_k$, output $\bar X$, and {\bf stop}.

\smallskip
\noindent \textbf{Step 4.} Compute
\begin{equation}\label{eqn:alpha-linesearch-appendix}
(\alpha_k,\omega_k)
\ \in\
\argmin\left\{
\phi(\alpha X_k+\omega X_k^F)
:\
\alpha\ge0,\ \omega\ge0
\right\}
\end{equation}
and set
\begin{equation}\label{eqn:fw-short-step-appendix}
\tilde X_k
\ \defi\
\alpha_{k}X_k+\omega_{k}X_k^F;
\end{equation}

\smallskip
\noindent \textbf{Step 5.} Set $k\gets k+1$ and return to Step~\textbf{1}.

\noindent\rule[0.5ex]{1\columnwidth}{1pt}
\vspace{0.75em}

Step~1 allows any choice of $X_k\in\mathbb S_+^n$ whose objective value does not exceed that of the current iterate $\tilde X_{k-1}$. Step~2 evaluates the gradient $G_k=\nabla\phi(X_k)$ and computes a minimum eigenpair of $G_k$ to construct $X_k^F$ and check the stopping conditions \eqref{eq:RUFW_stop_test} in Step~3. If these conditions are not satisfied, Step~4 determines the coefficients $(\alpha_k,\omega_k)$ by solving the two-dimensional optimization problem \eqref{eqn:alpha-linesearch-appendix}, following an approach similar to that of Laue \cite{laue2012hybrid}. The next iterate $\tilde X_k$ is then formed as the conic combination of $X_k$ and $X_k^F$ given in \eqref{eqn:fw-short-step-appendix}. Since $(\alpha,\omega)=(1,0)$ is feasible in \eqref{eqn:alpha-linesearch-appendix}, Step~4 does not increase the objective value. Together with \eqref{eq:RUFW_correction_step}, this ensures that $\phi(X_k)\le\phi(\tilde X_0)$ for every $k$. Thus, whenever RUFW terminates in Step~3, its output $\bar X$ satisfies both \eqref{SolutionTypeRUFW} and \eqref{descentRUFW}.

The following proposition states the main convergence
result of the RUFW framework. Its proof is given in the following subsubsection.
In contrast to classical Frank--Wolfe complexity bounds, the bound in
\eqref{eq:RUFW_two_tolerance_complexity} below does not depend on the initial
objective gap $\phi(\tilde X_0)-\phi_*$. The initial point $\tilde X_0$ enters
the bound only through the radius $\bar\tau$ of its sublevel set.

\begin{proposition}\label{Prop:RUFW}
Assume that \textbf{(A1)}--\textbf{(A2)} hold
and define
\begin{equation}\label{eq:RUFW_level_set}
\bar\tau
\ \defi\
\max\left\{1,\
2\sup \left\{\operatorname{tr}(X):
X \in
\mathbb{S}^{n}_{+} , \ \phi(X) \le \phi(\tilde X_0) \right \}\right\},
\qquad
\epsilon_{k}(\bar \tau)\ \defi\ \max_{X\in\Delta_{\bar \tau}}\inner{G_k}{X_k-X}.
\end{equation}
Then, the following statements hold for RUFW:
\begin{itemize}
    \item[(a)] for every iteration index $k\geq 1$,
    it holds that $\phi(X_k)\leq \phi(\tilde X_{0})$;
    \item[(b)] RUFW terminates with
    a point $\bar X$ satisfying
    \eqref{SolutionTypeRUFW}
    and \eqref{descentRUFW}
    in at most \begin{equation}\label{eq:RUFW_two_tolerance_complexity}
    \left\lfloor
    2+12\bar\tau^2L\left(\frac{1}{\delta}+\frac{2}{\epsilon}\right)
    \right\rfloor
    \end{equation}
    iterations.
\end{itemize}
\end{proposition}

\subsubsection{Proof of Proposition~\ref{Prop:RUFW}}
This subsubsection proves Proposition~\ref{Prop:RUFW}.
Recall the trace-bounded spectrahedron $\Delta_\tau$ and its linear-oracle
value from \eqref{eq:Delta_tau_def} and
\eqref{eq:LMO_value_on_Delta_taut}, respectively. Let $X_*$ be a minimizer of
\eqref{eq:RUFW_problem}. Since $\phi(X_*)\le\phi(\tilde X_0)$, it follows from
\textbf{(A2)} and \eqref{eq:RUFW_level_set} that
$X_*\in\Delta_{\bar\tau/2}\subseteq\Delta_{\bar\tau}$. We now begin with a lemma describing key properties of the RUFW iterates.

\begin{lemma}\label{lem:RUFW_factor_monotonicity}
For every iteration index $k\ge1$, the following statements hold:
\begin{itemize}
\item[(a)] it holds
\begin{equation}\label{eq:RUFW_factor_monotonicity}
\phi(\tilde X_k) \leq \phi(X_k) \le \phi(\tilde X_{k-1}),
\qquad
\operatorname{tr}(X_k)\le \bar \tau/2,
\qquad
\operatorname{tr}(\tilde X_k)\le \bar \tau/2;
\end{equation}

\item[(b)] the FW gap $\epsilon_k(\bar\tau)$ satisfies
\begin{equation}\label{First impt rel}
|\inner{G_k}{X_k}|
\ \le\
\epsilon_k(\bar\tau),
\qquad
\lambda_k^{\min}(G_k)
\ \ge\
-\frac{2\epsilon_k(\bar\tau)}{\bar\tau}.
\end{equation}
\end{itemize}
\end{lemma}

\begin{proof}
(a) Relation \eqref{eq:RUFW_correction_step} gives
$\phi(X_k)\leq\phi(\tilde X_{k-1})$. Moreover, the choice of
$(\alpha_k,\omega_k)$ in \eqref{eqn:alpha-linesearch-appendix} and the update
\eqref{eqn:fw-short-step-appendix} give
$\phi(\tilde X_k)\leq\phi(X_k)$. Applying these inequalities recursively
yields
\[
\phi(X_k)\leq\phi(\tilde X_0),
\qquad
\phi(\tilde X_k)\leq\phi(\tilde X_0).
\]
The trace bounds now follow from the definition of $\bar\tau$ in
\eqref{eq:RUFW_level_set}.

(b) By part~(a) and Step~1, $X_k\in\Delta_{\bar\tau/2}$. Moreover, the
definition of the FW gap and \eqref{eq:trace_ball_gap_def} imply that
\[
\epsilon_k(\bar\tau)
=
\max_{X\in\Delta_{\bar\tau}}
\inner{G_k}{X_k-X}
=
\epsilon_{\bar\tau}(X_k;G_k).
\]
Hence, the argument proving
Proposition~\ref{prop:trace_ball_gap_residual_bounds}\emph{(b)},
applied with $(\tau,X,G)=(\bar\tau,X_k,G_k)$, and the identity
$\lambda_k^{\min}(G_k)=-\lambda_{\max}(-G_k)$ yield
\eqref{First impt rel}.
\end{proof}

The following lemma characterizes key properties of $\epsilon_k(\bar \tau)$ and $X_k^{F}$.

\begin{lemma}
For every iteration index $k\geq 1$, the following relations hold:
\begin{equation}\label{eq:RUFW_scaled_oracle}
\bar \tau X_k^F
\ \in\
\argmin_{X\in\Delta_{\bar \tau}}\inner{G_k}{X},
\qquad
\epsilon_k(\bar \tau)
\ =\
\inner{G_k}{X_k-\bar \tau X_k^F},
\end{equation}
\begin{equation}\label{epsilonlowerbound}
 \epsilon_k(\bar \tau)
\ \ge\
\phi(X_k)-\phi_{*}.
\end{equation}
\end{lemma}
\begin{proof}
It is easy to see from \eqref{eq:RUFW_oracle_step}, the unit eigenvector convention, and \eqref{eq:LMO_value_on_Delta_taut} with $\tau=1$ that $X_k^F\in\argmin_{X\in\Delta_1}\inner{G_k}{X}$: if $\lambda_k^{\min}<0$, then $X_k^F=v_k^{\min}(G_k)(v_k^{\min}(G_k))^\top$ attains $\lambda_k^{\min}(G_k)$, while if $\lambda_k^{\min}(G_k)\ge0$, then $X_k^F=0$ attains the zero minimum. Scaling gives the first relation in~\eqref{eq:RUFW_scaled_oracle}. The second relation in~\eqref{eq:RUFW_scaled_oracle} follows immediately from the first relation and the definition of $\epsilon_k(\bar\tau)$ in~\eqref{eq:RUFW_level_set}. Finally, convexity, $X_*\in\Delta_{\bar\tau}$, and \eqref{eq:RUFW_level_set} give $\phi_*=\phi(X_*)\ge\phi(X_k)+\inner{G_k}{X_*-X_k}=\phi(X_k)-\inner{G_k}{X_k-X_*}\ge\phi(X_k)-\epsilon_k(\bar\tau)$, which is \eqref{epsilonlowerbound}.
\end{proof}

The following proposition gives the one-step estimate on which the analysis
is based. It uses the fact that, for every $\alpha\in[0,1]$, the conic search
in Step~4 includes the point $(1-\alpha)X_k+\alpha\bar\tau X_k^F$.

\begin{proposition}\label{prop:RUFW_one_step_gap_descent}
For every iteration index $k\ge1$ and every $\alpha\in[0,1]$, we have
\begin{equation}\label{upperboundepsilon}
\phi(\tilde X_k)
\ \le\
\phi(X_k)-\alpha\,\epsilon_k(\bar \tau)
+\frac{3\bar\tau^2L}{2}\,\alpha^2,
\end{equation}
where $\epsilon_k(\bar \tau)$ is as in \eqref{eq:RUFW_level_set}.
\end{proposition}

\begin{proof}
Let $\alpha\in[0,1]$ be given. Since the pair $(1-\alpha,\alpha\bar\tau)$ is feasible in~\eqref{eqn:alpha-linesearch-appendix}, relations \eqref{eqn:alpha-linesearch-appendix}--\eqref{eqn:fw-short-step-appendix} imply
\(\phi(\tilde X_{k})\le\phi((1-\alpha)X_k+ \alpha\bar \tau X_k^{F})=\phi(X_k+\alpha(\bar \tau X_k^{F}-X_k))\).
Observe also that the first relation in \eqref{eq:RUFW_scaled_oracle} and
\eqref{eq:RUFW_factor_monotonicity} imply that
\begin{equation}\label{boundXk}
   \|\bar \tau X_k^{F}-X_k\|\leq  \|\bar \tau X_k^{F}\|+\|X_k\|\leq \frac{3\bar \tau}{2}.
\end{equation}
It then follows from the above observations, \eqref{boundXk}, $G_k=\nabla \phi(X_k)$, and \eqref{eq:RUFW_factor_curvature}
with $X'=X_k+\alpha(\bar \tau X_k^{F}-X_k)\in\mathbb S_+^n$ that
\begin{align*}
    \phi(\tilde X_k)&\leq \phi\left(X_k+\alpha(\bar \tau X_k^{F}-X_k) \right)\leq
    \phi(X_k)+\alpha\inner{\nabla \phi(X_k)
    }{\bar \tau X_k^{F}-X_k}+\frac{L\alpha^2}{2}\|\bar \tau X_k^{F}-X_k\|^2\\
&\overset{\eqref{boundXk}}{\leq} \phi(X_k)+\alpha\inner{G_k
    }{\bar \tau X_k^{F}-X_k}+\frac{9L\alpha^2\bar \tau^2}{8}
\overset{\eqref{eq:RUFW_scaled_oracle}}{=}\phi(X_k)-\alpha \epsilon_k(\bar \tau)+\frac{9L\alpha^2\bar \tau^2}{8}.
\end{align*}
Relation \eqref{upperboundepsilon} then follows immediately from the above inequalities and the fact that $9/8\le3/2$.
\end{proof}

The next lemma establishes $\mathcal O(1/k)$ bounds on the objective gap and the minimum Frank--Wolfe gap in \eqref{epsilonlowerbound}, both independent of the initial objective gap $\phi(\tilde X_0)-\phi_*$. By the two relations in \eqref{First impt rel}, the latter bound also controls the quantities appearing in the stopping conditions \eqref{eq:RUFW_stop_test}.

\begin{lemma}\label{lem:RUFW_gap_complexity}
Let $k\ge1$ be an iteration index such that Step~4 is performed at
iterations $1,\ldots,k$. Then the following statements hold:
\begin{itemize}
\item[(a)] we have
\begin{equation}\label{eq:RUFW_objective_gap_rate}
\phi(\tilde X_k)-\phi_*
\ \le\
\frac{6\bar\tau^2L}{k+3};
\end{equation}
\item[(b)] if $k\ge2$, then
\begin{equation}\label{eq:RUFW_gap_complexity}
\min_{1\le j\le k}\epsilon_j(\bar\tau)
\ \le\
\frac{12\bar\tau^2L}{k}.
\end{equation}
\end{itemize}
\end{lemma}

\begin{proof}
(a) Let $j\in\{1,\ldots,k\}$ and $\alpha\in[0,1]$ be given. Using
relations \eqref{upperboundepsilon} and \eqref{epsilonlowerbound}, the fact
that $1-\alpha\ge0$, and the second inequality in
\eqref{eq:RUFW_factor_monotonicity}, we obtain
\begin{align}
\phi(\tilde X_j)-\phi_*
&\overset{\eqref{upperboundepsilon}}{\le}
\phi(X_j)-\phi_*-\alpha\,\epsilon_j(\bar\tau)
+\frac{3\bar\tau^2L}{2}\alpha^2
\overset{\eqref{epsilonlowerbound}}{\le}
(1-\alpha)\left[\phi(X_j)-\phi_*\right]
+\frac{3\bar\tau^2L}{2}\alpha^2
\notag\\
&\overset{\eqref{eq:RUFW_factor_monotonicity}}{\le}
(1-\alpha)\left[\phi(\tilde X_{j-1})-\phi_*\right]
+\frac{3\bar\tau^2L}{2}\alpha^2.
\label{eq:RUFW_gap_recursion}
\end{align}
We now prove \eqref{eq:RUFW_objective_gap_rate} by induction.
Relation \eqref{eq:RUFW_gap_recursion} with $j=1$ and $\alpha=1$ gives
\[
\phi(\tilde X_1)-\phi_*
\ \le\
\frac{3\bar\tau^2L}{2}
\ =\
\frac{6\bar\tau^2L}{1+3},
\]
which is \eqref{eq:RUFW_objective_gap_rate} with $k=1$. Observe that this
bound does not depend on $\tilde X_0$ since the coefficient $1-\alpha$
multiplying $\phi(\tilde X_0)-\phi_*$ in \eqref{eq:RUFW_gap_recursion}
vanishes. Assume now that $2\le j\le k$ and that
$\phi(\tilde X_{j-1})-\phi_*\le 6\bar\tau^2L/(j+2)$. Relation
\eqref{eq:RUFW_gap_recursion} with $\alpha=2/(j+2)\in[0,1]$ and this
assumption then imply that
\[
\phi(\tilde X_j)-\phi_*
\ \le\
\frac{j}{j+2}\cdot\frac{6\bar\tau^2L}{j+2}
+\frac{3\bar\tau^2L}{2}\cdot\frac{4}{(j+2)^2}
\ =\
\frac{6\bar\tau^2L(j+1)}{(j+2)^2}
\ \le\
\frac{6\bar\tau^2L}{j+3},
\]
where the last inequality is due to the fact that $(j+1)(j+3)\le(j+2)^2$.
Relation \eqref{eq:RUFW_objective_gap_rate} then follows by induction.

(b) Let $k\ge2$ be given and note that $\lfloor k/2\rfloor\ge1$. Let
$\alpha\in(0,1]$ and $j\in\{1,\ldots,k\}$ be given. Relation
\eqref{upperboundepsilon} and the second inequality in
\eqref{eq:RUFW_factor_monotonicity} imply that
\begin{equation}\label{eq:RUFW_gap_telescoping_step}
\alpha\,\epsilon_j(\bar\tau)
\ \overset{\eqref{upperboundepsilon}}{\le}\
\phi(X_j)-\phi(\tilde X_j)+\frac{3\bar\tau^2L}{2}\alpha^2
\ \overset{\eqref{eq:RUFW_factor_monotonicity}}{\le}\
\phi(\tilde X_{j-1})-\phi(\tilde X_j)+\frac{3\bar\tau^2L}{2}\alpha^2.
\end{equation}
Summing \eqref{eq:RUFW_gap_telescoping_step} over
$j=\lfloor k/2\rfloor+1,\ldots,k$, a range with $\lceil k/2\rceil$
indices, and using the fact that $\phi(\tilde X_k)\ge\phi_*$ and
\eqref{eq:RUFW_objective_gap_rate} with $k$ replaced by $\lfloor k/2\rfloor$, we obtain
\begin{align*}
\alpha\left\lceil\frac k2\right\rceil
\min_{1\le j\le k}\epsilon_j(\bar\tau)
&\le
\alpha\sum_{j=\lfloor k/2\rfloor+1}^{k}\epsilon_j(\bar\tau)
\overset{\eqref{eq:RUFW_gap_telescoping_step}}{\le}
\phi(\tilde X_{\lfloor k/2\rfloor})-\phi(\tilde X_k)
+\left\lceil\frac k2\right\rceil\frac{3\bar\tau^2L}{2}\alpha^2
\\
&\overset{\eqref{eq:RUFW_objective_gap_rate}}{\le}
\frac{6\bar\tau^2L}{\lfloor k/2\rfloor+3}
+\left\lceil\frac k2\right\rceil\frac{3\bar\tau^2L}{2}\alpha^2.
\end{align*}
Taking $\alpha=2/(\lfloor k/2\rfloor+3)\in(0,1]$ in the above inequality and
dividing the result by $\alpha\lceil k/2\rceil$, we conclude that
\[
\min_{1\le j\le k}\epsilon_j(\bar\tau)
\ \le\
\frac{3\bar\tau^2L}{\lceil k/2\rceil}
+\frac{3\bar\tau^2L}{\lfloor k/2\rfloor+3}
\ \le\
\frac{6\bar\tau^2L}{k}
+\frac{6\bar\tau^2L}{k}
\ =\
\frac{12\bar\tau^2L}{k},
\]
where the last inequality is due to the fact that both $\lceil k/2\rceil$ and
$\lfloor k/2\rfloor+3$ are bounded below by $k/2$.
\end{proof}

We are now ready to prove Proposition~\ref{Prop:RUFW}.

\begin{proof}[Proof of Proposition~\ref{Prop:RUFW}]
(a) Applying the first relation in
\eqref{eq:RUFW_factor_monotonicity} recursively implies that
$\phi(X_k)\le\phi(\tilde X_0)$ for every iteration index $k\geq 1$.

(b)
Let $k$ denote the integer in \eqref{eq:RUFW_two_tolerance_complexity}.
Then $k\ge2$ and $k>12\bar\tau^2L(1/\delta+2/\epsilon)$, and hence
\begin{equation}\label{eq:FW_q_threshold}
\frac{12\bar\tau^2L}{k}
\ <\
\left(\frac{1}{\delta}+\frac{2}{\epsilon}\right)^{-1}
\ \le\
\min\left\{\delta,\frac{\epsilon}{2}\right\}.
\end{equation}
Suppose now for contradiction, that RUFW does not terminate at any iteration index $j\le k$, so that it performs Step~4 at iterations $1,\ldots,k$ and the stopping test \eqref{eq:RUFW_stop_test} fails at every iteration index $j\le k$. Both relations in \eqref{First impt rel}, the fact that $\bar\tau\ge1$, and
\eqref{eq:FW_q_threshold}
imply that there exists an iteration index $i\in\{1,\ldots,k\}$ such that
\begin{align*}
&\left|\inner{G_i}{X_i}\right|\overset{\eqref{First impt rel}}{\leq} \epsilon_i(\bar \tau) \ \overset{\eqref{eq:RUFW_gap_complexity}}{\le}\ \frac{12\bar\tau^2L}{k}\ \overset{\eqref{eq:FW_q_threshold}}{<}\ \delta,\\
&\lambda_i^{\min}(G_i) \overset{\eqref{First impt rel}}{\geq} \frac{-2\epsilon_i(\bar \tau)}{\bar \tau}\ \overset{\eqref{eq:RUFW_gap_complexity}}{\ge}\ -\frac{24\bar\tau L}{k}\ \ge\ -2\cdot\frac{12\bar\tau^2L}{k}\ \overset{\eqref{eq:FW_q_threshold}}{>}\ -\epsilon.
\end{align*}
Thus, the stopping test \eqref{eq:RUFW_stop_test} holds at the iteration index
$i\le k$, which contradicts the assumption that RUFW does not terminate at any iteration index $j\le k$. Hence, RUFW terminates at some iteration index $j\le k$ with output $\bar X=X_j$, which satisfies \eqref{SolutionTypeRUFW} in view of the stopping test \eqref{eq:RUFW_stop_test} and satisfies \eqref{descentRUFW} in view of part~(a).
\end{proof}

\section{\texorpdfstring{$U$}{U}-Space and Factorized Machinery}
\label{app:factorized-machinery}

This appendix collects the material used to implement the descent step
\eqref{eq:RUFW_correction_step} of the RUFW framework of
Appendix~\ref{sec:RUFW} in the factorized $U$-space. It is divided into two
parts. Subsection~\ref{sec:adap-fista} presents the ADAP-FISTA method while
Subsection~\ref{app:aipp} presents the ADAP-AIPP method, whose subproblems are
solved by ADAP-FISTA.

\subsection{ADAP-FISTA}\label{sec:adap-fista}
Let $\mathbb E$ be a finite-dimensional inner product space and let $\Omega\subseteq \mathbb E$ be a nonempty closed convex set.
This section presents the ADAP-FISTA method for the problem
\begin{equation}\label{NonconvexOptProb}
\min\{\psi(U): U \in \Omega\}
\end{equation}
where $\psi$ and $\Omega$ are assumed to satisfy the following assumptions:

\begin{itemize}
\item[\textbf{(B1)}]
$\psi$ is a real-valued function that is differentiable on $\mathbb E$ and
whose gradient is locally Lipschitz continuous in the sense of
\eqref{eq:locally_lipschitz_gradient};
\item[\textbf{(B2)}] $\Omega$ is closed, possibly unbounded, convex set that admits an easily computable projection; 
\item[\textbf{(B3)}] 
a point $X_0\in \Omega$ is given such that the sublevel set
$\mathcal S(X_0;\psi)$ has finite diameter $D(X_0;\psi)$,
where $\mathcal S(\cdot;\cdot)$ and $D(\cdot;\cdot)$
are as in \eqref{SublevelDef} and \eqref{DiameterDef}, respectively.
\end{itemize}

% \paragraph{Remark.}
% Using the optimality of $y^+(M)$ with the test point $u=\tilde x$ and
% Cauchy--Schwarz gives the weaker estimate
% $M\|y^+(M)-\tilde x\|\le 2\|g(\tilde x)\|$. We do not use this weaker bound in
% the sequel because the projection argument above gives the sharper threshold
% appearing in \eqref{def:M0}.

\vgap

\noindent{\bf Problem:}
Given $(\psi, \Omega, X_0)$ satisfying \textbf{(B1)-(B3)}, the goal of ADAP-FISTA is to find a pair $(\bar Y,\bar V)\in \Omega\times\mathbb E$ such that
\begin{equation}\label{eq:adap_fista_goal}
\bar V \in \nabla \psi(\bar Y)+\mathcal{N}_{\Omega}(\bar Y),\qquad \|\bar V\|\leq \sigma\|\bar Y-X_0\|,\qquad \psi(\bar Y)\leq \psi(X_0).
\end{equation}
The formal description of ADAP-FISTA is now given below.

\noindent\begin{minipage}[t]{1\columnwidth}%
\rule[0.5ex]{1\columnwidth}{1pt}

\noindent \textbf{ADAP-FISTA Method}

\noindent \rule[0.5ex]{1\columnwidth}{1pt}%
\end{minipage}
\noindent \textbf{Universal Parameters}: $\sigma>0$ and $\chi \in (0,1)$.

\noindent \textbf{Input}: initial point $\bar Y_0 \in \Omega$,  pair $\mu>0$, $L_0>\mu$, and a function $\psi$ and a set $\Omega$ as in \textbf{(B1)}-\textbf{(B3)}
\begin{itemize}
\item[{\bf 0.}] set $Y_0=X_0=\bar Y_0$, $ A_0=0 $, $\tau_0=1$, and $i=1$;
\item[{\bf 1.}] set $L=L_{i-1}$;
\item[{\bf 2.}]	compute
\begin{equation}\label{def:ak-sfista1Tenative}
a=\frac{\tau_{i-1}+\sqrt{\tau_{i-1}^2+4\tau_{i-1} A_{i-1}(L-\mu)}}{2(L-\mu)}, \quad \tilde X=\frac{A_{i-1}Y_{i-1}+a X_{i-1}}{A_{i-1}+a},
\end{equation}
\begin{equation}
Y:=\underset{U\in \Omega}\argmin\left\lbrace \ell_{\psi}(U;\tilde X) + \frac{L}{2}\|U-\tilde X\|^2\right\rbrace.
\label{eq:ynext-sfista1Tenative}
\end{equation}
% \jacob{I may be missing the definition, but could we recall $\ell_{\psi}(\cdot;\tilde X)$ here? I do not see it introduced in this section.}
% \arnesh{I have now added the definition of $\ell (\cdot;\cdot)$ in \eqref{eq:defell} and I just decided to recall its definition after \eqref{eq:ynext-sfista1} instead of the algorithm so the algorithm isn't too crowded.}

If the inequality
\begin{equation}\label{ineq check}
\ell_{\psi}(Y;\tilde X)+\frac{(1-\chi) L}{2}\|Y-\tilde X\|^2\geq \psi(Y)
\end{equation}
holds then set
\begin{equation}\label{setIterates}
    a_{i-1}=a, \quad \tilde X_{i-1}=\tilde X, \quad Y_i=Y, \quad L_i=L
\end{equation}
and go to step~3;
else set $L \leftarrow 2 L $ and repeat step~2;
\item[{\bf 3.}] compute
\begin{align}
A_{i}&=A_{i-1}+a_{i-1}, \quad \tau_{i}= \tau_{i-1} + a_{i-1}\mu,  \label{eq:taunext-sfista1} \\
\quad X_{i}&= \frac{1}{\tau_{i}} \left[\mu a_{i-1}Y_{i} + \tau_{i-1} X_{i-1}-a_{i-1}(L_{i}-\mu)(\tilde X_{i-1}-Y_{i}) \right]; \notag
\end{align}
	
\item[{\bf 4.}]  if the following three inequalities hold
\begin{equation}\label{ineq: ineq 1}
\|Y_{i}-\bar Y_0\|^{2} \geq \chi A_{i}L_{i} \|Y_{i}-\tilde X_{i-1}\|^2, \quad  \psi(Y_{i})\leq \psi (\bar Y_0)
\end{equation}
\begin{equation}\label{ineq 2}
   \|Y_{i}-X_{i}\|\leq \|Y_{i}-\bar Y_0\|
\end{equation}
hold, then  go to step 5; otherwise, stop with {\bf failure};

\item[{\bf 5.}] compute 
\begin{equation}\label{def:uk}
V_{i}=\nabla \psi(Y_{i})-\nabla \psi(\tilde X_{i-1})+L_{i}(\tilde X_{i-1}-Y_{i});  \end{equation}
if the inequality
\begin{equation}\label{u sigma criteria}
\|V_{i}\| \leq \sigma \|Y_{i}-\bar Y_0\|
\end{equation}
holds then stop with {\bf success} and output $(\bar Y,\bar V,\bar L):=(Y_{i},V_{i}, L_{i})$; otherwise,
 $ i \leftarrow i+1 $ and go to step~1.
\end{itemize}
\noindent \rule[0.5ex]{1\columnwidth}{1pt}
It is easy to see that the update for $Y_i$ in \eqref{eq:ynext-sfista1} is equivalent to projecting the point $\tilde X_{i-1}-(1/L_i)\nabla \psi(\tilde X_{i-1})$ onto the set $\Omega$. ADAP-FISTA takes a parameter $\mu>0$ as input, even when we do not know in advance whether $\psi$ is $\mu$-strongly convex. The method therefore checks the inequalities in \eqref{ineq: ineq 1} and \eqref{ineq 2} at each iteration. If any of these inequalities is violated, the method stops with failure. When $\psi$ is $\mu$-strongly convex, these inequalities always hold, and ADAP-FISTA returns a triple satisfying \eqref{u sigma criteria} within a logarithmic number of iterations. To accommodate objectives with locally Lipschitz continuous gradients, we use different inequalities from those in earlier adaptive variants of FISTA~\cite{LinQihang,sujanani2025efficient}.

Before stating the main complexity result of ADAP-FISTA, we introduce the quantities which appear in the complexity result:
\begin{equation}\label{def:D0_G0_Lip0}
D_0\defi 6D(X_0;\psi)+1,\quad G_0:=\nabla (D_0/2;\psi), \quad ({\rm lip})_0:=[{\rm LIP}](D_0;\psi)
\end{equation}
where $D(\cdot; \cdot)$, $\nabla(\cdot; \cdot)$, and 
$[{\rm LIP}]$ are as in \eqref{DiameterDef}, \eqref{GradDef}, and \eqref{SmoothnessDef}, respectively.
% We also let $({\rm LIP})_0$ be any Lipschitz constant of
% $\nabla\psi$ on $B[D_0]$, i.e, 
% \begin{equation}\label{def:barL}
% \|\nabla\psi(X')-\nabla\psi(X'')\|
% \ \le\
% ({\rm LIP})_0\|X'-X''\|
% \qquad \forall X',X''\in B[D_0].
% \end{equation}
% and define {\bf Move to te main result}
% \begin{equation}\label{ND0 Def}
% N_0:=\max\left\{L_0, G_0,({\rm LIP})_0\right\}.
% \end{equation} 

The quantities above are finite. Indeed, \textbf{(B3)} implies that $D_0$ is finite. Since $\nabla\psi$ is continuous, it follows that $G_0<\infty$. Since
$\nabla\psi$ is locally Lipschitz, a finite-cover argument on the compact ball
$B[D_0]$ yields that the Lipschitz constant $({\rm lip})_0$ is finite.
The main proposition which states ADAP-FISTA's main complexity result is now given. Its proof is given in the next two subsections.

\begin{proposition}\label{prop:nest_complex1} 
% ({\bf Result translated from Proposition~B.2 of \cite{monteiro2024low}})
Let \begin{equation}\label{ND0 Def}
N_0:=\max\left\{L_0, G_0,({\rm lip})_0\right\}.
\end{equation} 
where $L_0$ is input to ADAP-FISTA and $G_0$ and $({\rm lip})_0$ are as in \eqref{def:D0_G0_Lip0}.
The following statements about ADAP-FISTA hold:
\begin{itemize}
\item[(a)]
% $L_0 = {\cal O}(({\rm lip})_0)$, 
it
always stops (with either success or failure) in at most 
\[\mathcal O_1\left(\sqrt{\frac{N_0}{\mu}}\log_2^+ \left(N_0\right)+ \log_2^+\left(\frac{N_0}{L_0}\right)\right)\] 
iterations/resolvent evaluations;
% where $({\rm lip})_0$ and $N_0$
% are as in \eqref{def:D0_G0_Lip0} and \eqref{ND0 Def}, respectively;
\item[(b)] if it stops successfully, then the pair $(\bar Y,\bar V)\in\Omega\times\mathbb E$ output by ADAP-FISTA satisfies the goal \eqref{eq:adap_fista_goal}, i.e.,
\begin{align}
&  \bar V \in \nabla \psi(\bar Y)+\mathcal{N}_{\Omega}(\bar Y), \quad \|\bar V\|\leq \sigma\|\bar Y-\bar Y_0\|,\label{u inclusion 1}\\
&\psi(\bar Y)\leq \psi(\bar Y_0), \label{secondOutput}
\end{align}
and the Lipschitz estimate $\bar L$ that it outputs satisfies
\begin{equation}\label{Lbound}
L_0\leq \bar L\leq \max\left\{L_0,12 G_0,\ \frac{2({\rm lip})_0}{1-\chi} \right\};
\end{equation}
% where $L_0$ is input to ADAP-FISTA and $G_0$ is as in \eqref{def:D0_G0_Lip0};
\item[(c)] if $\psi$ is $\mu$-strongly convex on $B[D_0]$, then ADAP-FISTA always terminates with success.
% and, in addition to \eqref{u inclusion 1}--\eqref{Lbound}, its output also satisfies
% \jacob{Would it be better to say ``$\mu$-strongly convex'' here, or to define ``$\mu$-convex'' explicitly?}
% \arnesh{Yes, you are right. I have fixed this now to say $\mu$-strongly convex.}
% \begin{equation}\label{KeyInequality}
% \psi (X_0) - \psi(\bar Y) \ge \inner{\bar V}{X_0-\bar Y}.
% \end{equation}
\end{itemize}
% where $D_0$ is as in \eqref{def:D0_G0_Lip0}.
\end{proposition}

We now make one brief remark about the above proposition. In general, the assumption in $(c)$ that $\psi$ is $\mu$-strongly convex on $B[D_0]$ is not directly verifiable since we do not know $D_0$. However, we will see in the next section that it will play an important role in guaranteeing that an inner loop of our other method ADAP-AIPP always terminates.
% {\bf Add a remark that the assumption of (c) is not directly verifiable since we do not know $D_0$}

\subsubsection{Proof of Proposition~\ref{prop:nest_complex1}(a)-(b)}
This subsection is dedicated to proving 
Proposition~\ref{prop:nest_complex1}(a)-(b).
We first define several quantities which arise in our technical analysis:
\begin{equation}\label{omega zeta}
\zeta:=({\rm lip})_0+\max\left\{L_0,12 G_0,\ \frac{2 ({\rm lip})_0}{1-\chi} \right\}, \quad Q:=  \sqrt{\frac{4}{\mu}} \sqrt{\max\left\{L_0,12G_0, \frac{2({\rm lip})_0}{1-\chi}\right\}}.
\end{equation}
Also note that since $\nabla \psi$ is $({\rm lip})_0$-Lipschitz on $B[D_0]$ it follows from 
that
\begin{equation}\label{adapFistaUpperCurvature}
\psi(X'')
\ \le\
\psi(X')+\inner{\nabla\psi(X')}{X''-X'}
+\frac{({\rm lip})_0}{2}\|X''-X'\|^2
\qquad
\forall\,X',X''\in B[D_0].
\end{equation}
We start off with a key lemma which characterizes properties of the output of the loop in step 2 of ADAP-FISTA.

\begin{lemma}
    Let $i\geq 1$ be an iteration index of ADAP-FISTA. If the line-search in the loop in step 2 terminates, then the following statements then hold:
\begin{itemize}
    \item[(a)]
\begin{equation}\label{def:ak-sfista1}
a_{i-1}=\frac{\tau_{i-1}+\sqrt{\tau_{i-1}^2+4\tau_{i-1} A_{i-1}(L_{i}-\mu)}}{2(L_{i}-\mu)}, \quad \tilde X_{i-1}=\frac{A_{i-1}Y_{i-1}+a_{i-1} X_{i-1}}{A_{i-1}+a_{i-1}},
\end{equation}
\begin{equation}
Y_{i}:=\underset{U\in \Omega}\argmin\left\lbrace \ell_{\psi}(U;\tilde X_{i-1}) + \frac{L_{i}}{2}\|U-\tilde X_{i-1}\|^2\right\rbrace;
\label{eq:ynext-sfista1}
\end{equation}
\item[(b)] it also holds that
\begin{align}
       &\tau_{i-1} = 1+A_{i-1}\mu, \quad \frac{\tau_{i-1}  A_{i}}{ a_{i-1}^2}=L_{i}-\mu\quad \label{tauproperty}\\
       & A_iL_i \geq \max \left\{\frac{i^2}{4} ,  \left(1+Q^{-1}\right)^{2(i-1)}
 \right\}\label{AiLiBound}.
    \end{align}
\end{itemize}

% \jacob{Could this be $\ell_{\psi}$ rather than $\ell_{\psi_s}$, and $\|U-\tilde X_{i-1}\|^2$ rather than $\|u-\tilde X_{i-1}\|^2$?}
% \arnesh{Yes, you are right. I have fixed these typos now.}
\end{lemma}
\begin{proof}
    (a) The proof of (a) is immediate from relations \eqref{def:ak-sfista1Tenative} and \eqref{eq:ynext-sfista1Tenative} and the way $a_{i-1}$, $\tilde X_{i-1}$, $Y_i$, and $L_i$ are set in \eqref{setIterates} at the end of the loop in step 2 of ADAP-FISTA.

    (b) The first equality in \eqref{tauproperty} follows directly from both of the relations in \eqref{eq:taunext-sfista1}. The second equality in \eqref{tauproperty} follows immediately from the definition of $a_{i-1}$ in \eqref{def:ak-sfista1} and the first relation in  \eqref{eq:taunext-sfista1}. The proof of the lower bound on $A_iL_i$ in \eqref{AiLiBound} follows from a near identical argument to the proof of Lemma A.4 in \cite{sujanani2023adaptive}.
\end{proof}

% \begin{equation}\label{Q}
% Q:=  2\sqrt{2} \sqrt{\frac{\max\{L_0,\omega {({\rm lip})_0}\}}{\mu}}.
% \end{equation}
We now start off with a key technical lemma which is instrumental for showing that the line-search performed in step 2 of ADAP-FISTA always terminates.

\begin{lemma}\label{lem:stepTwoTermination}
Let $\theta>0$, $D_0$ be as in \eqref{def:D0_G0_Lip0}, and $\hat X\in B[D_0/2]$. For $\hat L>0$, define
\begin{equation}\label{ProjYM}
Y^+_{\hat L} \defi \arg\min_{U\in\Omega}
\left\{\ell_{\psi}(U; \hat X) + \frac{\hat L}{2}\|U-\hat X\|^2 \right\}.
\end{equation}
If
\begin{equation}\label{assumpM}
\hat L \;\ge\; \max\left\{ 6\|\nabla \psi(\hat X)\|,\ \frac{({\rm lip})_0}{\theta} \right\}\geq\max\left\{ \frac{6\|\nabla \psi(\hat X)\|}{D_0},\ \frac{({\rm lip})_0}{\theta} \right\},
\end{equation}
then $Y^+_{\hat L} \in B[D_0]$ and
\begin{equation}\label{Upcurvature}
\ell_{\psi}(Y^+_{\hat L};\hat X)
+ 0.5\theta \hat L\|Y^+_{\hat L} - \hat X\|^2
\;\ge\;
\psi(Y^+_{\hat L}).
\end{equation}
\end{lemma}

\begin{proof}
It follows from the fact that $X_0\in \mathcal S(X_0;\psi)$, \textbf{(B3)}, and the definition of $D_0$ in \eqref{def:D0_G0_Lip0} that $\|X_0\|\leq D_0/6$.
It is easy to see from this fact, \eqref{ProjYM}, the nonexpansivity of the projection operator defined in \eqref{ProjOperator}, the assumption on $\hat L$ in \eqref{assumpM}, the fact that $D_0\geq 1$, and the fact that $\hat X \in B[D_0/2]$ that:
\begin{align*}
  \|Y_{\hat L}^{+}\|&\overset{\eqref{ProjYM}}{=} \left\|\mbox{proj}_{\Omega}\left(\hat X-\frac{1}{\hat L}\nabla \psi(\hat X)\right)-\mbox{proj}_{\Omega}(X_0)+\mbox{proj}_{\Omega}(X_0)\right\|\leq \left\|\hat X-\frac{1}{\hat L}\nabla \psi(\hat X)-X_0\right\|+\|X_0\|\\
  &\leq \|\hat X-X_0\|+\frac{1}{\hat L}\left\|\nabla \psi(\hat X)\right\|+\|X_0\|\overset{\eqref{assumpM}}{\leq} \frac{4D_0}{6}+\frac{D_0}{6}+\frac{D_0}{6}=D_0.
\end{align*}
Hence, $Y_{\hat L}^{+}\in B[D_0]$.
Since $\hat X$ and $Y_{\hat L}^{+}$ belong to $B(D_0)$, it follows from \eqref{adapFistaUpperCurvature} with $X''=Y_{\hat L}^{+}$ and $X'=\hat X$
that 
\[\psi(Y_{\hat L}^{+})\leq \ell_{\psi}(Y_{\hat L}^{+};\hat X)+\frac{({\rm lip})_0}{2}\|Y_{\hat L}^{+}-\hat X\|^2\]
which in view of the assumption on $\hat L$ in \eqref{assumpM}
implies the second conclusion of the lemma.
\end{proof}

The following lemma shows that the line-search in step 2 of ADAP-FISTA always terminates and provides a lower and upper bound on the sequence of Lipschitz estimates $\left\{L_i\right\}$.
\begin{lemma}\label{LineSearchLemma}
For any iteration index $i\geq1$ generated by ADAP-FISTA, the line-search in step 2 always terminates, the iterates $X_{i-1}$, $Y_{i-1}$, and $\tilde X_{i-1}$ all lie in $ B[D_0/2]$, and 
    \begin{equation}
        L_0\leq L_{i}\leq \max\left\{L_0,12 G_0,\ \frac{2 ({\rm lip})_0}{1-\chi} \right\}.\label{upper bound}
        \end{equation}
\end{lemma}
\begin{proof}
We proceed by induction. We do not include the argument for the $i=1$ case since it follows from the same argument as below and by using the fact that $X_0=Y_0\in B[D_0/6]$.

Assume then that the statement holds for an iteration index $i\geq 1$ generated by ADAP-FISTA. We show that these facts hold also for the $(i+1)$-st iteration index. Since the method generates the $(i+1)$-st iteration index, the inequalities checked in step 4 of the $i$-th iteration are satisfied, i.e.,
\begin{equation}\label{SublevelBound}
\psi(Y_{i})\leq \psi(X_0), \quad \|Y_{i}-X_{i}\|\leq \|Y_{i}-X_0\|. 
\end{equation}
Hence, $Y_i \in \mathcal S(X_0;\psi)$ and by \eqref{def:D0_G0_Lip0},
% \begin{equation}\label{eq:Y_norm_bound}
$\|Y_{i}\|\ \le D_0/6.$
% \end{equation}
Using this fact, \eqref{SublevelBound}, and $\|X_0\|\le D_0/6$, we obtain
\begin{equation}\label{eq:X_norm_bound}
\|X_{i}\|
\ \le\
\|Y_{i}\|+\|Y_{i}-X_{i}\|
\overset{\eqref{SublevelBound}}{\leq}
\|Y_{i}\|+\|Y_{i}-X_0\|
\le\
2\|Y_{i}\|+\|X_0\|
\leq
\frac{D_0}{2}.
\end{equation}
Thus, $X_i, Y_i \in B[D_0/2]$, and every tentative point $\tilde X$ computed in \eqref{def:ak-sfista1Tenative} of the $(i+1)$-st iteration also lies in $B[D_0/2]$.

This fact, the way $Y$ is computed in \eqref{eq:ynext-sfista1Tenative}, the last conclusion in Lemma~\ref{lem:stepTwoTermination}, and the check in \eqref{ineq check}
imply that if the Lipschitz estimate $L$ generated in step 2 of ADAP-FISTA becomes large enough, then the line-search in step 2 of the $(i+1)$-st iteration must terminate. This fact and \eqref{setIterates} imply that $a_i,\tilde X_{i},Y_{i+1}$, and $L_{i+1}$ are all generated and that $\tilde X_{i}$ lies in $B[D_0/2]$. We now show that the inequalities in \eqref{upper bound} hold for $L_{i+1}$. It is immediate from the way that $L$ is updated in step 2 and the way $L_{i+1}$ is set in \eqref{setIterates} that $L_0\leq L_{i+1}$. We now show by contradiction that the second relation in \eqref{upper bound} holds for $L_{i+1}$.
Suppose then that the line search in step 2 of ADAP-FISTA's $(i+1)$-st iteration terminates with
    $L_{i+1}$ that satisfies 
    \begin{equation}\label{Contradiction}
    L_{i+1}> \max\left\{L_0,12 G_0,\ \frac{2 ({\rm lip})_0}{1-\chi} \right\}.
    \end{equation}
It then follows that \eqref{ineq check} was not satisfied with 
$L=L_{i+1}/2$ and the $\tilde X$ and $Y$ in \eqref{def:ak-sfista1Tenative} and \eqref{eq:ynext-sfista1Tenative} corresponding to this $L$. It then follows from \eqref{Contradiction}, the fact that $\tilde X\in B[D_0/2]$, and the definition of $G_0$ in \eqref{def:D0_G0_Lip0} that
\[L=\frac{L_{i+1}}{2}> \max\left\{6 G_0,\ \frac{({\rm lip})_0}{1-\chi} \right\}\overset{\eqref{def:D0_G0_Lip0}}{\geq}\max\left\{ 6\|\nabla \psi(\tilde X)\|,\ \frac{({\rm lip})_0}{1-\chi} \right\}.\]
Hence, it follows from this inclusion, the update for $Y$ in \eqref{eq:ynext-sfista1Tenative}, and Lemma~\ref{lem:stepTwoTermination} with $\theta=(1-\chi)$, $\hat X=\tilde X$, $\hat L=L$, and $Y_{\hat L}^{+}=Y$
that
$Y\in B[D_0]$ and
\begin{equation}\label{ContraIneq}
\ell_{\psi}(Y;\tilde X)+\frac{(1-\chi) L}{2}\|Y-\tilde X\|^2\geq \psi(Y)
\end{equation}
which contradicts that \eqref{ineq check} was not satisfied with $L=L_{i+1}/2$, $\tilde X$, and $Y$.
\end{proof}

\begin{lemma}
    For any iteration index $i\geq 1$ where ADAP-FISTA does not fail in its step 4, it holds that
    \begin{align}
      &Y_i\in B[D_0]\label{Yinclusion}\\
      &V_{i} \in \nabla \psi(Y_{i}) + \mathcal{N}_{\Omega}(Y_{i}), \quad  \|V_{i}\|\leq \zeta \|Y_{i}-\tilde X_{i-1}\|\label{ubound-1}.
    \end{align}
\end{lemma}
\begin{proof}
Let $i\geq 1$ be an iteration index where ADAP-FISTA does not fail in its step 4. It follows that the second relation in \eqref{ineq: ineq 1} must hold. This observation and relation \eqref{eq:ynext-sfista1} imply that $Y_i\in \mathcal S(X_0;\psi)$. It then follows from this fact and \eqref{def:D0_G0_Lip0} that $Y_i\in B[D_0/6]\subseteq B[D_0]$, which shows \eqref{Yinclusion}.

It follows from the fact that ADAP-FISTA did not fail in its step 4 that it reaches step 5 and $V_i$ in \eqref{def:uk} is computed. It follows from \eqref{eq:ynext-sfista1} and its associated optimality condition that
$0 \in \nabla \psi(\tilde X_{i-1}) + \mathcal{N}_{\Omega}(Y_{i})-L_{i}(\tilde X_{i-1}-Y_{i})$, which in view of the definition of $V_{i}$ in \eqref{def:uk} implies the inclusion in \eqref{ubound-1}.

Lemma~\ref{LineSearchLemma} and \eqref{Yinclusion} imply that $\tilde X_{i-1}$ and $Y_i$ both lie in $B[D_0]$. It follows from these observations, triangle
inequality, 
the fact that $\nabla \psi$ is $({\rm lip})_0$-continuous, the bound \eqref{upper bound} on $L_{i}$, and the
definition of $\zeta$ in \eqref{omega zeta} that
\begin{align*}
    \|V_i\|&\overset{\eqref{def:uk}}{\leq} \|\nabla \psi(Y_i)-\nabla \psi(\tilde X_{i-1})\|+L_i\|\tilde X_{i-1}-Y_i\|\\
    &\leq ({\rm lip})_0\|Y_i-\tilde X_{i-1}\|+L_i\|Y_i-\tilde X_{i-1}\|\overset{\eqref{upper bound}}{\leq} \zeta\|Y_i-\tilde X_{i-1}\|
\end{align*}
which immediately implies that the inequality in \eqref{ubound-1} holds.
\end{proof}

The following proposition states that ADAP-FISTA terminates with success or failure in a logarithmic number of iterations.
\begin{proposition}\label{exact ACG complexity}
ADAP-FISTA
always stops (with either success or failure)
and does so by performing
at most 
\begin{equation}\label{eq:eq1}
\left\lceil\left(1+Q\right)\log_2^+\left(\frac{\zeta^2}{\chi \sigma^{2}}\right)+1\right\rceil+\left\lceil \log_2^+\left(\max\left\{\frac{6G_0}{L_0},\ \frac{({\rm lip})_0}{(1-\chi)L_0} \right\}\right)\right\rceil
\end{equation}
iterations/resolvent evaluations,
where $D_0$, $G_0$, and $({\rm lip})_0$ are as in \eqref{def:D0_G0_Lip0} and $\zeta$ and $Q$ are as in
\eqref{omega zeta}
\end{proposition}
\begin{proof}
Let $l$ denote the first quantity in \eqref{eq:eq1}. Using 
this definition and the
inequality $\log (1+ \alpha) \ge \alpha/(1+\alpha)$ for
any $\alpha\geq 0$, it is easy to verify that
\begin{equation}\label{Q bound}
\left(1+Q^{-1}\right)^{2(l-1)} \ge 
\frac{\zeta^2}{\chi \sigma^2}.
\end{equation}
We claim that ADAP-FISTA terminates with success or failure in at most $l$ iterations.
% Indeed, assume for contradiction that
Indeed, it suffices to show that
if ADAP-FISTA has not stopped with success or failure before the
$l$-th iteration and does not fail during the $l$-th iteration, then
it must stop successfully
at the $l$-th iteration.
% \jacob{If ADAP-FISTA stops successfully before iteration $l$, are $Y_l$ and $V_l$ still available for the argument below?}
% \arnesh{We are assuming that ADAP-FISTA has not stopped with success or failure up to the l-th iteration. We show then it must stop successfully in the l-th iteration if it does not stop with failure in the l-th iteration. I have rewritten the above sentence to make this more clear to address your point.}
So, assume that
ADAP-FISTA has not stopped with failure
up to (and including) the
$l$-th iteration.
In view of step~4 of ADAP-FISTA,
it follows that
\eqref{ineq: ineq 1} holds with $i=l$.

This observation, inequality in \eqref{ubound-1} with $i=l$, the first inequality in \eqref{ineq: ineq 1} with $i=l$, inequality \eqref{AiLiBound} with $i=l$,
and \eqref{Q bound}
then imply that
\begin{equation}\label{bound diff-1}
\|Y_{l}-\bar Y_{0}\|^{2} \overset{\eqref{ineq: ineq 1}}{\geq} \chi A_{l}L_{l} \|Y_{l}-\tilde X_{l-1}\|^2
\overset{\eqref{ubound-1}}{\ge}
\frac{\chi}{\zeta^2}A_lL_l \|V_{l}\|^2\overset{\eqref{AiLiBound}}{\geq}\frac{\chi}{\zeta^2}
\left(1+Q^{-1} \right)^{2(l-1)}
\|V_l\|^2 \overset{\eqref{Q bound}}{\geq} \frac{1}{\sigma^2}\|V_l\|^2,
\end{equation}
and hence that
\eqref{u sigma criteria} is satisfied.
In view of Step 5 of
ADAP-FISTA, the method
must
successfully stop at the
end of the $l$-th iteration.
We have thus shown that
the above claim holds.
Moreover, in view of \eqref{upper bound},
it follows that
the second term in 
\eqref{eq:eq1} is a bound
on the total number of times $L_i$ is multiplied by $2$ and step 2 is repeated. Since exactly one resolvent evaluation occurs every time step 2 is executed, the desired conclusion follows.
% \jacob{Could we justify a bit more carefully why this second term bounds the total number of doublings across all warm-started line searches?}
% \arnesh{I think it is fine to not add more details as it is a standard argument. I have read other papers that make similar observation about extra grad evaluations from the line-search without going into too many details.}
\end{proof}

\begin{proof}[Proof of Proposition~\ref{prop:nest_complex1}(a)-(b)]
(a) It follows from the definition of $\zeta$ in \eqref{omega zeta} that $\mathcal O(\log_{2}^{+}(\zeta^2))=\mathcal O(\log_{2}^{+}(\zeta))=\mathcal O\left(\log_{2}^{+}\left(\max\left\{L_0,({\rm lip})_0, G_0\right\}\right)\right)$.
The result immediately follows from this observation, Proposition~\ref{exact ACG complexity},
and the definitions of 
$Q$ and $N_0$ in \eqref{omega zeta} and \eqref{ND0 Def}, respectively.

(b) It follows immediately from the termination criterion \eqref{u sigma criteria} in step 5 of ADAP-FISTA and the inclusion in \eqref{ubound-1} that \eqref{u inclusion 1} must hold if ADAP-FISTA stops successfully.

To see \eqref{secondOutput}, first observe that since ADAP-FISTA stops successfully, the second relation in \eqref{ineq: ineq 1} held in ADAP-FISTA's last iteration. Hence, this observation together with relation \eqref{eq:ynext-sfista1} imply that $Y\in \Omega$ and $\psi(\bar Y)\leq \psi(\bar Y_0)$. The bound \eqref{Lbound} follows immediately from the fact that relation \eqref{upper bound} holds for all iteration indices $i\geq 1$ generated by ADAP-FISTA.
\end{proof}

\subsubsection{Proof of Proposition~\ref{prop:nest_complex1}(c)}
This subsection proves Proposition~\ref{prop:nest_complex1}(c). Throughout, we assume that $\psi$ is $\mu$-strongly convex on $B[D_0]$. The following lemma proves key facts of the $i$-th iteration of ADAP-FISTA.

% \jacob{Since Lemma~\ref{lm:hysub-4-fista} seems to drive part~(c), would it help to include a proof here rather than omitting it?}
% \arnesh{The proof is identical to previous papers so I do not believe it is worth to rehash it.}

\begin{lemma} \label{lm:hysub-4-fista}
% ({\bf Result translated from Lemma~A.7 of \cite{sujanani2023adaptive}})
For every iteration $i \ge 1$ of ADAP-FISTA, the following statements hold:
\begin{itemize}
\item[(a)] for every $X \in  \Omega\cap 
B[{D_0}]$, we have
\begin{equation*}
\left(A_{i-1} [ \psi(Y_{i-1}) - \psi(X) ] + \frac{\tau_{i-1}}{2} \|X-X_{i-1}\|^2\right) - \left(A_i [ \psi(Y_i) - \psi(X) ] + \frac{\tau_i}{2} \|X-X_i\|^2\right) \ge   \frac{\chi A_{i}L_{i}}{2} \|Y_{i} - \tilde X_{i-1} \|^2;
\end{equation*}
\item[(b)] for every $X \in  \Omega\cap 
B[{D_0}]$, we have
\begin{align}\label{key convergence inequality 2}
 A_{i} [\psi(Y_{i}) - \psi(X) ] + \frac{\tau_{i}}2 \|X-X_{i}\|^2 \le
 \frac12 \|X-\bar Y_0\|^2 - \frac{\chi}{2} \sum_{j=1}^{i} A_{j}L_{j} \| Y_{j} - \tilde X_{j-1} \|^2.
\end{align}
\item[(c)] the method does not fail at iteration $i$.
\end{itemize}
\end{lemma}

\begin{proof}
(a) 
For every
iteration $i\geq 1$ of ADAP-FISTA and $X\in \Omega$,
define
\begin{align}
		\gamma_i(X)&:=\ell_{\psi}(Y_{i},\tilde X_{i-1})+\inner{(L_i-\mu)(\tilde X_{i-1}-Y_i)}{X - Y_{i}}+\frac{\mu}{2}\|Y_{i}-\tilde X_{i-1}\|^2
		+ \frac{\mu}2 \|X-Y_{i}\|^2. \label{def:gamma-sfista}
	\end{align}
It is now easy to see from similar arguments as in Lemma A.7 of \cite{sujanani2023adaptive} that
for every iteration $i\geq 1$ and $X \in \Omega$, we have
\begin{align}\label{ineq:recur}
	A_{i-1}\gamma_i(Y_{i-1}) &+ a_{i-1}\gamma_i(X) + \frac{\tau_{i-1}}2 \|X_{i-1} - X\|^2 - \frac{\tau_{i}}2 \|X_{i} - X\|^2 \nonumber \\
	&\ge A_{i}\psi(Y_{i})+\frac{\chi A_{i} L_{i}}{2}\|Y_{i}-\tilde X_{i-1}\|^2.  
	\end{align}
To get an inequality that telescopes, we now show that
$\gamma_{i}(X)\leq \psi(X)$ for all $X\in \Omega\cap B[D_0]$.
To show the claim for every iteration $i\geq 1$,
define
    \begin{equation}\tilde \gamma_i(X):= \ell_{\psi}(X;\tilde X_{i-1}) +\frac{\mu}{2}\|X-\tilde X_{i-1}\|^2 \label{def:tgamma-sfista}.
    \end{equation}
Now observe that Lemma~\ref{LineSearchLemma} implies that $\tilde X_{i-1}\in B[D_0/2]\subseteq B[D_0]$.
It follows from this observation and the assumption that
$\psi$ is $\mu$-strongly convex on $B[D_0]$  that $\tilde \gamma_i(X) \leq \psi(X)$ for all $X\in B[D_0]$. The rest of the argument showing that $\gamma_{i}(X)\leq \psi(X)$ for all $X\in \Omega\cap B(D_0)$ is identical to the one in Lemma A.8 in \cite{sujanani2023adaptive}.

Lemma~\ref{LineSearchLemma} imply that $Y_{i-1}\in B[D_0/2]\cap \Omega \subseteq B[D_0]\cap \Omega$. This observation, the fact that $\gamma_{i}(X)\leq \psi(X)$ for all $X\in B[D_0]\cap \Omega$, and \eqref{ineq:recur} imply that 
\begin{align}\label{ineq:recur2}
	A_{i-1}\psi(Y_{i-1}) &+ a_{i-1}\psi(X) + \frac{\tau_{i-1}}2 \|X_{i-1} - X\|^2 - \frac{\tau_{i}}2 \|X_{i} - X\|^2 \nonumber \\
	&\ge A_{i}\psi(Y_{i})+\frac{\chi A_{i} L_{i}}{2}\|Y_{i}-\tilde X_{i-1}\|^2
	\end{align}
holds for all $X\in B[D_0]\cap \Omega$.
The relation in part (a) then immediately follows from subtracting $A_i\psi(X)$ from both sides of \eqref{ineq:recur2} and using the fact that \eqref{eq:taunext-sfista1} implies
$A_i=A_{i-1}+a_{i-1}$.

(b) It follows from summing the inequality of in part (a) from $1$ to $i$, and using the facts that $X_0=\bar Y_0$, $A_0=0$ and $\tau_0=1$ that
\begin{align*}
A_{i} [\psi(Y_{i}) - \psi(X) ] + \frac{\tau_{i}}2 \|X-X_{i}\|^2 \le
 \frac12 \|X-\bar Y_0\|^2 - \frac{\chi}{2} \sum_{j=1}^{i} A_{j}L_{j} \| Y_{j} - \tilde X_{j-1} \|^2.
\end{align*}  
which immediately implies that \eqref{key convergence inequality 2} holds.

(c)  To show that ADAP-FISTA does not fail at its iteration $i$, we show that all three inequalities checked in \eqref{ineq: ineq 1} and \eqref{ineq 2} in its step 4 must hold.
We first show that $\psi(Y_i)\leq \psi(\bar Y_0)$. First observe that it is easy to see that $\bar Y_0\in B[D_0/6]\cap \Omega\subseteq B[D_0]\cap \Omega$ in view of the definition of $D_0$ in \eqref{def:D0_G0_Lip0}. Thus, using \eqref{key convergence inequality 2} with $X=\bar Y_0$ and the facts that the sequences $\tau_{i}$, $A_{i}$, and $L_i$ are nonnegative,
we see that 
\[A_{i}[\psi(Y_{i}) - \psi(X_0) ]\leq \frac{1}{2}\|\bar Y_0-\bar Y_0\|^2=0\]
Hence, diving both sides by $A_{i}$,
we see that $\psi(Y_{i})\leq \psi(\bar Y_0)$.
To conclude the proof, we now show that
\begin{equation}\label{YiDist}
    \|Y_i-X_i\|\leq \|Y_i-X_0\|, \quad \|Y_{i}-X_{0}\|^{2} \geq \chi A_{i}L_{i}\|Y_{i}-\tilde X_{i-1}\|^2.
\end{equation}
It follows from above conclusion, \eqref{eq:ynext-sfista1}, and the definition of $D_0$ in \eqref{def:D0_G0_Lip0} that $Y_{i} \in B[D_0]\cap \Omega$.
    Hence, using \eqref{key convergence inequality 2} with $X=Y_{i}$ and the fact that $\tau_{i}\geq 1$, we see that $\|Y_{i}-X_{i}\|\leq \|Y_{i}-\bar Y_0\|$.
    Likewise also using \eqref{key convergence inequality 2} with $X=Y_{i}$ but keeping the last term on the right-hand side of \eqref{key convergence inequality 2},
    we see that
    $\|Y_{i}-\bar Y_{0}\|^{2} \geq \chi A_{i}L_{i}\|Y_{i}-\tilde X_{i-1}\|^2$. Hence, the three inequalities checked in \eqref{ineq: ineq 1} and \eqref{ineq 2} must hold at iteration $i$, which implies that ADAP-FISTA does not fail in its $i$-ith iteration.
\end{proof}

We are now ready to prove Proposition~\ref{prop:nest_complex1}(c).

\begin{proof}[Proof of Proposition~\ref{prop:nest_complex1}(c)]
It follows from Lemma~\ref{lm:hysub-4-fista}(c) that, for any iteration $i\geq 1$, ADAP-FISTA does not terminate with failure at step 4 when applied to a function $\psi$ that is $\mu$-strongly convex on $B[D_0]$ and satisfies assumptions \textbf{(B1)}–\textbf{(B3)}.
Hence, it follows from Proposition~\ref{prop:nest_complex1}(a)-(b)
that ADAP-FISTA must always terminate successfully and return a pair
$(\bar Y,\bar V) \in \Omega \times \mathbb E$ satisfying \eqref{u inclusion 1} and \eqref{secondOutput} and a Lipschitz estimate $\bar L$ satisfying \eqref{Lbound}.
% To show \eqref{KeyInequality}, 
% first observe that since ADAP-FISTA must terminate with success, its output $\bar Y$ satisfies the second inequality in \eqref{ineq: ineq 1}. This fact and and the definition $D_0$ in \eqref{def:D0_G0} imply that both
% $X_0$ and output $\bar Y$ both lie in $B[D_0]$. 
% It follows from this observation together with the fact $\psi$ is $\mu$-strongly convex on $B[D_0]$
% that
% \begin{equation}\label{ConvBound}
%     \psi(X_0)\geq \psi(\bar Y)+\inner{\nabla \psi(\bar Y)}{X_0-\bar Y}.
% \end{equation}
% It follows from the inclusion in \eqref{u inclusion 1} that the output $(\bar Y,\bar V)$ satisfies $\bar V-\nabla\psi(\bar Y)\in \mathcal{N}_{\Omega}(\bar Y)$. Hence using the definition of normal cone and the fact that $X_0\in \Omega$,
% we have that 
% \begin{equation}\label{NormalConeBound}
% \inner{\nabla \psi(\bar Y)}{X_0-\bar Y}\geq \inner{\bar V}{X_0-\bar Y}.
% \end{equation}
% Inequality \eqref{KeyInequality} then follows immediately from combining relations \eqref{ConvBound} and \eqref{NormalConeBound}.
\end{proof}

\subsection{ADAP-AIPP Method and Complexity}\label{app:aipp}

Let $\mathbb E$ be a finite-dimensional inner product space. This section considers
\begin{equation*}
g_{*}:=\min\left\{g(U):\ U\in \Omega\right\},
\end{equation*}
under the following assumptions:

\begin{itemize}
\item[\textbf{(C1)}]
$g$ is a real-valued function that is differentiable on $\mathbb E$ and
whose gradient is locally Lipschitz continuous in the sense of
\eqref{eq:locally_lipschitz_gradient};

\item[\textbf{(C2)}]
$g$ is bounded below on $\Omega$;

\item[\textbf{(C3)}] $\Omega$ is closed, possibly unbounded, convex set that admits an easily computable projection; 

\item[\textbf{(C4)}] 
a point $W_0\in \Omega$ is given such that the sublevel set
$\mathcal S(W_0;g)$ has finite diameter $D(W_0;g)$,
where $\mathcal S(\cdot;\cdot)$ and $D(\cdot;\cdot)$
are as in \eqref{SublevelDef} and \eqref{DiameterDef}, respectively.
\end{itemize}
\noindent{\bf Problem:}
Given $(g, \Omega, W_0)$ satisfying \textbf{(C1)-(C4)} and a tolerance $\bar\rho\in\R_{++}$, ADAP-AIPP seeks
$(W,R)\in \Omega\times\mathbb E$ such that
\begin{equation*}
R\in \nabla g(W)+\mathcal{N}_{\Omega}(W),\qquad \|R\|\le \bar\rho.
\end{equation*}
At its $j$-th iteration, calls ADAP-FISTA several times until it sucessfully finds
an approximate stationary point of
\begin{equation}\label{eq:aipp_prox_subproblem}
\min_{U\in\Omega}\ \left\{\psi_{\eta}(U):=\eta g(U)+\tfrac12\|U-W_{j-1}\|^2\right\}
\end{equation}
where $\eta>0$ is a prox stepsize that is halved after an unsuccessful ADAP-FISTA call.

\begin{minipage}[t]{\linewidth}
\rule[0.5ex]{\linewidth}{1pt}

\noindent \textbf{ADAP-AIPP Method}

\noindent \rule[0.5ex]{\linewidth}{1pt}
\end{minipage}

\noindent \textbf{Universal Parameters (used by ADAP-FISTA)}: $\sigma\in(0,1/2)$, $\chi\in(0,1)$.

\smallskip
\noindent \textbf{Input}: an initial point $W_0\in \Omega$, a function $g$ satisfying \textbf{(C1)}-\textbf{(C4)}, a set $\Omega$, prox stepsize $\eta_0>0$, an initial Lipschitz estimate $M_0$, and tolerance $\bar\rho>0$.

\smallskip
\noindent \textbf{Step 0.} Set $j\gets1$ and $\eta\gets\eta_0$.

\smallskip
\noindent \textbf{Step 1.} 
Call ADAP-FISTA with input
\begin{equation}\label{eq:aipp_adap_fista_call}
\bar Y_0=W_{j-1}\,,\quad \mu=\frac{1}{2}\,,\quad L_0=M_{j-1}\,,\quad \psi(\cdot) =\eta g(\cdot)+\frac{1}{2}\|\cdot-W_{j-1}\|^2, \quad \Omega=\Omega.
\end{equation}

\noindent\smallskip
\textbf{Step 2.}
\textbf{If} the above ADAP-FISTA call terminates successfully, then let $(W_j, V_j,M_j)$ be its output and set
\begin{equation}\label{eq:aipp_Rj_def}
\begin{aligned}
\eta_j\gets \eta,\quad R_j\ \gets\ \frac{V_j+W_{j-1}-W_j}{\eta_j};
\end{aligned}
\end{equation}
{\bf else}, set $\eta\gets \eta/2$ and return to Step~{\bf 1}.

\smallskip
\noindent \textbf{Step 3.} \textbf{If} $\|R_j\|\le \bar\rho$, then output $(\bar{W},\bar{R})\gets(W_j,R_j)$ and {\bf stop}; {\bf else}, set $j\gets j+1$ and return to Step~{\bf 1}.

\noindent\rule[0.5ex]{\linewidth}{1pt}

We now define several quantities which arise in ADAP-AIPP's technical analysis and complexity result.

Let
\begin{equation}\label{DiamDefAIPP}
    \mathcal D_{j}:=6D(W_{j-1};g)+1, \quad \mathcal D_0:=\mathcal D_1
\end{equation}
where $j\geq 1$ is an iteration index and $D(\cdot; \cdot)$ is as in \eqref{DiameterDef}.
For every $j\geq 0$, we define
\begin{equation}\label{SUPLIP}
    \nabla_j:=\nabla(\mathcal D_j/2;g),  \quad ({\rm Lip})_j:=[{\rm LIP}](\mathcal D_j;g)
\end{equation}
where $\nabla(\cdot; \cdot)$, and 
$[{\rm LIP}](\cdot;\cdot)$ are as in \eqref{GradDef} and \eqref{SmoothnessDef}, respectively.
We also let
\begin{equation}\label{Csigma}
    C_{\sigma}=\frac{2(1-\sigma)^2}{(1-2\sigma)^2}
\end{equation}
where $\sigma\in(0,1/2)$ is an universal input to ADAP-AIPP.

The next lemma states useful facts about the sublevel sets of the prox function in \eqref{eq:aipp_prox_subproblem}.

\begin{lemma}\label{KeyBounds}
    Let $j\geq 1$ be an iteration index generated by ADAP-AIPP and consider $\psi_{\eta}(\cdot)$ as in \eqref{eq:aipp_prox_subproblem}, i.e.,
    $\psi_{\eta}(\cdot):=\eta g(\cdot)+0.5\|\cdot-W_{j-1}\|^2$ where $\eta>0$ is a prox stepsize. The following facts then hold:
   \begin{itemize}
       \item[(a)] the diameter
       $D(W_{j-1};\psi_{\eta})\leq D(W_{j-1};g)$;
       \item[(b)]  $\nabla \psi_{\eta}(\cdot)$ is $1+\eta \Lip{j}$ Lipschitz continuous on $B[\mathcal D_{j}]$;
       \item[(c)] it holds that
       \begin{equation}\label{SupBound}
       \sup\left\{\|\nabla \psi_{\eta}(U)\|:\ U\in B[\mathcal D_j/2]\right\}\leq \eta \nabla_j+\mathcal D_{j}.
       \end{equation}
   \end{itemize}
\end{lemma}

\begin{proof}
(a) Let $U'\in \mathcal S(W_{j-1};\psi_{\eta})$. Hence, it follows from the definition of $\mathcal S(\cdot;\cdot)$ in \eqref{SublevelDef}
that
$\psi_{\eta}(U')\leq \psi_{\eta}(W_{j-1})$.
It then follows from the definition of $\psi_{\eta}(\cdot)$ that
$\eta g(U')\leq\eta g(U')+0.5\|U'-W_{j-1}\|^2 \leq \eta g(W_{j-1})$, which implies that $U'\in \mathcal S(W_{j-1};g)$.
Thus, it follows that $\mathcal S(W_{j-1},\psi_{\eta})\subseteq \mathcal S(W_{j-1};g)$
which implies that $D(W_{j-1};\psi_{\eta})\leq D(W_{j-1};g)$
in view of the definition of 
$D(W_{j-1};g)$ in \eqref{DiameterDef}.

(b)  It is easy to see that 
\begin{equation}\label{FormGradientPsi}
\nabla \psi_{\eta}(U)=\eta \nabla g(U)+(U-W_{j-1}).
\end{equation}
The result then follows from \eqref{FormGradientPsi} and the fact that $\nabla g$ is $\Lip{j}$-Lipschitz continuous on $B[\mathcal D_j]$.

(c) It is immediate from \eqref{FormGradientPsi} and the definition of $\nabla_j$ in \eqref{SUPLIP}
that
\begin{align}
\sup\left\{\|\nabla \psi_{\eta}(U)\|:\ U\in B[\mathcal D_j/2]\right\}&\leq \eta \nabla_j+\sup_{U\in B[\mathcal D_{j}/2]}\|U\|+\|W_{j-1}\|\\
&\leq \eta\nabla_j+\frac{\mathcal D_j}{2}+\frac{\mathcal D_j}{6}\leq \eta\nabla_j+\mathcal D_{j}
\end{align}
where the second inequality follows from the definition of $\mathcal D_j$ 
in \eqref{DiamDefAIPP} and the fact that $W_{j-1}\in \mathcal S(W_{j-1};g)$.
\end{proof}

\begin{lemma}\label{Step1Lemma}
    Let $j\geq 1$ be an iteration index of ADAP-AIPP, and define \begin{equation}\label{aippConstants}
\begin{aligned}
\Theta_j
&\defi
\max\left\{
M_{j-1},\ \eta_0\nabla_j+\cD_j,\
1+\eta_0 \Lip{j}
\right\}.
\end{aligned}
\end{equation}
Then the following statements hold:
    \begin{itemize}
        \item[(a)] the loop in steps 1 and 2 of the $j$-th iteration is repeated at most
        a logarithmic number of
        times before terminating with a triple $(W_j,V_j,M_j)\in\Omega\times \mathbb E\times \mathbb R_{++}$ and $\eta_j$ satisfying
\begin{align}
& \eta_0\geq\eta_j\geq \min\left\{\eta_0,1/(4\Lip{j}) \right\} , \quad
R_j\in \nabla g(W_j)+\mathcal{N}_{\Omega}(W_j),\qquad
\label{Output1}\\
&\|V_j\|\le\sigma\|W_j-W_{j-1}\|, \quad \eta_{j}\left(g(W_{j-1})-g(W_j)\right)\geq \frac{1}{2}\|W_j-W_{j-1}\|^2,\label{Output2}
\end{align}
and
\begin{equation}\label{AIPPMBound}
M_{j-1}\leq M_j
\le
\max\left\{
M_{j-1},12(\eta_j\nabla_j+\cD_j),
\frac{2(1+\eta_j \Lip{j})}{1-\chi}
\right\};
\end{equation}

\item[(b)]  every ADAP-FISTA call made in step 1 within this loop terminates with either success or failure, in at most
\begin{equation}\label{eq:aipp_complexity}
\mathcal O_1\left(
\sqrt{\Theta_j}\log_2^+\!\left(\Theta_j\right)
+
\log_2^+\!\left(
\frac{\Theta_j}{M_{j-1}}\right)
\right)
\end{equation}
ACG iterations/resolvent evaluations;
\item[(c)] the inequality holds:
\begin{align}
g(W_{j-1})-g(W_j)\ \ge\ \frac{\eta_j}{C_\sigma}\|R_j\|^2.\label{FuncBound}
\end{align}
where $C_{\sigma}$ is as in \eqref{Csigma}.
    \end{itemize}
\end{lemma}

\begin{proof}
   (a) Let $j\geq 1$ be an iteration index of ADAP-AIPP and consider its $j$-th iteration. It follows from step 1 and Lemma~\ref{KeyBounds}(b) that ADAP-AIPP calls ADAP-FISTA with a function $\psi(\cdot) =\eta g(\cdot)+\frac{1}{2}\|\cdot-W_{j-1}\|^2$ whose gradient is $1+\eta \Lip{j}$ Lipschitz continuous on $B[\mathcal D_j]$. Thus, this function is $(1-\eta \Lip{j})$-strongly convex on $B[\mathcal D_j]$. Thus, it follows from this observation, the fact that ADAP-FISTA is called with $\mu=1/2$, and Proposition~\ref{prop:nest_complex1}(c) that if ADAP-FISTA is called with $\psi(\cdot) =\eta g(\cdot)+\frac{1}{2}\|\cdot-W_{j-1}\|^2$ where $\eta\leq 1/(2\Lip{j})$, then the call to ADAP-FISTA must terminate successfully. 
   The first conclusion of part (a) and the first relation in \eqref{Output1} then immediately follow from this observation, the fact that $\eta$ is halved after every unsuccessful ADAP-FISTA call, and the way $\eta_j$ is set in \eqref{eq:aipp_Rj_def}.

 It then follows from Proposition~\ref{prop:nest_complex1}(b), the fact that ADAP-AIPP calls ADAP-FISTA with input as in \eqref{eq:aipp_adap_fista_call}, and the way $\eta_j$ is set in \eqref{eq:aipp_Rj_def} that ADAP-FISTA terminates with a triple $(W_j,V_j,M_j)$ and $\eta_j$ satisfying
\begin{align}
&  V_j\in \eta_j\nabla g(W_j)+(W_j-W_{j-1})+\mathcal{N}_{\Omega}(W_j), \quad \|V_j\|\leq \sigma\|W_j-W_{j-1}\|\label{OutputFISTA1}\\
&\eta_{j}\left(g(W_{j-1})-g(W_j)\right)\geq \frac{1}{2}\|W_j-W_{j-1}\|^2\label{OutputFISTA2}
\end{align}
which immediately implies both relations in \eqref{Output2}. The second relation in \eqref{Output1} is immediate from the inclusion in \eqref{OutputFISTA1} and the definition of $R_j$ in \eqref{eq:aipp_Rj_def}.
Relation \eqref{AIPPMBound} follows from the fact that ADAP-FISTA is called with input as in \eqref{eq:aipp_adap_fista_call}, relation \eqref{Lbound}, Lemma~\ref{KeyBounds}(a)-(c), and the way that $\eta_{j}$ is set in \eqref{eq:aipp_Rj_def}.

(b) 
Let $j\geq 1$ and consider $\psi_{\eta}:=\eta g(\cdot)+0.5\|\cdot-W_{j-1}\|^2$ where $\eta>0$.
It follows from Lemma~\ref{KeyBounds}(a)
and the definition of $\mathcal D_{j}$ in \eqref{DiamDefAIPP}
that 
$D(W_{j-1};\psi_{\eta})\leq D(W_{j-1};g)\leq \mathcal D_{j}$.
This relation together with \eqref{SupBound} then imply that
\[\sup\left\{\|\nabla \psi_{\eta}(U)\|:\ U\in B[D(W_{j-1};\psi_{\eta})/2]\right\}\leq \sup\left\{\|\nabla \psi_{\eta}(U)\|:\ U\in B[\mathcal D_j/2]\right\}\leq \eta \nabla_j+\mathcal D_{j}.\]
It then follows from these conclusions, Proposition~\ref{prop:nest_complex1}(a), the fact that ADAP-FISTA is called with input as in \eqref{eq:aipp_adap_fista_call}, the fact that the way $\eta$ is updated in step 2 implies that $\eta\leq \eta_0$ for all $\eta$ generated by ADAP-AIPP, and the definition of $\Theta_{j}$ in 
\eqref{aippConstants}
that every ADAP-FISTA call within the loop in steps 1 and 2 performs at most
\eqref{eq:aipp_complexity} 
ACG iterations/resolvent evaluations. 

(c) Let $\tilde \sigma=\sigma/(1-\sigma)$ where $\sigma\in(0,1/2)$ is an universal input to ADAP-AIPP.
It follows from the first relation in \eqref{Output2} and the definition of $R_j$ in \eqref{eq:aipp_Rj_def} that
        $\|V_j\|\leq \sigma\|W_j-W_{j-1}\|\leq \sigma\|V_j\|+\sigma\|\eta_jR_j\|$. This relation together with the definition of $\tilde \sigma$ then imply that
\begin{equation}\label{V_jBound}
        \left\|\frac{V_j}{\eta_{j}}\right\|\leq \tilde \sigma \|R_j\|.
    \end{equation}
It then follows from the definition of $R_j$ in \eqref{eq:aipp_Rj_def}, relation \eqref{V_jBound},
and multiplying both sides of the second inequality in
\eqref{Output2} by $2/\eta_j^2$ that
\begin{align}
\frac{2}{\eta_j}[g(W_{j-1})-g(W_j)]
\overset{\eqref{Output2}}{\ge}
\left\|\frac{W_{j-1}-W_{j}}{\eta_j}\right\|^2&\overset{\eqref{eq:aipp_Rj_def}}{=} \left\|R_j-\frac{V_j}{\eta_j} \right\|^2 \nonumber\\
&\overset{\eqref{V_jBound}}{\ge}
\min_{\|z\|\le \tilde \sigma\|R_j\|}\|R_j-z\|^2 \nonumber\\
&=
\min_{0\le t\le \tilde\sigma\|R_j\|}
(\|R_j\|-t)^2 \nonumber\\
&=
(1-\tilde\sigma)^2\|R_j\|^2 \label{FinalRBound} .
\end{align}
Relation \eqref{FuncBound} then follows from multiplying both sides of \eqref{FinalRBound} by $\eta_j/2$, the fact that $\tilde \sigma=\sigma/(1-\sigma)$, and the definition of $C_{\sigma}$ in \eqref{Csigma}.
\end{proof}

\begin{lemma}
The sequence
$\left\{\mathcal D_{j}\right\}$
is nonincreasing, i.e., for any iteration index $j\geq 1$, it holds that
\begin{align}
    \mathcal D_{j}\leq \mathcal D_{j-1} \leq \mathcal D_{0}\label{DiamDec}
\end{align}
where recall $\mathcal D_{j}$ is as in \eqref{DiamDefAIPP}.
\end{lemma}

\begin{proof}
    Let $j\geq 1$ be an iteration index of ADAP-AIPP and let $U'\in \mathcal S(W_{j};g)$ where $\mathcal S(\cdot;\cdot)$ is as in \eqref{SublevelDef}. It follows from this inclusion, \eqref{FuncBound}, and the facts that $\eta_j$, $C_{\sigma}$, and $\|R_j\|^2$ are all nonnegative quantities that 
    \[g(U')\leq g(W_{j})\overset{\eqref{FuncBound}}{\leq} g(W_{j-1})-\frac{\eta_{j}}{C_{\sigma}}\|R_j\|^2\overset{\eqref{Csigma}}{\leq} g(W_{j-1}).\]
Hence, it follows from the above relation that $U'\in \mathcal S(W_{j-1};g)$ and thus $\mathcal S(W_{j};g) \subseteq \mathcal S(W_{j-1};g)$.
It then follows from the definition of $D(\cdot;\cdot)$ in \eqref{DiameterDef} that 
$D(W_{j};g)\leq D(W_{j-1};g)$. The first inequality in \eqref{DiamDec} then follows from this relation and the definition of $\mathcal D_{j}$ in \eqref{DiamDefAIPP}. The second inequality in \eqref{DiamDec} then follows from the first inequality in \eqref{DiamDec}, a simple induction argument, and the fact that $\mathcal D_1=\mathcal D_0$.
\end{proof}

\begin{lemma}
For every iteration index $j\geq 1$, it holds that
\begin{equation}\label{MUniformBound}
   M_0\leq M_j
\le
\max\left\{
M_{0},12(\eta_0\nabla_0+\cD_0),
\frac{2(1+\eta_0 \Lip{0})}{1-\chi}
\right\}.
\end{equation}
\end{lemma}
\begin{proof}
The first inequality in \eqref{MUniformBound} is immediate from the first inequality in \eqref{AIPPMBound} together with a simple induction argument. To show the second inequality in \eqref{MUniformBound}, we also proceed by induction. We do not include the argument for the $j=1$ case since it follows from a similar argument as below.
     
Suppose then that the second inequality in \eqref{MUniformBound}
holds for an iteration index $j\geq 1$ generated by ADAP-AIPP. We show that it also must hold for the $(j+1)$-st iteration index. It follows from 
the second inequality in
\eqref{AIPPMBound} with $j=j+1$, the bounds on $\mathcal D_{j+1}$ and $\eta_{j+1}$ in \eqref{DiamDec} and \eqref{Output1}, respectively, and the definitions of $\nabla_j$ and $\Lip{j}$ in \eqref{SUPLIP} that
\begin{align}
    M_{j+1}
&\overset{\eqref{AIPPMBound}}{\le}\max\left\{
M_{j},12(\eta_{j+1}\nabla_{j+1}+\cD_{j+1}),
\frac{2(1+\eta_{j+1} \Lip{j+1})}{1-\chi}
\right\}\nonumber\\
&\overset{\eqref{Output1}}{\leq} \max\left\{
M_{j},12(\eta_{0}\nabla_{j+1}+\cD_{j+1}),
\frac{2(1+\eta_{0} \Lip{j+1})}{1-\chi}
\right\}\nonumber\\
&\overset{\eqref{DiamDec}}{\leq} \max\left\{
M_{j},12(\eta_{0}\nabla_0+\cD_{0}),
\frac{2(1+\eta_{0} \Lip{0})}{1-\chi}
\right\}\label{BoundMJD}.
\end{align}
The fact that \eqref{MUniformBound} holds for all iteration indices $j\geq 1$ then follows immediately from \eqref{BoundMJD}
and the induction hypothesis.
\end{proof}

The following lemma uniformly bounds the quantity $\Theta_j$
in \eqref{aippConstants}. This is useful for uniformly bounding the number of gradient evaluations made during each ADAP-FISTA call in step 1.
\begin{lemma}\label{BarThetaLemma}
    Let
\begin{equation}\label{barThetaDef}
    \bar \Theta:=
\max\left\{
M_{0},\eta_0\nabla_0+\mathcal D_0,
1+\eta_0 \Lip{0}
\right\}.
\end{equation}
For every iteration index $j\geq 1$, it holds that
$\Theta_j=\mathcal O(\bar \Theta)$.
\end{lemma}
\begin{proof}
     Let $j\geq 1$ be an iteration index generated by ADAP-FISTA. It follows from \eqref{DiamDec} and the definitions of $\Lip{j}$ and $\nabla_{j}$ in \eqref{SUPLIP} that $\mathcal D_{j}\leq \mathcal D_0$, $\Lip{j}\leq \Lip{0}$, and $\nabla_j\leq \nabla_0$. It is then immediate from this bound, the bound on $M_{j}$ in \eqref{MUniformBound}, the fact that $\chi\in(0,1)$, and \eqref{barThetaDef} that
\begin{align}
    \Theta_{j}&\overset{\eqref{aippConstants}}{=}\max\left\{
M_{j-1},\ \eta_0\nabla_j+\cD_j,\
1+\eta_0 \Lip{j}
\right\} \nonumber\\
&\overset{\eqref{DiamDec}}{\leq} \max\left\{
M_{j-1},\ \eta_0\nabla_0+\cD_0,\
1+\eta_0 \Lip{0}
\right\}\\
&\overset{\eqref{MUniformBound}}{\leq} \max\left\{
M_{0},12\left(\eta_0\nabla_0+\mathcal D_0\right),
\frac{2\left(1+\eta_0 \Lip{0}\right)}{1-\chi}
\right\}\overset{\eqref{barThetaDef}}{=}\mathcal O(\bar \Theta)
\label{ThetaJBound}
\end{align}
which immediately implies the result of the lemma.
\end{proof}

The following proposition now states the main complexity results of ADAP-AIPP.

\begin{proposition}\label{prop:aipp_generic_complexity}
The following statements about ADAP-AIPP hold:
\begin{itemize}
\item[(a)] it terminates with a pair $(\bar W,\bar R)\in \Omega\times\mathbb E$ satisfying
\begin{align}
    &\bar R\in \nabla g(\bar W)+\mathcal{N}_{\Omega}(\bar W),\qquad \|\bar R\|\le \bar\rho\label{eq:aipp_stationary_output}\\
    & g(\bar W)\leq g(W_0);\label{MonoADAPAIPP}
\end{align}

\item[(b)] if $\ell$ denotes its final iteration index, then 
its total number of ADAP-AIPP iterations/ADAP-FISTA calls is bounded by
\begin{equation}\label{FISTACalls}
    1+\frac{C_{\sigma}\max\left\{\eta_0^{-1},\,4\Lip{0}\right\}}{\bar\rho^2}\left[g(W_0)-g(W_{\ell})\right]+ \left\lceil\log_2^+\left(4\eta_0\Lip{0}\right)\right\rceil
\end{equation}
where $\Lip{0}$ and $C_{\sigma}$ are as in \eqref{SUPLIP} and \eqref{Csigma}, respectively, and $\eta_0$, $\bar \rho$, and $W_0$ are inputs to ADAP-AIPP.
\end{itemize}
\end{proposition}

\begin{proof}
    (a) We first show that ADAP-AIPP must terminate
    and that its final iteration index $\ell\geq 1$
    satisfies
    \begin{equation}\label{FinalIterationIndex}
        \ell\leq  1+\frac{C_{\sigma}\max\left\{\eta_0^{-1},\,4\Lip{0}\right\}}{\bar\rho^2}\left[g(W_0)-g(W_{\ell})\right].
    \end{equation}
Suppose by contradiction that ADAP-AIPP generates some index $\bar j$  such that
\begin{equation}\label{iteration count}
\bar j>1+\frac{C_{\sigma}\max\left\{\eta_0^{-1},\,4\Lip{0}\right\}}{\bar\rho^2}\left[g(W_0)-g(W_{\bar j})\right].
\end{equation}
 As ADAP-AIPP generates $\bar j$ as an iteration index, it does not terminate at the $(\bar j-1)$-th iteration. In view of step 3 of ADAP-AIPP, this implies that
$\|R_j\| > \bar \rho$ for every $j=1,\ldots,\bar j-1$.
Using this conclusion, summing inequality \eqref{FuncBound} from $1$ to $\bar j$, and using the lower bound on $\eta_j$ in \eqref{Output1}, we conclude that
\begin{align*}
(\bar j-1) \bar \rho^2 < \sum_{j=1}^{\bar j-1} \|R_j\|^2 \leq \sum_{j=1}^{\bar j} \|R_j\|^2 \overset{\eqref{FuncBound}}{\leq} C_{\sigma}\sum_{j=1}^{\bar j}
	\frac{g(W_{j-1})-g(W_j)}{\eta_j}
\overset{\eqref{Output1}}{\leq} \frac{C_{\sigma}}{\left[\min\left\{\eta_0,1/(4\Lip{0}) \right\}\right]} \left[ g(W_0)-g(W_{\bar j})\right]
\end{align*}
which contradicts \eqref{iteration count}.   
Hence, ADAP-AIPP only generates iteration indices $j\geq 1$ satisfying 
\begin{equation}\label{jBound}
j\leq 1+\frac{C_{\sigma}\max\left\{\eta_0^{-1},\,4\Lip{0}\right\}}{\bar\rho^2}\left[g(W_0)-g(W_{j})\right].
\end{equation}
In particular, since $g(W_j)\geq g_*$ for every generated index $j$, we have
\[
j\leq1+\frac{C_{\sigma}\max\left\{\eta_0^{-1},\,4\Lip{0}\right\}}{\bar\rho^2}\left[g(W_0)-g_{*}\right].
\]
Thus ADAP-AIPP cannot generate infinitely many iteration indices, and hence it must terminate. Applying \eqref{jBound} with $j=\ell$, where $\ell$ is ADAP-AIPP's final iteration index, gives \eqref{FinalIterationIndex}.
    
The fact that ADAP-AIPP's output $(\bar W,\bar R)$ satisfies \eqref{eq:aipp_stationary_output} is immediate from ADAP-AIPP's termination criterion in step 3 and the fact that the inclusion in \eqref{Output1} holds for every generated iteration index $j\geq 1$. Moreover, it follows from \eqref{FuncBound} and a simple induction argument that $g(W_j)\leq g(W_0)$ for all iteration indices $j\geq 1$, which implies that output $(\bar W,\bar R)$
satisfies \eqref{MonoADAPAIPP}.

(b) The bound on the total number of ADAP-FISTA calls in \eqref{FISTACalls} is immediate from the bound on ADAP-AIPP's final iteration index $\ell$ in \eqref{FinalIterationIndex}, Lemma~\ref{Step1Lemma}(a), and the fact that \eqref{DiamDec} implies that $4\Lip{j}\leq 4\Lip{0}$.
\end{proof}

\begin{remark}
Each ADAP-FISTA call performs at most
\begin{equation}\label{eq:aipp_total_resolvent_bound}
{\mathcal{O}}\left(
\sqrt{\bar \Theta}\log_2^+(\bar \Theta)+\log_2^+(\bar \Theta)
\right).
\end{equation}
gradient evaluations and projections onto $\Omega$ where $\bar \Theta$ is as in \eqref{barThetaDef}.    
\end{remark}

\begin{proof}
The result follows immediately from the complexity bound in \eqref{eq:aipp_complexity} fact that Lemma~\ref{BarThetaLemma} implies that $\Theta_j=\mathcal O(\bar \Theta)$, and the first inequality in \eqref{MUniformBound}.
\end{proof}

\bibliographystyle{siam}
\bibliography{references}
\end{document}